\documentclass[12pt,leqno]{amsart}
\usepackage{amsmath,amssymb,amscd, amsthm}
\usepackage[ansinew]{inputenc}
\usepackage[all]{xy}
\usepackage{hyperref}
\usepackage{mathrsfs}
 \let\mathcal\mathscr

\def\CyrillicGuillemets{\DeclareFontEncoding{OT2}{}{}%
     \DeclareFontSubstitution{OT2}{wncyr}{m}{n}%
     \DeclareTextCommand{\guillemotleft}{OT1}{%
        {\fontencoding{OT2}\fontfamily{wncyr}\selectfont\char60}}%
     \DeclareTextCommand{\guillemotright}{OT1}{%
        {\fontencoding{OT2}\fontfamily{wncyr}\selectfont\char62}}}

\def\LasyGuillemets{%
     \DeclareTextCommand{\guillemotleft}{OT1}{\hbox{%
        \fontencoding{U}\fontfamily{lasy}\selectfont(\kern-0.20em(}}%
     \DeclareTextCommand{\guillemotright}{OT1}{\hbox{%
        \fontencoding{U}\fontfamily{lasy}\selectfont)\kern-0.20em)}}}
  \IfFileExists{ot2wncyr.fd}{\CyrillicGuillemets}{\LasyGuillemets}
   \DeclareTextSymbolDefault{\guillemotleft}{OT1}
   \DeclareTextSymbolDefault{\guillemotright}{OT1}
   \def\guill@spacing{\penalty\@M\hskip.8\fontdimen2\font
                               plus.3\fontdimen3\font
                               minus.8\fontdimen4\font}

\newtheorem{theorem}[subsubsection]{Theorem}
\newtheorem{proposition}[subsubsection]{Proposition}

\newtheorem{lemma}[subsubsection]{Lemma}

\newtheorem{corollary}[subsubsection]{Corollary}

\newtheorem{thm}[subsubsection]{Theorem}
\newtheorem{cor}[subsubsection]{Corollary}
\newtheorem{lem}[subsubsection]{Lemma}

\newtheorem{prop}[subsubsection]{Proposition}

\theoremstyle{definition}
\newtheorem{ex}[subsubsection]{Example}

\newtheorem{rem}[subsubsection]{Remark}
\newtheorem{example}[subsubsection]{Example}
\numberwithin{equation}{subsubsection}
\newtheorem{definition}[subsubsection]{Definition}
\newtheorem{remark}[subsubsection]{Remark}

\newtheorem{assume}[subsubsection]{Assumption}

  \DeclareMathSymbol{A}{\mathalpha}{operators}{`A}%
  \DeclareMathSymbol{B}{\mathalpha}{operators}{`B}%
  \DeclareMathSymbol{C}{\mathalpha}{operators}{`C}%
  \DeclareMathSymbol{D}{\mathalpha}{operators}{`D}%
  \DeclareMathSymbol{E}{\mathalpha}{operators}{`E}%
  \DeclareMathSymbol{F}{\mathalpha}{operators}{`F}%
  \DeclareMathSymbol{G}{\mathalpha}{operators}{`G}%
  \DeclareMathSymbol{H}{\mathalpha}{operators}{`H}%
  \DeclareMathSymbol{I}{\mathalpha}{operators}{`I}%
  \DeclareMathSymbol{J}{\mathalpha}{operators}{`J}%
  \DeclareMathSymbol{K}{\mathalpha}{operators}{`K}%
  \DeclareMathSymbol{L}{\mathalpha}{operators}{`L}%
  \DeclareMathSymbol{M}{\mathalpha}{operators}{`M}%
  \DeclareMathSymbol{N}{\mathalpha}{operators}{`N}%
  \DeclareMathSymbol{O}{\mathalpha}{operators}{`O}%
  \DeclareMathSymbol{P}{\mathalpha}{operators}{`P}%
  \DeclareMathSymbol{Q}{\mathalpha}{operators}{`Q}%
  \DeclareMathSymbol{R}{\mathalpha}{operators}{`R}%
  \DeclareMathSymbol{S}{\mathalpha}{operators}{`S}%
  \DeclareMathSymbol{T}{\mathalpha}{operators}{`T}%
  \DeclareMathSymbol{U}{\mathalpha}{operators}{`U}%
  \DeclareMathSymbol{V}{\mathalpha}{operators}{`V}%
  \DeclareMathSymbol{W}{\mathalpha}{operators}{`W}%
  \DeclareMathSymbol{X}{\mathalpha}{operators}{`X}%
  \DeclareMathSymbol{Y}{\mathalpha}{operators}{`Y}%
  \DeclareMathSymbol{Z}{\mathalpha}{operators}{`Z}%

\DeclareMathOperator{\Id}{Id}
\DeclareMathOperator{\Hom}{\mathrm{Hom}}
\DeclareMathOperator{\Ext}{\mathrm{Ext}}

\DeclareMathOperator{\Rhom}{\mathcal{R}\it{hom}}

\DeclareMathOperator{\Ker}{\mathrm{ker}}

\DeclareMathOperator{\RHom}{\mathrm{Rhom}}
\DeclareMathOperator{\ext}{\mathcal{E}\it{xt}}
\DeclareMathOperator{\D}{{D}}
\DeclareMathOperator{\homo}{\mathcal{H}om}
\DeclareMathOperator{\Hu}{\mc{H}}

\DeclareMathOperator{\Otimes}{\stackrel{\mathbf L}{\otimes}}
 
\DeclareMathOperator{\U}{\mathcal{U}}
\DeclareMathOperator{\ra}{\rightarrow}

 \newcommand{\NN}{\mathbf{N}}
\DeclareMathOperator{\Spec}{\mathrm{Spec}}

\newcommand{\og}{\guillemotleft}
\newcommand{\fg}{\guillemotright}

\DeclareMathOperator{\hocolim}{\mathrm{hocolim}}

\renewcommand{\O}{\mathcal{O}}
\newcommand{\FF}{\mathbf{F}}

\newcommand{\CC}{\mathbf{C}}

\newcommand{\X}{\mathcal{X}}

\newcommand{\F}{\mathcal{F}}
\newcommand{\G}{\mathcal{G}}
\newcommand{\Y}{\mathcal{Y}}
\newcommand{\T}{\mathcal{T}}

\newcommand{\isom}{\xrightarrow{\sim}}
\newcommand{\si}{\textup{ if }}
\renewcommand{\le}{\textup{lis-\'et}}
\newcommand{\LE}{{\textup{Lisse-Et}}}
\newcommand{\ET}{{\textup{\'Etale}}}

\newcommand{\et}{\textup{\'et}}
\newcommand{\opp}{\textup{opp}}

\renewcommand{\!}{\textup{!`}}
\newcommand{\mc}{\mathcal}
\newtheorem{pg}[subsubsection]{}
\newcommand{\s}{{ \ \ \ }}

\begin{document}

\title{The six operations for sheaves on Artin stacks I: Finite Coefficients}

\author{Yves Laszlo and Martin Olsson}
\address{\'Ecole Polytechnique CMLS UMR 7640 CNRS F-91128 Palaiseau Cedex France}
 \email{laszlo@math.polytechnique.fr}
 \address{University of Texas at Austin
Department of Mathematics 1 University Station C1200 Austin, TX
78712-0257, USA}\email{molsson@math.utexas.edu}

\begin{abstract}In this paper we develop a theory of Grothendieck's six operations
of lisse-\'etale constructible sheaves on Artin stacks locally of finite
type over an affine regular noetherian scheme of dimension $\leq 1$. We also generalize the classical base change theorems and Kunneth formula to stacks, and  prove new results about cohomological descent for unbounded complexes. \end{abstract}

\maketitle

\section{Introduction}

 We denote by $\Lambda$ a Gorenstein local ring of dimension $0$
and characteristic $l$. Let $S$ be an affine regular, noetherian
 scheme of dimension $\leq 1$ and assume $l$ is invertible
on $S$. We assume that all $S$-schemes of finite type $X$ satisfy
$\mathrm{cd}_{l}(X)<\infty$ (see \ref{cohdimension} for more
discussion of this).   For an algebraic stack $\mc X$ locally of
finite type over $S$ and $*\in \{+, -, b, \emptyset, [a,b]\}$ we
write $\D_c^*(\mc X)$ for the full subcategory of the derived
category $\D^*(\mc X)$ of complexes of $\Lambda $--modules on the
lisse-\'etale site of $\mc X$ with constructible cohomology
sheaves.

In this paper we develop a theory of Grothendieck's six operations
of lisse-\'etale constructible sheaves on Artin stacks locally of finite
type over $S$\footnote{In fact our method could apply
to other situations like analytic stacks or non separated analytic
varieties.}. In  forthcoming papers, we will also develop a theory
of adic sheaves and perverse sheaves for Artin stacks. In addition to being of basic
foundational interest, we hope that the development of these six
operations for stacks will have a number of applications. Already
the work done in this paper (and the forthcoming ones) provides the
necessary tools needed in several papers on the geometric
Langland's program (e.g. \cite{LN}, \cite{Lau}, \cite{FGV}).  We
hope that it will also shed further light on the Lefschetz trace
formula for stacks proven by Behrend (\cite{Ber03}), and also to
versions of such a formula for stacks not necessarily of finite
type. We should also remark that recent work of Toen should
provide another approach to defining the six operations for
stacks, and in fact should generalize to a theory for $n$--stacks.

Let us describe more precisely the contents of this papers. For a
morphism $f:\mc X\rightarrow \mc Y$ of such $S$--stacks we define
functors
\begin{equation*}
Rf_*:\D_c^{+}(\mc X)\rightarrow \D_c^{+}(\mc Y), \ \
Rf_!:\D_c^-(\mc X)\rightarrow \D_c^-(\mc Y),
\end{equation*}
\begin{equation*}
Lf^*:\D_c(\mc Y)\rightarrow \D_c(\mc X), \ \ Rf^!:\D_c(\mc Y)\rightarrow \D_c(\mc X),
\end{equation*}
\begin{equation*}
\Rhom:\D_c^{-}(\mc X)^{\text{op}}\times \D_c^+(\mc X)\rightarrow \D_c^+(\mc X),
\end{equation*}
and
\begin{equation*}
(-){\Otimes}(-):\D_c^-(\mc X)\times \D_c^-(\mc X)\rightarrow
\D_c^-(\mc X)
\end{equation*}
satisfying all the usual adjointness properties that one has in
the theory for schemes\footnote{We will often write $f^*,f^!$, $f_*$, $f_!$ for
$Lf^*,Rf^!$, $Rf_*$, $Rf_!$.}.

The main tool is to define $Rf_!,f^!$, even for unbounded
constructible complexes, by duality. One of the key points is
that, as observed by Laumon,  the dualizing complex is a local
object of the derived category and hence  has to exist for stacks
by glueing (see~\ref{1.2}). Notice that this formalism
applies to non-separated schemes, giving a theory of cohomology
with compact supports in this case. Previously, Laumon and
Moret-Bailly constructed the truncations of dualizing complexes
for Bernstein-Lunts stacks (see~\cite{Lau-Mor2000}). Our
constructions reduces to theirs in this case. Another approach
using a dual version of cohomological descent has been suggested
by Gabber but seems to be technically much more complicated.

\begin{remark}\label{cohdimension} The cohomological dimension hypothesis on schemes of finite type over $S$ is achieved for
instance if $S$ is the spectrum of a finite field or of a
separably closed field. In dimension $1$, it will be achieved for
instance for the spectrum of a complete discrete valuation field
with residue field either finite or separably closed, or if $S$ is
a smooth curve over $\CC,\FF_q$ (cf. \cite{SGA43}, exp. X
and~\cite{Ser94}). In these situations, $\mathrm{cd}_l(X)$ is
bounded by a function of the dimension $\dim(X)$. Notice that, as
pointed out by Illusie, recent results of Gabber enables one to
dramatically weaken the hypothesis on $S$.
Unfortunately no written version of these results seems to be
available  at this time.\end{remark}

\subsection{Conventions}

Recall that for any ring $\O$ of a topos, the category of
complexes of $\O$-modules has enough $K$-injective (or
homotopically injective). Recall  that a complex $I$ is
$K$-injective if for any  acyclic complex of
$\O$-modules $A$, the complex $\mc Hom(A,I)$ is acyclic
(see~\cite{Spa88}). For instance, any injective resolution in the
sense of Cartan-Eileberg of a bounded below complex is
$K$-injective.  This result is due, at least for sheaves on a
topological space to~\cite{Spa88} and enables him to extend the
formalism of direct images and $\Rhom$ to unbounded complexes. But
this result is true for any Grothendieck category (\cite{Serp03}).
Notice that the category of $\O$-modules has enough K-flat
objects, enabling one to define $\Otimes$ for unbounded objects
(\cite{Spa88}).

All the stacks we will consider will be locally of finite type over
$S$.  As in \cite{Lau-Mor2000}, lemme 12.1.2, the lisse-\'etale
topos $\X_\le$ can be defined using the site $\LE(\X)$ whose
objects are $S$-morphisms $u:U\ra\X$ where $U$ is an algebraic
space which is \emph{separated and of finite type over $S$}. The
topology is generated by the pretopology such that the covering
families are finite families $(U_i,u_i)\ra (U,u)$ such that
$\bigsqcup U_i\ra U$ is surjective and \'etale (use the comparison
theorem~\cite{SGA41}, III.4.1 remembering $\X$ is locally of
finite type over $S$). Notice that products over $\X$ are
representable in $\LE(\X)$, simply because the diagonal morphism
$\X\ra\X\times_S\X$ is representable by definition
(\cite{Lau-Mor2000}).

If $C$ is a complex of sheaves and $d$ a locally constant valued
function $C(d)$ is the Tate twist and $C[d]$ the shifted complex.
We denote $C(d)[2d]$ by $C{\langle}d{\rangle}$. Let
$\Omega=\Lambda{\langle}\dim(S){\rangle}$ be the dualizing complex
of $S$ (\cite{SGA4.5}, "Dualit\'e").

\section{Homological algebra}

\subsection{Existence of $K$--injectives}

 Let $(\mc  S, {\O} )$ denote a ringed site, and let
$\mc  C$ denote a full subcategory of the category of ${\O}
$--modules on $\mc  S$. Let $M$ be a complex of ${\O} $--modules
on $\mc  S$. By (\cite{Spa88}, 3.7) there exists a morphism of
complexes $f:M\rightarrow I$ with the following properties:
\begin{enumerate}
\item [(i)] $I = \varprojlim I_n$ where each $I_n$ is a bounded
below complex of flasque ${\O} $--modules. \item [(ii)] The
morphism $f$ is induced by a compatible collection of
quasi--isomorphisms $f_n:\tau _{\geq -n}M\rightarrow I_n.$ \item
[(iii)] For every $n$ the map $I_n\rightarrow I_{n-1}$ is
surjective with kernel $K_n$ a bounded below complex  of flasque
${\O} $--modules. \item [(iv)] For any pair of integers $n$ and
$i$ the sequence
\begin{equation}
0\rightarrow K_n^i\rightarrow I_n^i\rightarrow
I_{n-1}^i\rightarrow 0
\end{equation}
is split.
\end{enumerate}

\begin{rem} In fact (\cite{Spa88}, 3.7) shows that we can choose $I_n$ and $K_n$ to be complexes of injective ${\O} $--modules (in which case (iv) follows from (iii)).  However, for technical reasons it is sometimes useful to know that one can work just with flasque sheaves.
\end{rem}

 We make the following finiteness
assumption, which is the analog of~\cite{Spa88}, 3.12 (1).

\begin{assume}\label{assumption1}  For any  object $U\in \mc  S$ there
exists a covering $\{U_i\rightarrow U\}_{i\in I}$ and an integer
$n_0$ such that for any sheaf of ${\O} $--modules $F\in \mc  C$ we
have $H^n(U_i, F) = 0$ for all $n\geq n_0$.
\end{assume}

\begin{ex}
Let $\mc  S=\LE(\mc X)$ be the lisse-\'etale site of an algebraic
$S$-stack locally of finite type $\mc  X$ and ${\O}$ a constant
local Artinian ring
 of characteristic $\Lambda$ invertible on $S$. Then the class $\mc
C$ of all ${\O}$-sheaves, cartesian or not, satisfies the
assumption. Indeed, if $U\in\mc  S$ is of finite type over $S$ and
$F\in\mc  S$, one has
$H^n(U,F)=H^n(U_{\et},F_U)$\footnote{Cf.~\ref{LETdata} below}
which is zero for $n$ bigger than a constant depending only on $U$
(and not on $F$). Therefore, one can take the trivial covering in
this case. We could also take ${\O}=\mc O_{\mc  X}$ and $\mc C$ to
be the class of quasi-coherent sheaves.
\end{ex}
With hypothesis \ref{assumption1}, one has the following criterion for $f$ being
 a quasi-isomorphism (cf. \cite{Spa88}, 3.13).
\begin{prop}\label{1.4} Assume that $\mc  H^j(M)\in \mc  C$ for all $j$.  Then the  map $f$ is a quasi--isomorphism.  In particular, if each $I_n$ is a complex of injective ${\O} $--modules then
by {\rm \cite{Spa88}, 2.5,} $f:M\rightarrow I$ is a $K$--injective
resolution of $M$.
\end{prop}
\begin{proof}
For a fixed integer $j$, the map $\mc  H^j(M)\rightarrow \mc
H^j(I_n)$ is an isomorphism for $n$ sufficiently big.  Since this
isomorphism factors as
\begin{equation}
\mc  H^j(M)\rightarrow \mc  H^j(I)\rightarrow \mc  H^j(I_n)
\end{equation}
it follows that the map $\mc  H^j(M)\rightarrow \mc  H^j(I)$ is
injective.

To see that $\mc  H^j(M)\rightarrow \mc  H^j(I)$ is surjective,
let $U\in \mc  S$ be an object and $\gamma \in \Gamma (U, I^j)$ an
element with $d\gamma = 0$ defining a class in $\mc  H^j(I)(U)$.
Since $I = \varprojlim I_n$ the class $\gamma $ is given by a
compatible collection of sections $\gamma _n\in \Gamma (U, I_n^j)$
with $d\gamma _n = 0$.

Let $(\mc U = \{U_i\ra U\},n_0)$ be the data provided by \ref{assumption1}.
Let $N$ be an integer greater than $n_0-j$. For $m>N$ and $U_i\in
\mc  U$ the sequence
\begin{equation}
\Gamma (U_i, K_m^{j-1})\rightarrow  \Gamma (U_i,
K_m^{j})\rightarrow \Gamma (U_i, K_m^{j+1})\rightarrow \Gamma
(U_i, K_m^{j+2})
\end{equation}
is exact.  Indeed $K_m$ is a bounded below complex with $\mc
H^j(K_m)\in \mc  C$ for every $j$ and $\mc  H^j(K_m) = 0$ for
$j\geq -m+2$. It follows that $H^j(U_i, K_m) = 0$ for $j\geq
n_0-m+2$.

Since the maps $\Gamma (U_i, I_m^r)\rightarrow \Gamma (U_i,
I_{m-1}^r)$ are also surjective for all $m$ and $r$, it follows
from (\cite{Spa88}, 0.11) applied to the system
\begin{equation}
\Gamma (U_i, I_m^{j-1})\rightarrow \Gamma (U_i, I_m^i)\rightarrow \Gamma (U_i, I_m^{j+1})\rightarrow \Gamma (U_i, I_m^{j+2})
\end{equation}
that the map
\begin{equation}
H^j(\Gamma (U_i, I))\rightarrow H^j(\Gamma (U_i, I_m))
\end{equation}
is an isomorphism.

Then since the map $\mc  H^j(M)\rightarrow \mc  H^j(I_m)$ is an
isomorphism it follows that for every $i$ the restriction of
$\gamma $ to $U_i$ is in the image of $\mc  H^j(M)(U_i)$.

\hfill\hfill$\square$\end{proof}

 Next consider a fibred topos
$\mc T\ra D$  with corresponding total topos $\mc T_\bullet$
(\cite{SGA42}, VI.7). We call  $\mc T_\bullet$  a
\emph{$D$-simplicial topos}. Concretely, this means that for each
$i\in D$ the fiber $\mc T_i$ is a topos and that any
$\delta\in\Hom_D(i,j)$ comes together with a morphism of topos
$\delta:\mc T_i\ra\mc T_j$ such that $\delta^{-1}$ is the inverse
image functor of the fibred structure. The objects of the total
topos are simply collections $(F_i\in E_i)_{i\in D}$ together with
functorial transition morphisms $\delta^{-1}F_j\ra F_i$ for any
$\delta\in \Hom_D(i,j)$. We assume furthermore that $\mc
T_\bullet$ is ringed by a ${\O} _\bullet$ and that for any
$\delta\in\Hom_D(i,j)$, the morphism $\delta:(\mc T_i,{\O}_i)\ra
(\mc T_j,{\O}_j)$ is flat.

\begin{example} Let $\Delta ^+$ be  the category whose objects are the ordered sets $[n] = \{0, \dots , n\}$ ($n\in \mathbb{N}$) and whose morphisms are injective order-preserving maps.  Let $D$ be the opposite category of $\Delta ^+$.  In this case $\mc T_.$ is
called a \emph{strict simplicial topos}.  For instance, if
$U\ra\mc X$ is a presentation, the simplicial algebraic space
$U_\bullet=\mathrm{cosq}_0(U/\mc X)$ defines a strict simplicial
topos $U_{\bullet\le}$ whose fiber over $[n]$ is $U_{n\le}$. For a
morphism $\delta :[n]\rightarrow [m]$ in $\Delta ^+$ the  morphism
$\delta:\mc T_{m}\ra \mc T_n$ is induced by the (smooth)
projection $U_m\ra U_n$ defined by $\delta\in
\Hom_{\Delta^{+\opp}}([m],[n])$.
\end{example}
\begin{example}
Let $\mathbf N$ be the natural numbers viewed as a category in
which $\Hom(n, m)$ is empty unless $m\geq n$ in which case it
consists of a unique element. For a topos $T$ we can then define
an $\mathbf N$-simplicial topos $T^{\mathbf N}$.  The fiber over
$n$ of $T^\NN$ is $T$ and the transition morphisms by the identity
of $T$. The topos $T^{\mathbf N}$ is the category of projective
systems in $T$. If ${\O} _\bullet $ is a constant projective
system of rings then the flatness assumption is also satisfied, or
more generally if $\delta ^{-1}{\O} _n\rightarrow {\O} _m$ is an
isomorphism for any morphism $\delta :m\rightarrow n$ in
$\mathbb{N}$ then the flatness assumption holds.
\end{example}

 Let $\mc  C_\bullet $ be a full subcategory of the category of
${\O} _\bullet $--modules on a ringed $D$-simplicial topos $(\mc
T_\bullet,{\O}_\bullet)$. For $i\in D$, let $e_i:\mc T_n\ra\mc
T_\bullet$ the morphism of topos defined by
$e_i^{-1}F_\bullet=F_n$ (cf. \cite{SGA42}, Vbis, 1.2.11). Recall
that the family $e_i^{-1},i\in D$ is conservative. Let ${\mc C}_i$
denote the essential image of ${\mc C}_\bullet$ under $e_{i}^{-1}$
(which coincides with $e_i^*$ on $\text{Mod}(\mc
T_\bullet,{\O}_\bullet)$ because $e_i^{-1}{\O}_\bullet={\O}_i$).
\begin{assume}\label{assume2}
For every $i\in D$ the ringed topos $(\mc T_i, {\O} _i)$ is
isomorphic to the topos of a ringed site satisfying
\ref{assumption1} with respect to $\mc  C_i$.
\end{assume}

\begin{example} Let $\mc T_\bullet$ be the topos $(\mc X_{\le})^{\mathbf
N}$ of a $S$-stack locally of finite type. Then,  the full
subcategory $\mc C_{\bullet}$ of $\text{Mod}(\mc
T_\bullet,{\O}_\bullet)$ whose objects are families $F_i$ of
\emph{cartesian} modules satisfies the hypothesis.\end{example}

Let $M$ be a complex of ${\O} _\bullet $--modules on $\mc
T_\bullet $. Again by (\cite{Spa88}, 3.7) there exists a morphism
of complexes $f:M\rightarrow I$ with the following properties:
\begin{enumerate}
\item [(S i)] $I = \varprojlim I_n$ where each $I_n$ is a bounded
below complex of injective modules. \item [(S ii)] The morphism
$f$ is induced by a compatible collection of quasi--isomorphisms
$f_n:\tau _{\geq -n}M\rightarrow I_n.$ \item [(S iii)] For every
$n$ the map $I_n\rightarrow I_{n-1}$ is surjective with kernel
$K_n$ a bounded below complex  of injective ${\O} $--modules.
\item [(S iv)] For any pair of integers $n$ and $i$ the sequence
\begin{equation}\label{extseq}
0\rightarrow K_n^i\rightarrow I_n^i\rightarrow
I_{n-1}^i\rightarrow 0
\end{equation}
is split.
\end{enumerate}

\begin{prop}\label{1.7} Assume that $\mc  H^j(M)\in \mc  C_\bullet $ for all
$j$.  Then the morphism $f$ is a quasi--isomorphism and
$f:M\rightarrow I$ is a $K$--injective resolution of $M$.
\end{prop}
\begin{proof}
By \cite{Spa88}, 2.5, it suffices to show that $f$ is a
quasi--isomorphism.  For this in turn it suffices to show that for
every $i\in D $ the restriction $e_i^*f:e_i^*M\rightarrow e_i^*I$
is a quasi--isomorphism of complexes of ${\O}_i$-modules since the
family $e_i^* = e_i^{-1}$ is conservative. But
$e_i^*:\text{Mod}(\mc T_\bullet,{\O}_\bullet)\ra \text{Mod}(\mc
T_i,{\O}_i)$ has a left adjoint $e_{i!}$ defined by
$$[e_{i!}(F)]_j=\oplus_{\delta\in\Hom_D(j,i)}\delta^*F$$
with the obvious transition  morphisms. It  is exact by the
flatness of the morphisms $\delta$. It follows that $e_i^*$ takes
injectives to injectives and commutes with direct limits. We can
therefore apply \ref{1.4} to  $e_i^*M\ra e_i^*I$ to deduce that
this map is a quasi--isomorphism. \hfill\hfill$\square$\end{proof}

In what follows we call a $K$--injective resolution $f:M\rightarrow I$ obtained from data (i)-(iv) as above a \emph{Spaltenstein resolution}.

The main technical lemma is the following.

\begin{lem}\label{key-truncation} Let $\epsilon :(\mc T_\bullet , {\O} _\bullet )\rightarrow (S, \Psi )$ be a morphism of ringed topos, and let $C$ be a complex of ${\O} _\bullet $--modules.
Assume that \begin{enumerate}
   \item $\mc  H^n(C)\in \mc  C_\bullet $ for
all $n$.
   \item There exists $i_0$ such that $R^i\epsilon_*\mc  H^n(C)=0$ for any $n$ and any $i> i_0$.
\end{enumerate}
Then, if $j\geq -n+i_0$, we have
$R^j\epsilon_*C=R^j\epsilon_*\tau_{\geq -n}C$.
\end{lem}
\begin{proof}
By \ref{1.7} and assumption (1), there exists a  Spaltenstein
resolution $f:C\rightarrow I$ of $C$.
Let $J_n := \epsilon _*I_n$ and $D_n :=\epsilon _*K_n$. Since the
sequences \ref{extseq} are split, the sequences
\begin{equation}
0\rightarrow D_n\rightarrow J_n\rightarrow J_{n-1}\rightarrow 0
\end{equation}
are exact.

 The exact sequence \ref{extseq} and property
(S ii) defines a distinguished triangle
$$K_n\rightarrow\tau_{\geq -n}C\rightarrow\tau_{\geq -n+1}C$$ showing that  $K_n$ is
quasi--isomorphic to $\mc  H^{-n}(C)[n]$. Because $K_n$ is a
bounded below complex of injectives, one gets $$R\epsilon_*\mc
H^{-n}(C)[n]=\epsilon_*K_n$$ and accordingly
$$R^{j+n}\epsilon_*\mc
H^{-n}(C)=\mc  H^j(\epsilon_*K_n)=\mc  H^j(D_n).$$

By assumption (2), we have therefore $$\mc  H^j(D_n) = 0\textup{
for }j> -n+i_0.$$ By (\cite{Spa88}, 0.11) this implies that
$$\mc  H^j(\varprojlim J_n) \rightarrow \mc  H^j(J_{n})$$ is an isomorphism for $j\geq
-n+i_0$. But, by adjunction, $\epsilon_*$ commutes with projective
limit. In particular, one has $$\varprojlim J_n=\epsilon_*I,$$
and by (S i) and (S ii)
$$R\epsilon_*C=\epsilon_*I\textup{ and }R\epsilon_*\tau_{\geq
-n}C=\epsilon_*J_n.$$ Thus for any $n$ such that $j\geq -n+i_0$
one has
\begin{equation}
R^j\epsilon _*C  = \mc  H^j(\epsilon_*I) = \mc  H^j(J_{n}) =
R^j\epsilon _*\tau_{\geq -n}C.
\end{equation}

\hfill\hfill$\square$\end{proof}

\subsection{The descent theorem}

 Let $(\mc T_\bullet ,{\O}_{\bullet})$ be a
simplicial or strictly simplicial \footnote{One could replace
simplicial by multisimplicial} ringed topos ($D=\Delta^\opp$ or
$D=\Delta^{+\opp}$), let $(S, \Psi)$ be another ringed topos, and
let $\epsilon :(\mc T_\bullet ,{\O}_{\bullet})\rightarrow (S, \Psi
)$ be an augmentation. Assume that $\epsilon $ is  a flat morphism
(i.e. for every $i\in D $, the morphism of ringed topos $(\mc T_i,
{\O}_i)\rightarrow (S, \Psi )$ is a flat morphism).

Let $\mc C$ be a full subcategory of the category of $\Psi
$--modules, and assume that $\mc C$ is closed under kernels,
cokernels and extensions (one says that $\mc C$ is a \emph{Serre
subcategory}). Let $\D(S)$ denote the derived category of $\Psi
$--modules, and let $\D_{\mc C}(S)\subset \D(S)$ be the full
subcategory consisting of complexes whose cohomology sheaves are
in $\mc C$. Let $\mc C_\bullet $ denote the essential image of
$\mc C$ under the functor $\epsilon ^*:\text{\rm Mod}(\Psi
)\rightarrow \text{\rm Mod}({\O} _\bullet )$.

We assume the following condition holds:

\begin{assume}\label{2.13}
Assumption \ref{assume2} holds (with respect to $\mc  C_\bullet
$), and $\epsilon^*:\mc C\ra\mc C_\bullet$ is an equivalence of
categories with quasi-inverse $R\epsilon_*$.
\end{assume}

\begin{lem}   The full subcategory $\mc C_{\bullet }\subset \text{\rm Mod}({\O} _\bullet )$ is closed
under extensions, kernels and cokernels.
\end{lem}
\begin{proof}
Consider an extension of sheaves of ${\O} _\bullet $--modules
\begin{equation}
\begin{CD}
0@>>> \epsilon ^*F_1@>>> E@>>> \epsilon ^*F_2@>>> 0,
\end{CD}
\end{equation}
where $F_1, F_2\in \mc C$.  Since $R^1\epsilon _*\epsilon ^*F_1 =
0$ and the maps $F_i\rightarrow R^0\epsilon _*\epsilon ^*F_i$ are
isomorphisms, we obtain by applying $\epsilon ^* \epsilon _*$ a
commutative diagram with exact rows
\begin{equation}
\begin{CD}
0@>>> \epsilon ^*F_1@>>> \epsilon ^*\epsilon _*E@>>> \epsilon
^*F_2@>>> 0\\
@. @V\text{id}VV @V\alpha VV @VV\text{id}V @. \\
0@>>> \epsilon ^*F_1@>>> E@>>> \epsilon ^*F_2@>>> 0.
\end{CD}
\end{equation}
It follows that $\alpha $ is an isomorphism.  Furthermore, since
$\mc C$ is closed under extensions we have $\epsilon _*E\in \mc
C$. Let $f\in\Hom(\epsilon^*F_1,\epsilon^*F_2)$. There exists a
unique $\varphi\in\Hom(F_1,F_2)$ such that $f=\epsilon^*\varphi$.
Because $\epsilon^*$ is exact, it maps the kernel and cokernel of
$\varphi$, which are objects of $\mc C$, to the kernel and
cokernel of $f$ respectively. Therefore, the latter are objects of
$\mc C_\bullet$. \hfill\hfill$\square$\end{proof}

Let $\D(\mc T_\bullet )$ denote the derived category of ${\O}
_\bullet $--modules, and let $\D_{\mc C_\bullet }(\mc T_\bullet
)\subset \D(\mc T_\bullet )$ denote the full subcategory of
complexes whose cohomology sheaves are in $\mc C_\bullet $.

Since $\epsilon $ is a flat morphism, we obtain a morphism of
triangulated categories (the fact that these categories are
triangulated comes precisely from the fact that both $\mc C$ and
$\mc C_\bullet$ are Serre categories \cite{Grivel}).
\begin{equation}\label{pullback}
\epsilon ^*:\D_{\mc C}(S)\rightarrow \D_{\mc C_\bullet }(\mc T_\bullet
).
\end{equation}

\begin{thm}\label{mainthm} The functor $\epsilon^*$ of \ref{pullback} is an equivalence
of triangulated categories with quasi--inverse given by $R\epsilon
_*$.
\end{thm}

\begin{proof}
Note first that if $M_\bullet \in \D_{\mc C_\bullet }(\mc T_\bullet )$,
then by lemma~\ref{key-truncation}, for any integer $j$ there
exists $n_0$ such that $R^j\epsilon _*M_\bullet = R^j\epsilon _*\tau
_{\geq n_0}M_\bullet $. In particular, we get by induction
$R^j\epsilon _*M_\bullet \in \mc C$. Thus $R\epsilon _*$ defines a
functor
\begin{equation}
R\epsilon _*:\D_{\mc C_\bullet }(\mc T_\bullet )\rightarrow \D_{\mc
C}(S).
\end{equation}
To prove \ref{mainthm} it suffices to show that for $M_\bullet \in
\D_{\mc C_\bullet }(\mc T_\bullet )$ and $F\in \D_{\mc C}(S)$ the
adjunction maps
\begin{equation}
\epsilon ^*R\epsilon _*M_\bullet \rightarrow M_\bullet , \ \ \
F\rightarrow R\epsilon _*\epsilon ^*F.
\end{equation}
are isomorphisms. For this note that for any integers $j$ and $n$
there are commutative diagrams
\begin{equation}
\begin{CD}
\epsilon ^*R^j\epsilon _*M_\bullet @>>> \mc  H^j(M_\bullet )\\
@VVV @VVV \\
\epsilon ^*R^j\epsilon _*\tau _{\geq n}M_\bullet @>>> \mc  H^j(\tau
_{\geq n}M_\bullet ),
\end{CD}
\end{equation}
and
\begin{equation}
\begin{CD}
\mc  H^j(F)@>>> R^j\epsilon _*\epsilon ^*F\\
@VVV @VVV \\
\mc  H^j(\tau _{\geq n}F)@>>> R^j\epsilon _*\epsilon ^*\tau _{\geq
n}F.
\end{CD}
\end{equation}
By the observation at the begining of the proof, there exists an
integer $n$ so that the vertical arrows in the above diagrams are
isomorphisms. This reduces the proof \ref{mainthm} to the case of
a bounded below complex. In this case one reduces by devissage to
the case when $M_\bullet \in \mc C_\bullet $ and $F\in \mc C$ in
which case the result holds by assumption.
\hfill\hfill$\square$\end{proof}

The Theorem applies in particular to the following examples.

\begin{ex}
Let $S$ be an algebraic space and $X_\bullet \rightarrow S$ a flat hypercover by algebraic spaces.   We then obtain an augmented simplicial topos $\epsilon :(X_{\bullet , \et}, \mc O_{X_{\bullet, \et}})\rightarrow (S_\et, \mc O_{\et}).$ Note that this augmentation is flat.  Let $\mc C$ denote the category of quasi--coherent sheaves on $S_\et $.  Then the category $\mc C_\bullet $ is the category of cartesian sheaves of $\mc O_{X_{\bullet, \et }}$--modules whose restriction to each $X_n$ is quasi--coherent.  Let $\D_{\text{qcoh}}(X_\bullet )$ denote the full subcategory of the derived category of $\mc O_{X_\bullet, \et }$--modules whose cohomology sheaves are quasi--coherent, and let $\D_{\text{qcoh}}(S)$ denote the full subcategory of the derived category of $\mc O_{S_\et }$--modules whose cohomology sheaves are quasi--coherent.  Theorem \ref{mainthm} then shows that the pullback functor
\begin{equation}
\epsilon ^*:\D_{\text{qcoh}}(S)\rightarrow \D_{\text{qcoh}}(X_\bullet )
\end{equation}
is an equivalence of triangulated categories with quasi--inverse $R\epsilon _*$.
\end{ex}
\begin{ex} Let $\mc X$ be an algebraic stack and let $U_\bullet \rightarrow \mc X$ be a smooth
hypercover  by algebraic spaces.  Let $\D(\mc X)$ denote the
derived category of sheaves of $\mc O_{\mc X_{{\le}}}$--modules in
the topos $\mc X_{{\le}}$, and let $\D_{\text{qcoh}}(\mc X)\subset
\D(\mc X)$ be the full subcategory of complexes with
quasi--coherent cohomology sheaves.

Let $U_{\bullet }^+$ denote the strictly simplicial algebraic space
obtained from $U_\bullet $ by forgetting the degeneracies.  Since
the Lisse-\'Etale topos is functorial with respect to smooth
morphisms, we therefore obtain a strictly simplicial topos
$U_{\bullet {\le}}$ and a flat morphism of ringed topos $$\epsilon
:(U_{\bullet {\le}}, \mc O_{U_{\bullet {\le}}})\rightarrow (\mc
X_{{\le}}, \mc O_{\mc X_{{\le}}}).$$  Then \ref{2.13} holds with
$\mc C$ equal to the category of quasi--coherent sheaves on $\mc
X$.  The category $\mc C_\bullet $ in this case is the category of
cartesian $\mc O_{U_{\bullet {\le}}}$--modules $M_\bullet $ such that
the restriction $M_n$ is a quasi--coherent sheaf on $U_n$ for all
$n$.  By \ref{mainthm} we then obtain an equivalence of
triangulated categories
\begin{equation}
\D_{\text{qcoh}}(\mc X)\rightarrow \D_{\text{qcoh}}(U_{\bullet ,
{\le}}),
\end{equation}
where the right side denotes the full subcategory of the derived
category of $\mc O_{U_{\bullet {\le}}}$--modules with cohomology
sheaves in $\mc C_{\bullet }$.

On the other hand, there is also a natural morphism of ringed
topos $$\pi :(U_{\bullet {\le}}, \mc O_{U_{\bullet {\le}}})\rightarrow
(U_{\bullet \et}, \mc O_{U_{\bullet \et }})$$ with $\pi _*$ and $\pi
^*$ both exact functors.  Let $\D_{\text{qcoh}}(U_{\bullet \et})$
denote the full subcategory of the derived category of $\mc
O_{U_{\bullet \et }}$--modules consisting of complexes whose
cohomology sheaves are quasi--coherent (i.e. cartesian and
restrict to a quasi--coherent sheaf on each $U_{n \et}$).  Then
$\pi $ induces an equivalence of triangulated categories
$\D_{\text{qcoh}}(U_{\bullet \et})\simeq \D_{\text{qcoh}}(U_{\bullet
{\le}})$.  Putting it all together we obtain an equivalence of
triangulated categories $\D_{\text{qcoh}}(\mc X_{{\le}})\simeq
\D_{\text{qcoh}}(U_{\bullet \et})$.
\end{ex}

\begin{ex}\label{etexample} Let  $\mc  X$ be an
algebraic stack locally of finite type over $S$ and ${\O} $ be a
constant local Artinian ring $\Lambda$ of characteristic
invertible on $S$. Let $U_\bullet \rightarrow \mc  X$ be a smooth
hypercover by algebraic spaces, and $\mc T_\bullet $ the localized
topos $\mc X_{\text{{\le}}}|_{U_\bullet }$.  Take $\mc  C$ to be
the category of constructible sheaves of ${\O} $--modules.  Then
\ref{mainthm} gives an equivalence $\D_c(\mc
X_{\text{{\le}}})\simeq \D_c(\mc T_\bullet , \Lambda )$.  On the
other hand, there is a natural morphism of topos $\lambda :\mc
T_\bullet \rightarrow U_{\bullet , et}$ and one sees immediately
that $\lambda _*$ and $\lambda ^*$ induce an equivalence of
derived categories $\D_c(\mc T_\bullet , \Lambda )\simeq
\D_c(U_{\bullet , et}, \Lambda )$.  It follows that $\D_c(\mc
X_{{\le}})\simeq \D_c(U_{\bullet , et})$.
\end{ex}

\subsection{The BBD gluing lemma}

The purpose of this section  is to explain how to modify the proof of  the gluing lemma \cite[3.2.4]{BBD82} for unbounded complexes.

Let $\widetilde {\Delta } $ denote the strictly simplicial category of finite ordered sets with injective order preserving maps, and let $\Delta ^+\subset \widetilde {\Delta } $ denote the full subcategory of nonempty finite ordered sets. For a morphism $\alpha $ in $\widetilde {\Delta } $ we write $s(\alpha )$ (resp. $b(\alpha )$) for its source (resp. target).

Let $T$ be a topos and $U_\cdot \rightarrow e$ a strictly simplicial hypercovering of the initial object $e\in T$. For
$[n]\in \widetilde {\Delta } $ write $U_n$ for the localized topos $T|_{U_n}$ where by definition we set $U_\emptyset = T$.  Then we obtain a strictly simplicial topos $U_\cdot $ with an augmentation $\pi :U_\cdot \rightarrow T$.

Let $\Lambda $ be a sheaf of rings in $T$ and write also $\Lambda $ for the induced sheaf of rings in $U_\cdot $ so that $\pi $ is a morphism of ringed topos.

Let $\mathcal C_\cdot $ denote a full substack of the fibered and cofibered category over $\widetilde {\Delta } $
$$
[n]\mapsto (\text{category of sheaves of $\Lambda $--modules in $U_n$})
$$
such that each $\mathcal C_n$ is a Serre subcategory of the category of $\Lambda $--modules in $U_n$.
For any $[n]$ we can then form the derived category $D_{\mathcal C}(U_n, \Lambda )$ of complexes of $\Lambda $--modules whose cohomology sheaves are in $\mathcal C_n$.  The categories $D_{\mathcal C}(U_n, \Lambda )$ form a fibered and cofibered category over $\widetilde {\Delta } $.

We make the following assumptions on $\mathcal C$:

\begin{assume}\label{assumption1b} (i)  For any $[n]$ the topos $U_n$ is equivalent to the topos associated to a site $\mathcal S_n$ such that for any object $V\in \mathcal S_n$ there exists an integer $n_0$ and a covering $\{V_j\rightarrow V\}$ in $\mathcal S_n$ such that for any $F\in \mathcal C_n$ we have $H^n(V_j, F) = 0$ for all $n\geq n_0$.

(ii) The natural functor
$$
\mathcal C_{\emptyset }\rightarrow (\text{cartesian sections of $\mathcal C|_{\Delta ^+}$ over $\Delta ^+$})
$$
is an equivalence of categories.

(iii) The category $D(T, \Lambda )$ is compactly generated.
\end{assume}

\begin{rem} The case we have in mind is when $T$ is the lisse-\'etale topos of an algebraic stack $\mathcal X$ locally of finite type over an affine regular, noetherian
scheme of dimension $\leq 1$,  $U_\cdot $ is given by a hypercovering of $\mathcal X$ by schemes, $\Lambda $ is  a Gorenstein local ring of dimension $0$
and characteristic $l$ invertible on $\mathcal X$, and $\mathcal C$ is the category of constructible $\Lambda $--modules.  In this case the category $D_{c}(\mathcal X_{\text{lis-et}}, \Lambda )$ is compactly generated.  Indeed a set of generators is given by sheaves $j_!\Lambda [i]$ for $i\in \mathbb{Z}$ and $j:U\rightarrow \mathcal X$ an object of the lisse-\'etale site of $\mathcal X$.
\end{rem}

There is also a natural functor
\begin{equation}\label{Apullback}
D_{\mathcal C}(T, \Lambda )\rightarrow (\text{cartesian sections of $[n]\mapsto D_{\mathcal C}(U_n, \Lambda )$ over $\Delta ^+$}).
\end{equation}

\begin{thm}\label{1.2} Let $[n]\mapsto K_n\in D_{\mathcal C}(U_n, \Lambda )$ be a cartesian section of $[n]\mapsto D_{\mathcal C}(U_n, \Lambda )$ over $\Delta ^+$ such that $\mathcal Ext^i(K_0, K_0) = 0$ for all $i<0$.  Then $(K_n)$ is induced by a unique object $K\in D_{\mathcal C}(T, \Lambda )$ via the functor \ref{Apullback}.
\end{thm}

The uniqueness is the easy part:

\begin{lem}\label{D1} Let $K, L\in D(T, \Lambda )$ and assume that $\mathcal Ext^i(K, L) = 0$ for $i<0$.  Then $U\mapsto \text{\rm Hom}_{D(U, \Lambda )}(K|_U, L|_U)$ is a sheaf.
\end{lem}
\begin{proof}
Let $\mathcal H$ denote the complex $\mathcal Rhom(K, L)$.  By assumption the natural map $\mathcal H\rightarrow \tau _{\geq 0}\mathcal H$ is an isomorphism.  It follows that $\text{Hom}_{D(U, \Lambda )}(K|_U, L|_U)$ is equal to the value of $\mathcal H^0(\mathcal H)$ on $U$ which implies the lemma.
\end{proof}

The existence part is more delicate.  Let $\mathcal A$ denote the fibered and cofibered category over $\widetilde {\Delta } $ whose fiber over $[n]\in \widetilde {\Delta } $ is the category of $\Lambda $--modules in $U_n$. For a morphism $\alpha :[n]\rightarrow [m]$, $F\in \mathcal A(n)$ and $G\in \mathcal A(m)$ we have
$$
\text{Hom}_\alpha (F, G) = \text{Hom}_{\mathcal A(m)}(\alpha ^*F, G) = \text{Hom}_{\mathcal A(n)}(F, \alpha _*G).
$$
We write $\mathcal A^+$ for the restriction of $\mathcal A$ to $\Delta ^+$.

Define a new category $\text{tot}(\mathcal A^+)$ as follows:
\begin{enumerate}
\item [$\bullet $] The objects of $\text{tot}(\mathcal A^+)$ are collections of objects $(A^n)_{n\geq 0}$ with $A^n\in \mathcal A(n)$.
\item [$\bullet $] For two objects $(A^n)$ and $(B^n)$ we define
$$
\text{Hom}_{\text{tot}(\mathcal A^+)}((A^n), (B^n)):= \prod _{\alpha }\text{Hom}_\alpha (A^{s(\alpha )}, B^{b(\alpha )}),
$$
where the product is taken over all morphisms in $\Delta ^+$.
\item [$\bullet $] If $f=(f_\alpha )\in \text{Hom}((A^n), (B^n))$ and $g=(g_\alpha )\in \text{Hom}((B^n), (C^n))$ are two morphisms then the composite is defined to be the collection of morphisms whose $\alpha $ component is defined to be
$$
(g\circ f)_\alpha := \sum _{\alpha = \beta \gamma }g_\beta f_\gamma
$$
where the sum is taken over all factorizations of $\alpha $.
\end{enumerate}

The category $\text{tot}(\mathcal A^+)$ is an additive category.




Let $(K, d)$ be a complex in $\text{tot}(\mathcal A^+)$ so for every degree $n$ we are given a family of objects $(K^n)^m\in \mathcal A(m)$.  Set
$$
K^{n, m}:= (K^{n+m})^n.
$$
For $\alpha :[n]\rightarrow [m]$ in $\Delta ^+$ let $d(\alpha )$ denote the $\alpha $--component of $d$ so
$$
d(\alpha )\in \text{Hom}_\alpha ((K^p)^n, (K^{p+1})^m) = \text{Hom}_\alpha (K^{n, p-n}, K^{m, m-p-1})
$$
or equivalently  $d(\alpha )$ is a map $K^{n, p}\rightarrow K^{m, p+n-m+1}$. In particular, $d(\text{id}_{[n]})$ defines a  map $K^{n, m}\rightarrow K^{n, m+1}$ and as explained in \cite[3.2.8]{BBD82} this map makes $K^{n, *}$ a complex.  Furthermore for any $\alpha $ the map $d(\alpha )$ defines  an $\alpha $--map of complexes $K^{n, *}\rightarrow K^{m, *}$ of degree $n-m+1$.
The collection of complexes $K^{n, *}$ can also be defined as follows.  For an integer $p$ let $L^pK$ denote the subcomplex with $(L^pK)^{n, m}$ equal to $0$ if $n<p$ and $K^{n, m}$ otherwise.  Note that for any $\alpha :[n]\rightarrow [m]$ which is not the identity map $[n]\rightarrow [n]$ the image of $d(\alpha )$ is contained in $L^{p+1}K$.  Taking the associated graded of $L$ we see that
$$
\text{gr}_L^nK[n] = (K^{n, *}, d^{\prime \prime })
$$
where $d^{\prime \prime }$ denote the differential $(-1)^nd(\text{id}_{[n]})$.
Note that the functor $(K, d)\mapsto K^{n, *}$ commutes with the formation of cones and with shifting of degrees.

As explained in \cite[3.2.8]{BBD82} a complex in $\text{tot}(\mathcal A^+)$ is completely characterized by the data of a complex $K^{n, *}\in C(\mathcal A^+)$ for every $[n]\in \Delta ^+$ and for every morphism $\alpha :[n]\rightarrow [m] $ an $\alpha $--morphism $d(\alpha ):K^{n, *}\rightarrow K^{m, *}$ of degree $n-m+1$, such that $d(\text{id}_{[n]})$ is equal to $(-1)^n$ times the differential of $K^{n, *}$ and such that for every $\alpha $ we have
$$
\sum _{\alpha = \beta \gamma }d(\beta )d(\gamma ) = 0.
$$
Via this dictionary, a morphism $f:K\rightarrow L$ in $C(\text{tot}(\mathcal A^+))$ is given by an $\alpha $--map $f(\alpha ):K^{n, *}\rightarrow K^{m, *}$ of degree $n-m$ for every morphism $\alpha :[n]\rightarrow [m]$ in $\Delta ^+$ such that for any morphism $\alpha $ we have
$$
\sum _{\alpha = \beta \gamma }d(\beta )f(\gamma ) = \sum _{\alpha = \beta \gamma }f(\beta )d(\gamma ).
$$


Let $K(\text{tot}(\mathcal A^+))$ denote the category whose objects are complexes in $\text{tot}(\mathcal A^+)$ and whose morphisms are homotopy classes of morphisms of complexes.  The category $K(\text{tot}(\mathcal A^+))$ is a triangulated category.  Let $L\subset K(\text{tot}(\mathcal A^+))$ denote the full subcategory of objects $K$ for which each $K^{n, *}$ is acyclic for all $n$.   The category $L$ is a localizing subcategory of $K(\text{tot}(\mathcal A^+))$ in the sense of \cite[1.3]{Bok-Nee93} and hence the localized category $D(\text{tot}(\mathcal A^+))$ exists. The category $D(\text{tot}(\mathcal A^+))$ is obtained from $K(\text{tot}(\mathcal A^+))$ by inverting quasi--isomorphisms. Recall that an object $K\in K(\text{tot}(\mathcal A^+))$ is called \emph{$L$--local} if for any object $X\in L$ we have $\text{Hom}_{K(\text{tot}(\mathcal A^+))}(X, K) = 0$.
Note that the functor $K\mapsto K^{n, *}$ descends to a functor
$$
D(\text{tot}(\mathcal A^+))\rightarrow D(U_n, \Lambda ).
$$
We define  $D^+(\text{tot}(\mathcal A^+))\subset D(\text{tot}(\mathcal A^+))$ to be the full subcategory of objects $K$ for which there exists an integer $N$ such that $\mathcal H^j(K^{n, *}) = 0$ for all $n$ and all $j\leq N$.

Recall \cite[4.3]{Bok-Nee93} that a \emph{localization} for an object $K\in K(\text{tot}(\mathcal A^+))$ is a morphism $K\rightarrow I$ with $I$ an $L$--local object such that for any $L$--local object $Z$ the natural map
\begin{equation}\label{des1}
\text{Hom}_{K(\text{tot}(\mathcal A^+))}(I, Z)\rightarrow \text{Hom}_{K(\text{tot}(\mathcal A^+))}(K, Z)
\end{equation}
is an isomorphism.


\begin{lem} A morphism $K\rightarrow I$ is a localization if and only if $I$ is $L$--local and for every $n$ the map $K^{n, *}\rightarrow I^{n, *}$ is a quasi--isomorphism.
\end{lem}
\begin{proof}
By \cite[2.9]{Bok-Nee93} the morphism \ref{des1} can be identified with the natural map
\begin{equation}\label{des2}
\text{Hom}_{D(\text{tot}(\mathcal A^+))}(I, Z)\rightarrow \text{Hom}_{D(\text{tot}(\mathcal A^+))}(K, Z).
\end{equation}
If $K\rightarrow I$ is a localization it follows that this map is a bijection for every $L$--local $Z$.  By Yoneda's lemma applied to the full subcategory of $D(\text{tot}(\mathcal A^+))$ of objects which can be represented by $L$--local objects, it follows that this holds if and only if $K\rightarrow I$ induces an isomorphism in $D(\text{tot}(\mathcal A^+))$ which is the assertion of the lemma.
\end{proof}


\begin{prop}\label{L-local} Let $K\in C(\text{\rm tot}(\mathcal A^+))$ be an object with each $K^{n, *}$ homotopically injective.  Then $K$ is $L$--local.
\end{prop}
\begin{proof}
Let $X\in L$ be an object.  We have to show that any morphism $f:X\rightarrow K$ in $C(\text{tot}(\mathcal A^+))$ is homotopic to zero. Such a homotopy $h$ is given by a collection of maps $h(\alpha )$ such that
$$
f(\alpha ) = - \sum _{\alpha = \beta \gamma }d(\beta )h(\gamma )+h(\beta )d(\gamma ).
$$
We usually write just $h$ for $h(\text{id}_{[n]})$.

We construct these maps $h(\alpha )$ by induction on $b(\alpha )-s(\alpha )$. For $s(\alpha ) = b(\alpha )$ we choose the $h(\alpha )$ to be any homotopies between the maps $f(\text{id}_{[n]})$ and the zero maps.

For the inductive step, it suffices to show that
$$
\Psi (\alpha ) = f(\alpha )+d(\alpha )h+hd(\alpha ) + \sum ^{}_{\alpha = \beta \gamma }{}^\prime d(\beta )h(\gamma )+h(\beta )d(\gamma )
$$
commutes with the differentials $d$, where $\Sigma _{\alpha = \beta \gamma }'$ denotes the sum over all possible factorizations with $\beta $ and $\gamma $ not equal to the identity maps.  For then $\Psi (\alpha )$ is homotopic to zero and we can take $ h(\alpha )$ to be a homotopy between $\Psi (\alpha )$ and $0$.

Define
$$
A(\alpha ) = \sum _{\alpha = \beta \gamma }{}^\prime d(\beta )h(\gamma )+h(\gamma )d(\beta )
$$
and
$$
B(\alpha ) = d(\alpha )h+hd(\alpha )+A(\alpha ).
$$

\begin{lem}
$$
\sum _{\alpha = \beta \gamma }{}^\prime A(\beta )d(\gamma )-d(\beta )A(\gamma ) = \sum _{\alpha = \beta \gamma }{}^\prime h(\beta )S(\gamma )-S(\beta )h(\gamma ),
$$
where $S(\alpha )$ denotes $\sum _{\alpha = \beta \gamma }{}^\prime d(\beta )d(\gamma ).$
\end{lem}
\begin{proof}
\begin{eqnarray*}
\sum _{\alpha = \beta \gamma }{}^\prime A(\beta )d(\gamma )-d(\beta )A(\gamma ) & = & \sum _{\alpha = \epsilon \rho \gamma }{}^\prime d(\epsilon )h(\rho )d(\gamma )+h(\epsilon )d(\rho )d(\gamma ) - d(\epsilon )h(\rho )d(\gamma )-d(\epsilon )d(\rho )h(\gamma )\\
& = & \sum _{\alpha = \beta \gamma }{}^\prime h(\beta )S(\gamma )-S(\beta )h(\gamma ),
\end{eqnarray*}
where $\Sigma _{\alpha = \epsilon \rho \gamma }'$ denotes the sum over all possible factorizations with $\epsilon $, $\rho $,  and $\gamma $ not equal to the identity maps.
\end{proof}

\begin{lem}\label{Alem1}
$$\sum _{\alpha = \beta \gamma }{}^\prime B(\beta )d(\gamma )-d(\beta )B(\gamma )= -h(d(\alpha )d+dd(\alpha ))+(d(\alpha )d+dd(\alpha ))h+\sum _{\alpha = \beta \gamma }{}^\prime h(\beta )S(\gamma )-S(\beta )h(\gamma ).
$$
\end{lem}
\begin{proof}
\begin{eqnarray*}
&&\sum _{\alpha = \beta \gamma }{}^\prime B(\beta )d(\gamma )-d(\beta )B(\gamma )\\
& = & \sum _{\alpha =\beta \gamma }{}^\prime d(\beta )hd(\gamma )+hd(\beta )d(\gamma )+A(\beta )d(\gamma )-d(\beta )d(\gamma )h-d(\beta )hd(\gamma ) - d(\beta )A(\gamma )\\
& = & -h(d(\alpha )d+dd(\alpha ))+(d(\alpha )d+dd(\alpha ))h+\sum _{\alpha = \beta \gamma }{}^\prime h(\beta )S(\gamma )-S(\beta )h(\gamma ).
\end{eqnarray*}
\end{proof}
We can now prove \ref{L-local}. We compute
\begin{eqnarray*}
&&dA(\alpha )-A(\alpha )d \\
& = & \sum _{\alpha = \beta \gamma }{}^\prime dd(\beta )h(\gamma )+dh(\beta )d(\gamma )-d(\beta )h(\gamma )d-h(\beta )d(\gamma )d\\
&=& \sum _{\alpha = \beta \gamma }{}^\prime dd(\beta )h(\gamma )+(-f(\beta )-B(\beta )-h(\beta )d)d(\gamma )-d(\beta )(-f(\gamma )-B(\gamma )-dh(\gamma ))-h(\beta )d(\gamma )d\\
&=& \sum _{\alpha = \beta \gamma }{}^\prime dd(\beta )h(\gamma )-f(\beta )d(\gamma )-B(\beta )d(\gamma )\\&-&h(\beta )dd(\gamma ) + d(\beta )f(\gamma )+d(\beta )B(\gamma )+d(\beta )dh(\gamma )-h(\beta )d(\gamma )d\\
& = & [\sum _{\alpha = \beta \gamma }{}^\prime (-S(\beta )h(\gamma ))+h(\beta )S(\gamma )-f(\beta )d(\gamma )+d(\beta )f(\gamma )]\\ & + &h(d(\alpha )d+dd(\alpha ))-(d(\alpha )d+dd(\alpha ))h-\sum _{\alpha = \beta \gamma }{}^\prime h(\beta )S(\gamma )-S(\beta )h(\gamma )\\
& = & f(\alpha )d-df(\alpha )+fd(\alpha )-d(\alpha )f+h(d(\alpha )d+dd(\alpha ))-(d(\alpha )d+dd(\alpha ))h.
\end{eqnarray*}
So finally
\begin{eqnarray*}
&& d\Psi (\alpha )-\Psi (\alpha )d \\
&=& df(\alpha )+dd(\alpha )h+dhd(\alpha )+dA(\alpha )-f(\alpha )d-d(\alpha )hd-hd(\alpha )d-A(\alpha )d\\
& = & df(\alpha )+dd(\alpha )h+dhd(\alpha )-f(\alpha )d -d(\alpha )hd-hd(\alpha )d\\
& + &f(\alpha )d-df(\alpha )+fd(\alpha )-d(\alpha )f+hd(\alpha )d
 + hdd(\alpha )-d(\alpha )dh-dd(\alpha )h\\
& = & 0.
\end{eqnarray*}
This completes the proof of \ref{L-local}.
\end{proof}

Let
$$
\epsilon ^*:C(\mathcal A(\emptyset ))\rightarrow C(\text{tot}(\mathcal A^+))
$$
be the functor sending a complex $K$ to the object of $C(\text{tot}(\mathcal A^+))$ with $\epsilon ^*K^{n, *} = K$ with maps $d(\text{id}_{[n]})$ equal to $(-1)^n$ times the differential, for $\partial _i:[n]\rightarrow [n+1]$ the map $d(\partial _i)$ is the canonical map of complexes, and all other $d(\alpha )$'s are zero.  The functor $\epsilon ^*$ takes quasi--isomorphisms to quasi--isomorphisms and hence induces a functor
$$
\epsilon ^*:D(\mathcal A(\emptyset))\rightarrow D(\text{tot}(\mathcal A^+)).
$$

\begin{lem} The functor $\epsilon ^*$ has a right adjoint $R\epsilon _*:D(\text{\rm tot}(\mathcal A^+))\rightarrow D(\mathcal A(\emptyset ))$ and $R\epsilon _*$ is a triangulated functor.
\end{lem}
\begin{proof}
We apply the adjoint functor theorem \cite[4.1]{Neeman}.  By our assumptions the category $D(\mathcal A(\emptyset ))$ is compactly generated.  Therefore it suffices to show that $\epsilon ^*$ commutes with coproducts (direct sums) which is immediate.
\end{proof}

More concretely, the functor $R\epsilon _*$ can be computed as follows.  If $K$ is $L$--local and there exists an integer $N$ such that for every $n$ we have $K^{n, m}=0$ for $m<N$, then $R\epsilon _*K$ is represented by the complex with
$$
(\epsilon _*K)^p = \oplus _{n+m=p}\epsilon _{n*}K^{n, m}
$$
with differential given by $\sum d(\alpha )$. This follows from Yoneda's lemma and the observation that for any $F\in D(\mathcal A(\emptyset))$ we have
\begin{eqnarray*}
\text{Hom}_{D(\mathcal A(\emptyset))}(F, R\epsilon _*K) & = & \text{Hom}_{D(\text{tot}(\mathcal A^+))}(\epsilon ^*F, K)\\
& = & \text{Hom}_{K(\text{tot}(\mathcal A^+))}(\epsilon ^*F, K) \text{ \ \ since $K$ is $L$-local}\\
& = & \text{Hom}_{K(\mathcal A(\emptyset ))}(F, \epsilon _*K) \text{ \ \ by \cite[3.2.12]{BBD82}}.
\end{eqnarray*}

\begin{lem}\label{A9} For any $F\in D^+(\mathcal A(\emptyset ))$ the natural map $F\rightarrow R\epsilon _*\epsilon ^*F$ is an isomorphism.
\end{lem}
\begin{proof} Represent $F$ by a complex of injectives.  Then $\epsilon ^*F$ is $L$--local by \ref{L-local}.  The result then follows from cohomological descent.
\end{proof}

\begin{prop}\label{A10} Let $K\in D^+_{\text{\rm cart}}(\text{\rm tot}(\mathcal A^+))$ be an object.  Then $\epsilon ^*R\epsilon _*K\rightarrow K$ is an isomorphism. In particular, $R\epsilon _*$ and $\epsilon ^*$ induce an equivalence of categories between $D^+_{\text{\rm cart}}(\text{\rm tot}(\mathcal A^+))$ and $D^+(\mathcal A(\emptyset ))$.
\end{prop}
\begin{proof}
For any integer $s$ and  system $(K^{n, *}, d(\alpha ))$ defining an object of $C(\text{tot}(\mathcal A^+))$ we obtain a new object by $(\tau _{\leq s}K^{n, *}, d(\alpha ))$ since for any $\alpha $ which is not the identity morphism the map $d(\alpha )$ has degree $\leq 0$.  We therefore obtain a functor $\tau _{\leq s}:C(\text{tot}(\mathcal A^+))\rightarrow C(\text{tot}(\mathcal A^+))$ which takes quasi--isomorphisms to quasi--isomorphisms and hence descends to a functor
$$
\tau _{\leq s}:D(\text{tot}(\mathcal A^+))\rightarrow D(\text{tot}(\mathcal A^+)).
$$
Furthermore, there is a natural morphism of functors $\tau _{\leq s}\rightarrow \tau _{\leq s+1}$ and we have
$$
K\simeq \text{hocolim }\tau _{\leq s}K.
$$
Note that the functor $\epsilon ^*$ commutes with homotopy colimits since it commutes with direct sums.  If we show the proposition for the $\tau _{\leq s}K$ then we see that the natural map
$$
\epsilon ^*(\text{hocolim}R\epsilon _*\tau _{\leq s}K)\simeq \text{hocolim}\epsilon ^*R\epsilon _*\tau _{\leq s}K\rightarrow \text{hocolim} \tau _{\leq s}K\simeq K
$$
is an isomorphism.  In particular $K$ is in the essential image of $\epsilon ^*$.  Write $K = \epsilon ^*F$.  Then by \ref{A9} $R\epsilon _*K\simeq F$ whence $\epsilon ^*R\epsilon _*K\rightarrow K$ is an isomorphism.

It therefore suffices to prove the proposition for $K$ bounded above.  Considering the distinguished triangles associated to the truncations $\tau _{\leq s}K$ we further reduce to the case when $K$ is concentrated in just a single degree.  In this case, $K$ is obtained by pullback from an object of $\mathcal A(\emptyset )$ and the proposition again follows from \ref{A9}.
\end{proof}

For an object $K\in K(\text{tot}(\mathcal A^+))$, we define $\tau _{\geq s}K$ to be the cone of the natural map $\tau _{\leq s-1}K\rightarrow K$.

Observe that the category $K(\text{tot}(\mathcal A^+))$ has products and therefore also homotopy limits. Let $K\in K_{\mathcal C}(\text{\rm tot}(\mathcal A^+))$ be an object. By \ref{A10}, for each $s$ we can find a bounded below complex of injectives $I_s\in C(\mathcal A(\emptyset ))$ and a quasi--isomorphism $\sigma _s:\tau _{\geq s}K\rightarrow \epsilon ^*I_s$.  Since $\epsilon ^*I_s$ is $L$--local and $\epsilon ^*:D^+(\mathcal A(\emptyset ))\rightarrow D(\text{tot}(\mathcal A^+))$ is fully faithful by \ref{A10}, the maps $\tau _{\geq s-1}K\rightarrow \tau _{\geq s}K$ induce a unique morphism $t_s:I_{s-1}\rightarrow I_s$ in $K(\mathcal A(\emptyset ))$ such that the diagrams
$$
\begin{CD}
\tau _{\geq s-1}K@>>> \tau _{\geq s}K\\
@V\sigma _{s-1}VV @VV\sigma _sV \\
\epsilon ^*I_{s-1}@>t_s>> \epsilon ^*I_s
\end{CD}
$$
commutes in $K(\text{tot}(\mathcal A^+))$.

\begin{prop} The natural map $K\rightarrow \text{\rm holim}\epsilon ^*I_s$ is a quasi--isomorphism.
\end{prop}
\begin{proof}
It suffices to show that for all $n$ the map $K^{n, *}\rightarrow \text{holim} \epsilon _n^*I_s$ is a quasi-isomorphism, where $\epsilon _n:U_n\rightarrow T$ is the projection.  Let $\mathcal S_n$ be a site inducing $U_n$ as in \ref{assumption1b}.
We show that for any integer $i$ the map of presheaves on the subcategory of $\mathcal S_n$ satisfying the finiteness assumption in \ref{assumption1b} (i)
$$
(V\rightarrow U_n)\mapsto H^i(V, K^{n, *})\rightarrow H^i(V, \text{holim} \epsilon _n^*I_s)
$$
is an isomorphism.  For this note that for every $s$ there is a distinguished triangle
$$
\mathcal H^s(K^{n, *})[s]\rightarrow \epsilon _n^*I_s\rightarrow \epsilon _n^*I_{s-1}
$$
and hence by the assumption \ref{assumption1b} (i) the map
\begin{equation}\label{transmap}
H^i(V, \epsilon _n^*I_s)\rightarrow H^i(V, \epsilon _n^*I_{s-1})
\end{equation}
is an isomorphism for $s<i-n_0$.  Since each $\epsilon _n^*I_s$ is a complex of injectives, the complex $\prod _s \epsilon _n^*I_s$ is also a complex of injectives.
Therefore
$$
H^i(V, \prod _s\epsilon _n^*I_s)=H^i(\prod _s\epsilon _n^*I_s(V)) = \prod _sH^i(\epsilon _n^*I_s(V)).
$$
It follows that there is a canonical long exact sequence
$$
\begin{CD}
\cdots @>>> \prod _sH^i(\epsilon _n^*I_s(V))@>1-\text{shift}>> \prod _sH^i(\epsilon _n^*I_s(V))@>>> H^i(V, \text{holim}\epsilon _n^*I_s)@>>> \cdots .
\end{CD}
$$
From this and the fact that the maps \ref{transmap} are isomorphisms for $s$ sufficiently big it follows that the cohomology group $H^i(V, \text{holim}\epsilon _n^*I_s)$ is isomorphic to $H^i(V, K^{n, *})$ via the canonical map.  Passing to the associated sheaves we obtain the proposition.
\end{proof}

\begin{cor}\label{A13} Every object $K\in D_{\mathcal C}(\text{\rm tot}(\mathcal A^+))$ is in the essential image of the functor
$$
\epsilon ^*:D_{\mathcal C}(\mathcal A(\emptyset ))\rightarrow D_{\mathcal C}(\text{\rm tot}(\mathcal A^+)).
$$
\end{cor}
\begin{proof}
Since $\epsilon ^*$ also commutes with products and hence also homotopy limits we find that
$
K\simeq \epsilon ^*(\text{holim }I_s)
$
in $D_{\mathcal C}(\text{tot}(\mathcal A^+))$ (note that $\mathcal H^i(\text{holim}I_s)$ is in $\mathcal C$ since this can be checked after applying $\epsilon ^*$).
\end{proof}

\begin{lem} Let $[n]\mapsto K^n$ be a cartesian section of $[n]\mapsto D(U_n, \Lambda )$ such that $\mathcal Ext^i(K^n, K^n) = 0$ for all $i<0$.  Then $(K^n)$ is induced by an object of $D(\text{\rm tot}(\mathcal A^+))$.
\end{lem}
\begin{proof}
Represent each $K^n$ by a homotopically injective complex (denoted by the same letter) in $C(U_n, \Lambda )$ for every $n$. For each morphism $\partial _i:[n]\rightarrow [n+1]$ (the unique morphism whose image does not contain $i$) choose a $\partial _i$-map of complexes $\partial _i^*:K^n\rightarrow K^{n+1}$ inducing the given map in $D(U_{n+1}, \Lambda )$ by the strictly simplicial structure.  The proof then proceeds by the same argument used to prove \cite[3.2.9]{BBD82}.
\end{proof}

Combining this with \ref{A13} we obtain \ref{1.2}. \qed

\section{Dualizing complex}

\subsection{Dualizing complexes on algebraic spaces} Let $W$ be an algebraic space and $w:W\ra
S$ be a separated\footnote{Probably one can assume only that $w$ quasi-separated,
cf.~\cite{SGA43}, XVII.7; but
we do not need this more general version.} morphism of finite type with
$W$ an algebraic space. We'll define ${\Omega}_w$ by glueing as
follows. By the comparison lemma (\cite{SGA41}, III.4.1), the
\'etale topos $W_\et$ can be defined using the site $\ET(W)$ whose
objects are \'etale morphisms $A:U\ra W$ where $a:U\ra S$ is affine
of finite type. The localized topos $W_{\et|U}$ coincides with
$U_\et$.

Unless otherwise explicitly stated, we will ring the various \'etale
or lisse-\'etale topos which will be appear by the constant
Gorenstein ring $\Lambda$ of dimension $0$ of the introduction.

 Notice that
this is not true for the corresponding lisse-\'etale topos. This
fact will cause some difficulties below. Let $\Omega $ denote the
dualizing complex of $S$, and let $\alpha :U\rightarrow S$ denote
the structural morphism. We define
\begin{equation}\label{def-omega}
    {\Omega}_A=\alpha^!\Omega\in {\D}(U_\et,\Lambda)=\D({W_\et}_{|U}).
\end{equation} which is the (relative) dualizing complex of $U$, and therefore one gets by
biduality (\cite{SGA4.5}, \og Th. finitude\fg\ 4.3)
\begin{equation}\label{OAOA}\Rhom({\Omega}_A,{\Omega}_A)=\Lambda\end{equation}
implying at once \begin{equation}\label{extW}
    \ext^i_{{W_{\et}}_{|U}}({\Omega}_A,{\Omega}_A)=0\si i<0.
\end{equation}

 We want to apply the glueing theorem~\ref{1.2}. Let us therefore consider a diagram
$$\xymatrix{V\ar[rr]^\sigma\ar@{..>}@/_1pc/[rdd]^
\beta\ar[rd]_B&&U\ar@{..>}@/^1pc/[ldd]_\alpha\ar[ld]^A\\
&W\ar[d]\\&S}$$ with a commutative triangle and
$A,B\in\ET(W)$.
\begin{lemma}\label{foncto-KW} There is a functorial
isomorphism
$$\sigma^{*}{\Omega}_A={\Omega}_B.$$
\end{lemma}

\begin{proof} Let $\tilde W=U\times_W V$ : it  is an affine scheme, of
finite type over $S$, and \'etale over both $U,V$. In fact, we have
a cartesian diagram
$$\xymatrix{U\times_WV\ar[r]\ar[d]_\delta&W\ar[d]^\Delta\\U\times_SV\ar[r]&W\times_S W}$$
where $\Delta$ is a closed immersion ($W/S$ separated) showing
that $\tilde W=U\times_WV$ is a closed subscheme of $U\times_SV$
which is affine. Looking at the graph diagram with cartesian
square
$$\xymatrix{\tilde W\ar[r]^b\ar[d]_a&U\ar[d]^A\\V\ar@/^1pc/[u]^s
\ar[r]_B\ar[ru]^\sigma&W},$$
we get that $a,b$ are \'etale and separated like $A,B$. One deduces
a commutative diagram
$$\xymatrix{\tilde W\ar[r]^b\ar[d]_a&U\ar[d]^\alpha\\V\ar@/^1pc/[u]^s
\ar[r]_\beta\ar[ru]^\sigma&S}.$$

We claim that

\begin{equation}
  b^*\alpha^!\Omega=a^*\beta^!\Omega.
\end{equation}
Indeed, $a,b$ being smooth of relative dimension $0$, one has
$$b^*\alpha^!\Omega=b^!\alpha^!\Omega$$
and analogously
$$a^*\beta^!\Omega=a^!\beta^!\Omega.$$
Because $\alpha b=\beta a$, one gets $b^!\alpha^!=a^!\beta^!$.
Pulling back by $s$ gives the
result.\hfill\hfill$\square$\end{proof}

Therefore $({\Omega}_A)_{A\in\ET(W)}$ defines locally an
object ${\Omega_w}$ of
 $\D(W)$ with vanishing negative $\ext$'s (recall that $w:W\rightarrow S$ is the structural morphism). By~\ref{1.2}, we get

 \begin{proposition}\label{KW} There
 exists a unique ${\Omega}_w\in\D(W_\et)$ such that ${\Omega}_{w|U}={\Omega}_A$.
 \end{proposition}

 We need functoriality for smooth morphisms.

 \begin{lemma}\label{fonctW} If $f:W_1\ra W_2$ is a smooth $S$-morphism
of relative dimension $d$ between algebraic space separated and of
finite type over $S$ with dualizing complexes
${\Omega}_1,{\Omega}_2$, then
$$f^*{\Omega}_2={\Omega}_1{\langle}-d{\rangle}.$$\end{lemma}

\begin{proof} Start with $U_2\ra W_2$ \'etale and surjective with $U_2$
affine say. Then, $\tilde W_1=W_1\times_{U_2}W_2$ is an algebraic
space separated and of finite type over $S$. Let $U_1\ra \tilde
W_1$ be a surjective \'etale morphism with $U_1$ affine and let
$g:U_1\ra U_2$ be the composition $U_1\ra \tilde W_1\ra U_2$. It
is a smooth morphism of relative dimension $d$ between affine
schemes of finite type from which follows the formula
$g^!(-)=g^*(-){\langle}d{\rangle}$. Therefore, the pull-backs  of
$L_1={\Omega}_1{\langle}-d{\rangle}$ and $f^*{\Omega}_2$ on $U_1$
are the same, namely ${\Omega}_{U_1}$. One deduces that these
complexes coincide on the covering  sieve ${W_{1\et}}_{|U_1}$ and
therefore coincide by~\ref{D1} (because the relevant negative
$\ext^i$'s vanish.\hfill\hfill$\square$\end{proof}

\subsection{\'Etale dualizing data}
Let $\X\ra S$ be an algebraic $S$-stack locally of finite type. Let $A:U\ra \X$ in
$\LE(\X)$ and $\alpha:U\ra S$ the composition $U\ra\X\ra S$. We
define
\begin{equation}\label{def-KA}
    K_A=\Omega_\alpha{\langle}-d_A{\rangle}\in {\D}_c(U_\et,\Lambda)
\end{equation} where $d_A$ is the relative dimension
of $A$ (which is locally constant). Up to shift and Tate torsion,
$K_A$ is the (relative) dualizing complex of $U$ and therefore one
gets by biduality
\begin{equation}\label{KAKA}\Rhom(K_A,K_A)=\Lambda\textup{ and }\ext^i_{U_\et}(K_A,K_A)=0\si i<0.
\end{equation}

We need again a functoriality property of $K_A$. Let us consider a
diagram
$$\xymatrix{V\ar[rr]^\sigma\ar@{..>}@/_1pc/[rdd]^
\beta\ar[rd]_B&&U\ar@{..>}@/^1pc/[ldd]_\alpha\ar[ld]^A\\
&\X\ar[d]\\&S}$$ with a $2$-commutative triangle and
$A,B\in\LE(\X)$.

\begin{lemma}\label{foncto-K} There is  a functorial identification
$$\sigma^{*}K_A=K_B.$$
\end{lemma}

\begin{proof} Let $W=U\times_\X V$ which is an algebraic space. One has a
commutative diagram with cartesian square

$$\xymatrix{W\ar[r]^b\ar[d]_a&U\ar[d]^A\\V\ar@/^1pc/[u]^s
\ar[r]_B\ar[ru]^\sigma&\X}.$$

In particular, $a,b$ are smooth and separated like $A,B$. One
deduces a commutative diagram
$$\xymatrix{W\ar[r]^b\ar[d]_a&U\ar[d]^\alpha\\V\ar@/^1pc/[u]^s\ar[r]_\beta\ar[ru]^\sigma&S}.$$

I claim that

\begin{equation}
  b^*K_A=a^*K_B=K_w.
\end{equation}
where $w$ denotes the structural morphism $W\ra S$.

Indeed, $a,b$ being smooth of relative dimensions $d_A,d_B$, one
has \ref{fonctW}
$$b^*K_A=b^*\Omega_\alpha{\langle}-d_A{\rangle}=\Omega_\alpha{\langle}-d_A-d_B{\rangle}$$
and analogously
$$a^*K_B=a^*\Omega{\langle}-d_B{\rangle}=\Omega_\alpha{\langle}-d_B-d_A{\rangle}.$$
Because $\alpha b=\beta a$, one gets $b^!\alpha^!=a^!\beta^!$.
Pulling back by $s$ gives the result.
\hfill\hfill$\square$\end{proof}

\begin{remark}\label{finitude-inj_W} Because all $S$-schemes of finite type satisfy
$\textrm{cd}_\Lambda(X)<\infty$, we know that $K_X$ is not only of
finite quasi-injective dimension but of finite injective dimension
(\cite{SGA5}, I.1.5). By construction this implies that $K_A$ is
of finite injective dimension  for $A$ as
above.\end{remark}
\subsection{Lisse-\'etale dualizing data} In order to define $\Omega_\X\in\D(\X_\le)$ by glueing,
we need  glueing data $\kappa_A\in\D(\X_{\le|U}), U\in\LE(\X)$.
The inclusion $$\ET(U)\hookrightarrow\LE(\X)_{|U}$$ induces a
continuous morphism of sites.  Since finite inverse limits exist in $\ET(U)$ and this morphism of sites preserves such limits,
it defines by (\cite{SGA41},
4.9.2) a morphism of topos (we abuse notation slightly and omit the dependence on $A$ from the notation)
$$\epsilon:\X_{\le|U}\ra U_\et.$$

\begin{pg}\label{LETdata} Let us describe more explicitely the morphism $\epsilon$. Let $\LE(\mc X)_{|U}$ denote the category of morphisms $V\rightarrow U$ in $\LE (\mc X)$.  The category $\LE(\mc X)_{|U}$ has  a Grothendieck topology induced by the topology on $\LE(\mc X)$, and the resulting topos is canonicallly isomorphic to the localized topos $\X_{\le|U}$.  Note that there is a natural inclusion $\LE(U)\hookrightarrow \LE(\mc X)_{|U}$ but this is not an equivalence of categories since for an object $(V\rightarrow U)\in \LE(\mc X)_{|U}$ the morphism $V\rightarrow U$ need not be smooth.
It follows that an element of $\X_{\le|U}$ is equivalent to giving
for every $U$--scheme of finite type $V\rightarrow U$, such that
the composite $V\rightarrow U\rightarrow \X$ is smooth, a sheaf
$\F_V\in V_{\et }$ together with morphisms $f^{-1}\F_V\rightarrow
\F_{V'}$ for $U$--morphisms $f:V'\rightarrow V$.  Furthermore,
these morphisms satisfy the usual compatibility with compositions.
Viewing $\X_{\le|U}$ in this way, the functor $\epsilon^{-1}$ maps
$\F$ on $U_\et$ to $\F_V=\pi^{-1}\F\in V_\et$ where $\pi:V\ra
U\in\LE(\X)_{|U}$. For a sheaf $F\in \X_{\le |U}$ corresponding to
a collection of sheaves $\F_V$, the sheaf $\epsilon _*F$ is simply
the sheaf $\F_U$.

In particular, the functor $\epsilon _*$ is exact and,
accordingly, that $H^*(U,F)=H^*(U_{\et},F_U)$ for any shaf of
$\Lambda$ modules of $\X$.
\end{pg}

\begin{pg}A morphism $f:U\ra V$ of $\LE(\X)$ induces a diagram
\begin{equation}\label{main-diag-loc}
\begin{CD}
\mc X_{\le}|_U@>\epsilon >> U_{\et}\\
@VfVV @VVV \\
\mc X_{\le}|_{V}@>\epsilon >> V_{\et}
\end{CD}
\end{equation} where $\X_{\le|U}\ra \X_{\le|V}$
is the localization morphism (\cite{SGA41}, IV.5.5.2) which we
still denote by $f$ slightly abusively. For a sheaf $\F\in V_{\et }$, the pullback $f^{-1}\epsilon^{-1}\F$ is the sheaf corresponding to the system which to any $p:U'\rightarrow U$ associates $p^{-1}f^{-1}\F$.  In particular, $f^{-1}\circ \epsilon ^{-1} = \epsilon ^{-1}\circ f^{-1}$ which implies that \ref{main-diag-loc} is a commutative diagram of topos.
 We define
\begin{equation}\label{def-dual-loc}
    \kappa_A=\epsilon^*K_A\in\D(\X_{\le|U}).
\end{equation}
By the preceding discussion, if
$$\xymatrix{U\ar[rr]^f\ar[rd]_A&&V\ar[ld]^B\\&\X}$$ is a morphism
in $\LE(\X)$, we get $$f^*\kappa_B=\kappa_A$$ showing that the
family $(\kappa_A)$ defines locally an object of
$\D(\X_\le)$.\end{pg}

\subsection{Glueing the local dualizing data} Let $A\in\LE(\X)$ and
 $\epsilon:\X_{\le|U}\ra U_\et$
be as above. We need first the vanishing of
$\ext^i(\kappa_A,\kappa_A), i<0$.

\begin{lemma}\label{tauto-lisse-etale} Let $\F,\G\in\D(U_\et)$.
One has\begin{itemize}
    \item[(i)]$\Ext^i({\epsilon}^{*}{\F},{\epsilon}^{*}{\G})=\Ext^i({\F},{\G})$.
    \item[(ii)] The \'etale
    sheaf
    $\ext^i({\epsilon}^{*}{\F},{\epsilon}^{*}{\G})_U$ on $U_\et$ is $\ext^i_{U_\et}({\F},{\G})$.
\end{itemize}
\end{lemma}
\begin{proof} Since $\epsilon _*$ is exact and for any sheaf $F\in U_\et $ one has $F=\epsilon _*\epsilon ^*F$, the adjunction map $F\rightarrow R\epsilon _*\epsilon ^*F$ is an isomorphism for any $F\in \D(U_\et )$.
By trivial duality, one gets
$${\epsilon}_*\Rhom({\epsilon}^{*}{\F},{\epsilon}^{*}{\G})=\Rhom({\F},{\epsilon}_*{\epsilon}^{*}{\G})=\Rhom({\F},{\G}).$$
Taking ${\mc H}^iR\Gamma$ gives \textit{(i)}.

By construction, $\ext^i(\epsilon^*\F,\epsilon^*\F)_U$ is the sheaf
associated to the presheaf on $U_{\et }$ which to any \'etale morphism $\pi :V\rightarrow U$ associates
$\Ext^i(\pi^*\epsilon^*\F,\pi^*\epsilon^*\G)$ where
$\pi^*$ is the the pull-back functor associated to the
localization morphism $$(\X_\le|_U)_{|V}=\X_{\le|V}\ra\X_{\le|U}$$
(\cite{SGA42}, V.6.1). By the commutativity of the
diagram~\ref{main-diag-loc}, one has
$\pi^*\epsilon^*=\epsilon^*\pi^*$. Therefore
$$\Ext^i(\pi^*\epsilon^*\F,\pi^*\epsilon^*\G)=\Ext^i(\epsilon^*\pi^*\F,\epsilon^*\pi^*\G)=
\Ext^i_{V_\et}(\pi^*\F,\pi^*\G),$$ the last equality is by
{\it(i)}. Since $\ext_{U_\et}(\F, \G)$ is also the sheaf
associated to this presheaf we obtain {\it
(ii)}.\hfill\hfill$\square$\end{proof}

 Using~\ref{KAKA}, one obtains
\begin{corollary}\label{kAkA} One has $\Rhom(\kappa_A,\kappa_A)=\Lambda$ and
therefore $\ext^i(\kappa_A,\kappa_A)=0$ if $i<0$.\end{corollary}

The discussion above shows that we
can apply~\ref{1.2} to $(\kappa_A)$ to get

\begin{proposition}\label{exis-dual}  There exists
${\Omega}_\X(p)\in {\D}^b(\X_\le)$ inducing $\kappa_A$ for all
${A\in\LE(\X)_{|X}}$. It is well defined up to unique
isomorphism.\end{proposition}

The independence of the presentation is straightforward and is
left to the reader :

\begin{lemma}\label{indep-K}Let $p_i:X_i\ra\X,i=1,2$ two presentations as above. There exists a
canonical, functorial isomorphism ${\Omega}_\X(p_1)\isom
{\Omega}_\X(p_2)$.
\end{lemma}

\begin{definition}\label{def-dualisant} The \emph{dualizing complex of $\X$} is the "essential" value $\Omega _{\X}\in \D^b(\X_{\le})$ of ${\Omega_\X}(p)$, where
$p$ runs over presentations of $\X$.
It is well defined up to canonical functorial isomorphism and is
characterized by $\Omega_{\X|U}=\epsilon^*K_A$ for any $A:U\ra\X$
in $\LE(\X)$.
\end{definition}

\subsection{Biduality} For $A,B$ any abelian complexes of some
topos, there is a biduality \textit{morphism}
\begin{equation}\label{bidual}
    A\ra\Rhom(\Rhom(A,B),B)
\end{equation} (replace $B$ by some homotopically
injective complex isomorphic to it in the derived catgory).

In general, it is certainly not an isomorphism.
\begin{lemma}\label{bidualnonborne} Let $u:U\ra S$ be a separated $S$-scheme (or algebraic space) of finite type and $A\in
{\D}_c(U_\et,\Lambda)$. Then the biduality morphism
$$A\ra\Rhom(\Rhom(A,K_U),K_U)$$ is an isomorphism
(where $K_U$ is -up to shift and twist- the dualizing complex of
$U_\et$).
\end{lemma}
\begin{proof} If $A$ is moreover bounded, it is the usual theorem of
~\cite{SGA4.5}. Let us denote by $\tau_n$ the two-sides truncation
functor $$\tau_{\geq -n}\tau_{\leq n}.$$ We know that $K_U$ is a
dualizing complex (\cite{SGA5}, exp. I), and is of \emph{finite
injective dimension} (\ref{finitude-inj_W}); the homology in
degree $n$ of the biduality morphism $A\ra DD(A)$ is therefore the
same as the homology in degree $n$ of the biduality morphism
$\tau_mA\ra DD(\tau_m A)$ for $m$ large enough and the lemma
follows.\hfill\hfill$\square$\end{proof}

We will be interested in a commutative diagram
$$\xymatrix{V\ar[rr]^f\ar[rd]_B&&\ar[ld]^AU\\&\X}$$ as above.

\begin{lemma}\label{tec1} Let $\F\in\D_c(\X_{\le})$ and let $\F_U\in D_c(U_\et )$ be the object obtained by restriction.\begin{itemize}
    \item[(i)] One has
    $f^*\Rhom(\F_U,K_A)=\Rhom(f^*\F_U,f^*K_A)=\Rhom(f^*\F_U,K_B)$.
    \item[(ii)] Moreover, $\Rhom(\F_U,K_A)$ is constructible.
\end{itemize}
\end{lemma}

\begin{proof} Let's prove (\textit{i}). By~\ref{foncto-K}, one has
$f^*K_A=K_B$, therefore one has a morphism
$$f^*\Rhom(\F_U,K_A)\ra\Rhom(f^*\F_U,K_B).$$

To prove that it is an isomorphism, consider first the case when $f$ is smooth.
Because both $K_A$ and
$K_B$ are of finite injective dimension (\ref{finitude-inj_W}),
one can assume  that $F$ is bounded where  it is obviously true by
reduction to $F$ the constant sheaf (or use~\cite{SGA5}, I.7.2). Therefore the result holds when $f$ is smooth.

From the case of a smooth morphism, one reduces the proof in general to the case when $\X$ is a scheme.  Let $\F_{\X}\in D_c(\X_\et)$ denote the complex obtained by restricting $\F$.  By the smooth case already considered, we have
\begin{eqnarray*}
f^*\Rhom(\F_U,K_A)& \simeq & f^*A^*\Rhom (\F_\X, K_\X)\\ & = & B^*\Rhom (\F_\X, K_\X)\\
& \simeq & \Rhom (B^*\F_X, B^*K_\X)\\ & \simeq &\Rhom (f^*\F_U, f^*K_A).
\end{eqnarray*}

For (\textit{ii}), one can also assume $\F$ bounded and one
uses~\cite{SGA5}, I.7.1.\hfill\hfill$\square$\end{proof}

\begin{lemma}\label{tec2} Let $\F\in\D_c(\X_\le)$. Then,
$$\epsilon^*\Rhom_{U_\et}(\F_U,K_A)=\Rhom(\F,{\Omega_\X})_{|U}$$
where $\F_U=\epsilon_*\F_{|U}$ is the restriction of $\F$ to
$\ET(U)$.\end{lemma} \begin{proof} By definition of
constructibility, ${\mc H}^i(\F)$ are cartesian sheaves. In other
words, $\epsilon_*$ being exact, the adjunction morphism
$$\epsilon^*\F_U=\epsilon^*\epsilon_*\F_{|U}\ra\F_{|U}$$ is an isomorphism. We therefore
have\begin{eqnarray*}
  \Rhom(\F,\Omega)_{|U} &=& \Rhom(\F_{|U},\Omega_{|U})\\
  &=& \Rhom(\epsilon^*\F_U,\epsilon^* K_A) \\
\end{eqnarray*}
Therefore, we get a morphism
$$\epsilon^*\Rhom_{U_\et}(\F_U,K_A)\ra\Rhom(\epsilon^*\F_U,\epsilon^* K_A)=\Rhom(\F,{\Omega_\X})_{|U}.$$
By~\ref{tauto-lisse-etale}, one has
$$\ext^i({\epsilon}^{*}{\F_U},{\epsilon}^{*}K_A)_V=
\ext^i_{V_\et}(f^{*}{\F_U},f^{*}K_A).$$ But, one has
$${\mc H}^i(\epsilon^*\Rhom_{U_\et}(\F_U,K_A))_V=f^*\ext^i_{U_\et}(\F_U,K_A))$$
and the lemma follows
from~\ref{tec1}.\hfill\hfill$\square$\end{proof}

 One gets immediately (cf.~\cite{SGA5}, I.1.4)

\begin{corollary}\label{diminffinie} ${\Omega_\X}$ is of finite quasi-injective
dimension.\end{corollary}

\begin{remark} It seems over-optimistic to think that ${\Omega_\X}$ would be
of finite injective dimension even if $\X$ is a scheme.
\end{remark}

\begin{lemma}\label{dula-construct} If $A\in D_c(\X)$, then $\Rhom(A,{\Omega_\X})\in
D_c(\X)$.\end{lemma} \begin{proof} Immediate consequence of
\ref{tec1} and~\ref{tec2}.\hfill\hfill$\square$\end{proof}

\begin{corollary}\label{Dinvolutif} The (contravariant) functor
$$D_\X:\left\{\begin{array}{ccc}{\D}_c(\X)&\ra&{\D}_c(\X)\\\F&\mapsto&\Rhom(\F,{\Omega_\X})\end{array}\right.$$
is an involution. More precisely, the morphism $$\iota:\Id\ra
{\D}_\X\circ {\D}_\X$$ induced by~\ref{bidual} is an isomorphism.
\end{corollary}

\begin{proof} We have to prove that the cone $C$ of the biduality
morphism is zero in the derived category, that is to say
$$C_U=\epsilon_*C_{|U}=0\text{ in }\D_c(U_\et).$$

But we have\begin{eqnarray*}
  \epsilon_*(\Rhom(\Rhom(\F,{\Omega_\X}),{\Omega_\X}))_{|U}
  &=&  \epsilon_*\Rhom(\Rhom(\F,{\Omega_\X})_{|U},\Omega_{\X|U}) \\
   &\stackrel{\ref{tec2}}{=}& \epsilon_*\Rhom(\epsilon^*\Rhom(\F_U,K_A),\Omega_{\X|U})\\
   &=&\Rhom(\Rhom(\F_U,K_A),\epsilon_*\epsilon^*K_A)\text{ by trivial duality}  \\
  &=&\Rhom(\Rhom(\F_U,K_A),K_A)\\
  &\stackrel{\ref{bidualnonborne}}{=}&\F_U
\end{eqnarray*}
\hfill\hfill$\square$\end{proof}

\begin{remark} Verdier duality $D_\X$ identifies $D^a_c$ and
$D_c^{-a}$ with $a=\varnothing,\pm 1,b$ and the usual conventions
$-\varnothing=\varnothing$ and $-b=b$.
\end{remark}

\begin{proposition}\label{Danti} One has a canonical (bifunctorial)
morphism
$$\Rhom(A,B)=\Rhom(D(B),D(A))$$ for all $A,B\in
D_c(\X)$.
\end{proposition}
\begin{proof} Let us prove first a well-known formula
\begin{lemma}\label{bilder} Let $A,B,C$ be complexes of $\Lambda$ modules on $\X_\le$.
 One has canonical identifications $$\Rhom(A,\Rhom(B,C))=\Rhom(A\Otimes
B,C)=\Rhom(B,\Rhom(A,C)).$$
\end{lemma}
\begin{proof} One can assume $A,B$ homotopically flat and $C$
homotopically injective. Let $X$ be an acyclic complex. One has
$$\Hom(X,\homo(B,C))=\Hom(X\otimes B,C).$$
Because $B$ est homotopically flat, $X\otimes B$ is acyclic.
Moreover, $C$ being homotopically injective, the abelian complex
$\Hom(X\otimes B,C)$ is acyclic. Therefore, $\homo(B,C)$
homotopically injective. One gets therefore
$$\Rhom(A,\Rhom(B,C))=\homo(A,\homo(B,C))=\homo(A\otimes
B,C)=\Rhom(A\Otimes B,C).$$\hfill\hfill$\square$\end{proof}

One gets then
$$\Rhom(D(B),D(A))=\Rhom(D(B),\Rhom(A,\Omega _\X))\stackrel{\ref{bilder}}{=}\Rhom(A,D\circ
D(B))=\Rhom(A,B).$$\hfill\hfill$\square$\end{proof}

\section{The $6$ operations}

\subsection{The functor $\Rhom (-, -)$}

Let $\X$ be an $S$--stack locally of finite type.  As in any
topos, one can define internal hom $\Rhom _{\X_{\le }}(F, G)$ for
any $F\in \D^-(\X)$ and $G\in \D^+(\mc X)$.

\begin{lem} Let $F\in \D_c^-(\X)$ and $G\in \D_c^+(\X)$, and let $j$ be an integer.  Then the restriction of the sheaf $\mc H^j(\Rhom _{\X _\le }(F, G))$ to the \'etale topos of any object $U\in \LE (\X)$ is canonically isomorphic to $\ext ^j_{U_\et }(F_U, G_U)$, where $F_U$ and $G_U$ denote the restrictions to $U_\et $.
\end{lem}
\begin{proof}
The sheaf $\mc H^j(\Rhom _{\mc X_{\le }}(F, G))$ is the sheaf
associated to the presheaf which to any smooth affine $\mc
X$--scheme $U$ associates $\text{Ext}^j_{\mc X_{\le |U}}(F, G)$,
where $\mc X_{\le |U}$ denotes the localized topos. Let $\epsilon
:\mc X_{\le |U}\rightarrow U_{\et}$ be the morphism of topos
induced by the inclusion of $\ET(U)$ into $\LE(\mc X)_{|U}$.  Then
since $F$ and $G$ have constructible cohomology,  the natural maps
$\epsilon ^*\epsilon _*F\rightarrow F$ and $\epsilon ^*\epsilon
_*G\rightarrow G$ are isomorphisms in $\D(\mc X_{\le |U})$.  By
the projection formula it follows that
$$
\text{Ext}^j_{\mc X_{\le |U}}(F, G)\simeq \text{Ext}^j_{\mc X_{\le |U}}(\epsilon ^*\epsilon _*F, \epsilon ^*\epsilon _*G)\simeq \text{Ext}^j_{U_{\et }}(\epsilon _*F, \epsilon _*G).
$$
Sheafifying this isomorphism we obtain the isomorphism in the
lemma. \hfill\hfill$\square$\end{proof}

\begin{cor} If $F\in \D^-_c(\X)$ and $G\in \D^+_c(\X)$, the complex $\Rhom _{\mc X_\le }(F, G)$ lies in $\D_c^+(\mc X)$.
\end{cor}
\begin{proof} By the previous lemma and the constructibility of
the cohomology sheaves of $F$ and $G$, it suffices to prove the
following statement: Let $f:V\rightarrow U$ be a smooth morphism
of schemes of finite type over $S$, and let $F\in \D_c^-(U_\et )$
and $G\in \D_c^+(U_\et )$. Then the natural map $f^*\Rhom _{U_\et
}(F, G)\rightarrow \Rhom _{V_{\et }}(f^*F, f^*G)$ is an
isomorphism as we saw in the proof of \ref{tec1}  (see
\cite{SGA5}, I.7.2). \hfill\hfill$\square$\end{proof}

\begin{prop}\label{homdescription} Let $X/S$ be an $S$--scheme locally of finite type and
$X\rightarrow \X$ be a smooth surjection. Let $X_\bullet
\rightarrow \X$ be the resulting strictly simplicial space.  Then
for $F\in \D_c^-(\X_{\le })$ and $G\in \D_c^+(\X_{\le })$ there is
a  canonical isomorphism
\begin{equation}
\Rhom _{\X_{\le }}(F, G)|_{X_\bullet, \et }\simeq \Rhom _{X_{\bullet \et }}(F|_{X_{\bullet, \et}}, G|_{X_{\bullet, \et }}).
\end{equation}
In particular, $\Rhom _{X_{\bullet \et }}(F|_{X_{\bullet, \et}}, G|_{X_{\bullet, \et }})$ maps under the equivalence of categories $\D_c(X_{\bullet, \et })\simeq \D_c(\mc X_{\le })$ to $\Rhom _{\X_{\le }}(F, G)$.
\end{prop}
\begin{proof} Let $\mc X_{\le |X_\bullet }$ denote the strictly simplicial localized topos and consider the morphisms of topos
\begin{equation}
\begin{CD}
\X_{\le }@<\pi << \X_{\le |X_\bullet }@>\epsilon >> X_{\bullet , \et}.
\end{CD}
\end{equation}
Let $F_\et := \epsilon _*\pi ^*F$ and $G_\et := \epsilon _*\pi
^*G$. Since $F,G\in \D_c(\mc X_{\le })$, the natural maps $F
\simeq R\pi _*\epsilon ^*F_\et $ and $G \simeq R\pi _*\epsilon
^*G$ are isomorphisms (\ref{mainthm}). Using the projection
formula we then obtain
\begin{eqnarray*}
\Rhom _{\X_{\le }}(F, G)|_{X_\bullet, \et } & \simeq  & \epsilon _*\pi ^*\Rhom _{\X_{\le }}(F, G)\\
& \simeq & \epsilon _*\pi ^*\pi _*\Rhom _{\X_{\le |X_\bullet }}(\epsilon ^*F_\et , \epsilon ^*G_\et )\\
& \simeq & \epsilon _*\Rhom _{\X_{\le |X_\bullet }}(\epsilon ^*F_\et , \epsilon ^*G_\et )\\
& \simeq & \Rhom _{X_{\bullet, \et }}(F_\et , G_\et ).
\end{eqnarray*}
\hfill\hfill$\square$\end{proof}

\subsection{The functor $f^*$} The lisse-\'etale
site is not functorial (cf.~\cite{Ber03}, 5.3.12): a morphism of
stacks does not induce a general a morphism between corresponding
lisse-\'etale topos. In~\cite{Ols05}, a functor $f^*$ is
constructed on $\D_c^+$ using cohomological descent. Using the
results of~\ref{mainthm} which imply that we have cohomological
descent also for unbounded complexes, the construction of
\cite{Ols05} can be used to define $f^*$ on the whole category
$\mc D_c$.

Let us review the construction here. Let $f:\X\ra\Y$ be a morphism
of  algebraic $S$--stacks locally of finite type. Choose a
commutative diagram
$$\xymatrix{X\ar[r]\ar[d]&\X\ar[d]\\Y\ar[r]&\Y}$$ where the
horizontal lines are presentations inducing a commutative diagram
of strict simplicial spaces
$$\xymatrix{X_{\bullet}\ar[r]^{\eta_X}\ar[d]_{f_\bullet}&
\X\ar[d]^f\\Y_{\bullet}\ar[r]^{\eta_Y}&\Y.}$$ We get a diagram of
topos

$$\xymatrix{X_{\bullet,\et}\ar[d]_{f_\bullet}&\mc X_{\le }|_{X_{\bullet }}\ar[r]^{\eta_X}\ar[l]_-{\Phi_X}&
\X_\le\\Y_{\bullet,\et}&\ar[l]_-{\Phi_Y}\mc Y_{\le }|_{Y_{\bullet }}\ar[r]^{\eta_Y}&\Y_\le.}$$

By \ref{etexample} the horizontal morphisms induce equivalences of topos
$$
\D_c(\mc X_{\le })\simeq \D_c(X_{\bullet , \et }), \ \ \ \D_c(\mc Y_{\le })\simeq \D_c(Y_{\bullet, \et}).
$$
We define the functor $f^*:\D_c(\mc Y_{\le })\rightarrow \D_c(\mc X_{\le })$ to be the composite
\begin{equation}
\begin{CD}
\D_c(\mc Y_{\le })\simeq \D_c(Y_{\bullet , \et })@>f_\bullet ^*>>  \D_c(X_{\bullet, \et})\simeq \D_c(\mc X_{\le }),
\end{CD}
\end{equation}
where $f_\bullet ^*$ denotes the derived pullback functor induced
by the morphism of topos $f_\bullet :X_{\bullet, \et }\rightarrow
Y_{\bullet , \et }$. Note that $f^*$ takes distinguished triangles
to distinguished triangles since this is true for $f_{\bullet
}^*$.

\begin{prop}\label{adjoint*} Let $A\in \D_c^-(\Y)$ and let $B\in \D_c^+(\X)$. Then there is a canonical isomorphism
\begin{equation}
f_*\Rhom (f^*A, B)\simeq \Rhom (A, f_*B).
\end{equation}
where we write $f_*$ for $Rf_*$.
\end{prop}
\begin{proof}
By \ref{homdescription} and \cite{Ols05}, we have
\begin{equation*}
Rf_*\Rhom (f^*A, B)|_{Y_{\bullet, \et }}\simeq Rf_{\bullet *}\Rhom
_{X_{\bullet , \et }}(f_{\bullet }^*A|_{Y_{\bullet, \et }},
B|_{X_{\bullet, \et }}).
\end{equation*}
The result therefore follows from the usual adjunction
\begin{equation}
Rf_{\bullet *}\Rhom _{X_{\bullet , \et }}(f_\bullet
^*(A|_{Y_{\bullet , \et }}), B|_{X_{\bullet \et }})\simeq \Rhom
_{Y_{\bullet, \et }}(A|_{\bullet, \et }, f_*B|_{X_{\bullet \et
}}).
\end{equation}

\hfill\hfill$\square$\end{proof}

\begin{remark} Its definitely hopeless to generalize \ref{adjoint*} to
$B\in\D_c(\X)$ because in general $Rf_*$ does not map $\D_c$ to itself
(for example consider $B\mathbb{G}_m\ra\Spec(k)$ and $B=\oplus_{i\geq
0}\Lambda[i]$).\end{remark}

\begin{remark} One can even show that \ref{adjoint*} still holds for arbitrary $A\in\D_c(\Y)$, but the geometric significance is unclear because it
is an equality of non constructible complexes.
\end{remark}

\subsection{Definition of $Rf_!,f^!$}\label{4.3}

Let $f:\X\ra\Y$ be a morphism of stacks (locally of finite type over $S$) of finite type. Recall (\cite{Lau-Mor2000}, corollaire
18.4.4) that $Rf_*$ maps $D_c^+(\X_\le)$ to $D_c^+(\Y_\le)$.

\begin{definition}\label{deff!}
We define $$Rf_!: D_c^-(\X_\le)\ra D_c^-(\Y_\le)$$ by the formula
$$Rf_!=D_\Y\circ Rf_*\circ  D_\X,$$  and
$$f^!: D_c^-(\Y_\le)\ra D_c^-(\X_\le)$$ by the formula
$$f^!=D_\X\circ  f^{*}\circ  D_\Y.$$
\end{definition}

By construction, one has\begin{equation}\label{fonct-ext-dual}
    f^!\Omega_\Y=\Omega_\X.
\end{equation}

\begin{proposition}\label{dual} Let $A\in D_c^-(\X_\le)$ and $B\in D_c^-(\Y_\le)$.
Then there is a  (functorial) adjunction formula
$$Rf_*\Rhom(A,f^!B)=\Rhom(Rf_!A,B).$$
\end{proposition}
\begin{proof} We write  $D$ for $D_\X,D_\Y$ and $A'=D(A)\in
D_c^+(\X)$. One has
\begin{align*}
\Rhom(Rf_!D(A'),B) & =\Rhom(D(Rf_*A'),B) \\
  & =\Rhom(D(B),Rf_*A')\ (\ref{Danti})\\
  &=Rf_*\Rhom(f^*D(B),A')\ (\ref{adjoint*})\\
  &=Rf_*\Rhom(D(A'),f^!B)\ (\ref{Danti})
\end{align*}
\hfill\hfill$\square$\end{proof}

\subsection{Projection formula}

\begin{lemma}\label{otimes=rhom} Let $A,B\in D_c(\X)$.
\begin{enumerate}
  \item [(i)] One has
$$\Rhom(A,B)=D_\X(A\Otimes D_\X(B)).$$
  \item [(ii)] If $A,B\in D^-_c(\X )$, then $A\Otimes B\in D^-_c(\X )$.
  \item [(iii)] If $A\in D^-_c(\X ),B\in D^+_c(\X )$, then $\Rhom(A,B)\in D^+_c(\X)$.
\end{enumerate}
\end{lemma}
\begin{proof} Let $\Omega _\X$ be the dualizing complex of $\X$.
\begin{align*}
  \Rhom(A,B) &=\Rhom(D_\X(B),\Rhom(A,\Omega _\X))\ (\ref{Danti}) \\
   & =\Rhom(D_\X(B)\Otimes A,\Omega _\X)\ (\ref{bilder})\\
   &=D_\X(A\Otimes D_\X(B) )
\end{align*}
proving (i). For (ii), using truncations, one can assume that
$A,B$ are sheaves : the result is obvious in this case.  Statement
(iii) follows from the two previous points.
\hfill\hfill$\square$\end{proof}

\begin{corollary}\label{corkunneth} Let $f:\X\rightarrow \Y$ be a morphism as in \ref{4.3}, and let  $B\in D^-_c(\Y), A\in D^-_c(\X)$.
One has the projection formula
$$R f_!(A\Otimes f^*B)=R f_!A\Otimes B.$$
\end{corollary}
\begin{proof} Notice that the left-hand side is well defined
by~\ref{otimes=rhom}. One has
\begin{align*}
  R f_!(A\Otimes f^*B) & =D_\Y \circ R f_*\circ D_\X (A\Otimes D_\X f^!D_\Y B)\\
   & =D_\Y \circ R f_*(\Rhom(A,f^!D_\Y B))\ (\ref{otimes=rhom})\\
   &=D_\Y (\Rhom(R f_!A,D_\Y B))\ (\ref{dual})\\
   &=R f_!A\Otimes B\ (\ref{otimes=rhom})\text{ and}\ (\ref{Dinvolutif}).
\end{align*}
\hfill\hfill$\square$\end{proof}

\begin{corollary}\label{f*f!} For all $A\in D_c^+(\Y),B\in D^-_c(\Y)$, one has
$f^!\Rhom(A,B)=\Rhom(f^*A,f^!B)$.\end{corollary} \begin{proof} By
lemma~\ref{otimes=rhom} and biduality, the formula reduces to the
formula $$f^*(A\Otimes D(B))=f^*A\Otimes f^*D(B).$$ Using suitable
presentation, one is reduced to the obvious formula
$$f_\bullet^*(A_\bullet\Otimes B_\bullet)=f^*_\bullet A_\bullet\Otimes f^*_\bullet B_\bullet$$ for
a morphism $f_\bullet$ of stricltly simplicial \'etale topos.
\hfill\hfill$\square$\end{proof}

\subsection{Computation of $f^!$ for $f$ smooth} Let
$f:\X\ra\Y$ be a smooth morphism of  stacks of relative
dimension $d$. Using~\ref{D1}, one gets immediately  the formula
$$f^*\Omega_\Y=\Omega_\X\langle -d\rangle$$
(choose a presentation of $Y\ra\Y$ and then a presentation
$X\ra\X_Y$; the morphism $X\ra Y$ being smooth, one checks that
these two complexes coincide on $\X_{\le|X}$ and have zero negative
$\ext$'s).

\begin{lemma}\label{Rhom-im-inv} Let $A\in D_c(\Y)$. Then, the canonical morphism $$f^*\Rhom(A,\Omega_\Y)\ra
\Rhom(f^*A,f^*\Omega_\Y)$$ is an isomorphism.
\end{lemma}
\begin{proof} Using~\ref{tauto-lisse-etale}, one is reduced to the usual
statement for \'etale sheaves on algebraic spaces. Because, in
this case, both $\Omega_\Y$ and $f^*\Omega_\Y$ are of finite
injective dimension, one can assume that $A$ is bounded or even a
sheaf. The assertion is well-known in this case (by d\'evissage,
one reduces to $A=\Lambda_\Y$ in which case the assertion is
trivial, cf.~\cite{SGA5}, exp. I).\hfill\hfill$\square$\end{proof}

\begin{corollary}\label{f!lisse} Let
$f:\X\ra\Y$ be a smooth morphism of stacks of relative
dimension $d$. One has $f^!=f^*\langle d\rangle$.\end{corollary}

Let $j:\U\ra\X$ be an open immersion.  Let us denote for a while
$j_{\!}$ the extension by zero functor : it is an exact functor on
the category sheaves preserving constructibility and therefore
passes to the derive category $\D_c$.

\begin{proposition}\label{!=!!}
One has $j^!=j^*$ and $j_!=j_{\!}$.
\end{proposition}
\begin{proof} The first equality is a particular case of~\ref{f!lisse}.
Because $j^*$ has a left adjoint $j_{\!}$ which is exact, it
preserves (homotopical) injectivity. Let $A,B$ be constructible
complexes on $\U,\X$ respectively and assume that $B$ is
homotopically injective. One has
\begin{eqnarray*}
  \RHom(j_{\!}A,B) &=& \Hom(j_{\!}A,B) \\
   &=& \Hom(A,j^*B)\textup{ (adjunction)} \\
   &=&\RHom(A,j^*B)
\end{eqnarray*}  Taking ${\mc H}^0$, one obtains that $j^*$ is the right
adjoint of $j_{\!}$ proving the lemma because $j^!=j^*$ is the
right adjoint of $j_!$.\hfill\hfill$\square$\end{proof}

\subsection{Computation  of $Ri_!$ for $i$ a closed immersion}
 Let $i:\X\hookrightarrow\Y$ be a closed immersion and
${\U}=\Y-\X\hookrightarrow \Y$ the open immersion of the
complement : both are representable. We define the cohomology with
support on $\X$ for any $F\in\X_{\le}$ as follows. First, for any
$Y\ra\Y$ in $\LE(\Y)$, the pull-back $Y_{\U}\ra{\U}$ is in
$\LE({\U})$ and $Y_{\U}\ra{\U}\ra\Y$ is in $\LE(\Y)$. Then, we
define $\underline{H}^0_\X(F)$
$$\Gamma(Y,\underline{H}^0_\X(F))=\Ker(\Gamma(Y,F)\ra\Gamma(Y_{\U},F))$$
and $R\Gamma_\X$ is the total derived functor of the left exact
functor $\Hu_\X^0$.

\begin{lemma}\label{dual-imm-ferme} One has $\Omega_\X=i^*R\Gamma_\X(\Omega_\Y)$.
\end{lemma}
\begin{proof} If $i$ is a closed immersion of schemes (or algebraic
spaces), one has a canonical (and functorial)  isomorphism, simply
because $i^*\underline{H}^0_\X$ is the right adjoint of $i_*$. If
$K$ denotes one of the objects on the two sides of the equality to
be proven, one has therefore $\ext^i(K,K)=0$ for $i<0$. Therefore,
these isomorphisms glue (use theorem 3.2.2 of \cite{BBD82} as
before).\hfill\hfill$\square$\end{proof}

\begin{proposition}\label{prop-dual-imm-ferme} The functor $B\mapsto i^*R\Gamma_\X(B)$
is the right adjoint of $i_*$, and therefore coincides with $i^!$.
More generally, one has
$$\Rhom(i_*A,B)=i_*\Rhom(A,i^*R\underline{H}^0_\X(B))$$ for all $A\in
D(\X),B\in D(\Y).$
 Moreover, one has
one has $i_!=i_*$ and has a right adjoint, the sections with
support on $\X$.

\end{proposition} \begin{proof} If $A,B$ are sheaves, one has the usual
adjunction formula
$$\hom(i_*A,B)=i_*\hom(A,i^*\underline{H}^0_\X(B)).$$ Because $i_*$ is
exact, it's right adjoint sends homotopically injective complexes to
homotopically injective complexes. The derived version follows. One gets
therefore

$$\begin{array}{ccll}
  i_!A &=&\Rhom(i_*\Rhom(A,\Omega_\X),\Omega_\Y ) \\
  &=& i_*\Rhom(\Rhom(A,\Omega_\X),i^*R\underline{H}^0_\X(\Omega_\Y )) \\
   &=& i_*\Rhom(\Rhom(A,\Omega_\X),\Omega_X)&(\ref{dual-imm-ferme} ) \\
   &=&i_*A &(\ref{bidualnonborne})
\end{array}$$
\hfill\hfill$\square$\end{proof}

\subsection{Computation of $f^!$ for a universal homeomorphism} By
universal homeomorphism we mean a representable, radiciel and
surjective morphism. By Zariski's main theorem, such a morphism is finite.

In the schematic situation, we know that such a morphism induces
an isomorphism of the \'etale topos (\cite{SGA42}, VIII.1.1). In
particular, $f^*$ is also a  right adjoint of $f_*$. Being exact,
one gets in this case an identification $f^*=f^!$. In particular,
$f^*$ identifies the corresponding dualizing complexes. Exactly as
in the proof of~\ref{dual-imm-ferme}, one gets

\begin{lemma}\label{dual-uni-homeo} Let $f:\X\ra\Y$ be a universal homeomorphism of stacks. One has $f^*\Omega_\X=\Omega_\Y$.\end{lemma}

One gets therefore
\begin{corollary}\label{cor-uni-homeo} Let $f:\X\ra\Y$ be a universal homeomorphism of
stacks. One has one has $f^!=f^*$ and $Rf_!=Rf_*$.\end{corollary}
\begin{proof} One has
$$\begin{array}{ccll}
  f^!A&=&\Rhom(f^*\Rhom(A,\Omega_\Y),\Omega_\X) \\
  &=& \Rhom(\Rhom(f^*A,f^*\Omega_\Y),\Omega_\X)&(\ref{Rhom-im-inv}) \\
   &=& \Rhom(\Rhom(f^*A,\Omega_\X),\Omega_\X)&(\ref{dual-uni-homeo} ) \\
   &=&f^*A &(\ref{bidualnonborne}).
\end{array}$$
The last formula follows by adjunction.
\hfill\hfill$\square$\end{proof}

\subsection{Computation of $Rf_!$ via hypercovers}

Let $Y$ be an $S$--scheme of finite type and $f:\X\rightarrow Y$ a morphism of finite type from an algebraic stack $\X$.  Let $X_\bullet \rightarrow \X$ be a smooth hypercover by  algebraic spaces, and for each $n$ let $d_n$ denote the locally constant function on $X_n$ which is the relative dimension over $\X$.  By the construction, the restriction of the dualizing complex $\Omega _{\X}$ to each $X_{n, \et}$ is canonically isomorphic to the dualizing complex $K_{X_n}= \Omega _{\X_n}\langle -d_n\rangle$ of $X_n$. Let $K _{X_\bullet }$ denote the restriction of $\Omega _{\X}$ to $X_{\bullet , \et}$.

Let $L\in D_c^-(\X)$, and let $L|_{X_\bullet }$ denote the restriction of $L$ to $X_{\bullet , \et}$.   Then $D_{\X}(L)|_{X_\bullet }$ is isomorphic to $D_{X_\bullet }(L|_{X_\bullet }):= \Rhom_{X_{\bullet , \et}}(L|_{X_\bullet }, K_{X_\bullet })$. In particular,  the restriction of $Rf_!L$ to $Y_{\et}$ is canonically isomorphic to
\begin{equation}
\Rhom _{Y_{\et}}(Rf_{\bullet *}D_{X_\bullet }(L|_{X_\bullet }), K_{Y})\in D_c(Y_{\et}),
\end{equation}
where $f_{\bullet }:X_{\et}\rightarrow Y_{\et }$ denotes the morphism of topos induced by $f$.

Let $Y_{\bullet , \et}$ denote the simplicial topos obtained by viewing $Y$ as a constant simplicial scheme.  Let $\epsilon :Y_{\bullet, \et}\rightarrow Y_{\et}$ denote the canonical morphism of topos, and let $\tilde f:X_{\bullet , \et}\rightarrow Y_{\bullet , \et}$ be the morphism of topos induced by $f$.  We have $f_\bullet = \epsilon \circ \tilde f$.  As in \cite{Ols05}, 2.7, it follows that there is a canonical spectral sequence
\begin{equation}
E_1^{pq} = R^qf_{p*}D_{X_p}(L|_{X_p})\implies R^{p+q}f_{\bullet *}D_{X_\bullet }(L_{X_\bullet }).
\end{equation}
On the other hand, we have
$$
R^qf_{p*}D_{X_p}(L|_{X_p}) = R^qf_{p*}\Rhom (L|_{X_p}, \Omega _{X_p}\langle -d_p\rangle ) \simeq \mc H^q(D_Y(Rf_{p!}(L|_{X_p}\langle d_p\rangle )),
$$
where the second isomorphism is by biduality \ref{Dinvolutif}. Combining all this we obtain

\begin{proposition}\label{5.16}
There is a canonical spectral sequence
\begin{equation}\label{5.16.1}
E_1^{pq} = \mc H^q(D_{Y_{\text{\rm et}}}(Rf_{p!}L|_{X_p}\langle d_p\rangle ))\implies \mc H^{p+q}(D_{Y_{\text{\rm et}}}(Rf_!L|_{Y_{\text{\rm et}}})).
\end{equation}
\end{proposition}

\begin{ex} Let $k$ be an algebraically closed field and $G$ a finite group.  We can then compute $H^*_c(BG, \Lambda )$ as follows.  We first compute $\Rhom (R\Gamma _!(BG, \Lambda ), \Lambda )$.  Let $\text{Spec}(k)\rightarrow BG$ be the surjection corresponding to the trivial $G$--torsor, and let $X_\bullet \rightarrow BG$ be the $0$--coskeleton.  Note that each $X_n$ isomorphic to $G^n$ and in particular is a discrete collection of points.  Therefore $Rf_{!p}\Lambda \simeq \text{Hom}(G^n, \Lambda )$. From this it follows that $\Rhom (R\Gamma _!(BG, \Lambda ), \Lambda )$ is represented by the standard cochain complex computing the group cohomology of $\Lambda $, and hence $R\Gamma _!(BG, \Lambda )$ is the dual of this complex.  In particular, this can be nonzero in infinitely many negative degrees. For example if $G = \mathbb{Z}/\ell $ for some prime $\ell $ and $\Lambda = \mathbb{Z}/\ell $ since in this case the group cohomology $H^i(G, \mathbb{Z}/\ell )\simeq \mathbb{Z}/\ell  $ for all $i\geq 0$.
\end{ex}

\begin{ex} Let $k$ be an algebraically closed field and $P$ the affine line $\mathbb{A}^1$ with the origin doubled. By definition $P$ is equal to two copies of $\mathbb{A}^1$ glued along $\mathbb{G}_m$ via the standard inclusions $\mathbb{G}_m\subset \mathbb{A}^1$.   We can then compute $R\Gamma _!(P, \Lambda )$ as follows.  Let $j_i:\mathbb{A}^1\hookrightarrow P$ ($i=1,2$) be the two open immersions, and let $h:\mathbb{G}_m\hookrightarrow P$ be the inclusion of the overlaps.  We then have an exact sequence
$$
0\rightarrow h_!\Lambda \rightarrow j_{1!}\Lambda \oplus j_{2!}\Lambda \rightarrow \Lambda \rightarrow 0.
$$
From this we obtain a long exact sequence
$$
\cdots \rightarrow H^i_c(\mathbb{G}_m, \Lambda )\rightarrow H^i_c(\mathbb{A}^1, \Lambda )\oplus H^i_c(\mathbb{A}^1, \Lambda )\rightarrow H^i_c(P, \Lambda )\rightarrow \cdots.
$$
From this sequence one deduces that $H^0_c(P, \Lambda )\simeq \Lambda $, $H^2_c(P, \Lambda )\simeq \Lambda (1)$, and all other cohomology groups vanish.  In particular, the cohomology of $P$ is isomorphic to the cohomology of $\mathbb{P}^1$.
\end{ex}

\subsection{Purity and the fundamental distinguished triangle}

We consider the usual situation of a closed immersion $i:\X\ra\Y$
of  stacks, the open immersion of the complement of $\Y$ being
$j:\U=\Y-\X\ra\Y$. For any (complex) of sheaves $A$ on $\Y$, one
has the exact sequence
$$0\ra j_{!}j^*A\ra A \ra i_*i^*A\ra 0.$$
Therefore, for any $A\in D_c(\Y)$, one has  the distinguished
triangle (\ref{!=!!})
\begin{equation}\label{T1}
    j_{!}j^*A\ra A \ra i_*i^*A
\end{equation}
which by duality gives the distinguished triangle

\begin{equation}\label{T2}
    i_{*}i^!A\ra A \ra j_*j^*A.
\end{equation}

Recall (\ref{prop-dual-imm-ferme}) the formula
$i^!=R\underline{H}^0_\X$. The usual purity theorem for
$S$-schemes gives

\begin{proposition}[Purity]\label{purity}
Assume moreover that $i$ is a closed immersion of smooth
$S$-stacks of codimension $c$ (a locally constant function on
$\Y$). Then, one has $i^!A=i^*A(-c)[-2c]$.
\end{proposition}

\begin{proof} Let $d$ be the relative dimension of $\Y\ra S$ and $s$ the
dimension of $S$. The relative dimension of $\X$ is (the
restriction to $\X$ of) $d-c$. By~\ref{f!lisse}, one has
$$\Omega_\Y=\Lambda(d+s)[2d+2s]\textup{ and
}\Omega_\X=\Lambda(d-c+s)[2d-2c-2e].$$ The identity
$i^!\Omega_\Y=\Omega_\X$ gives therefore the formula
\begin{equation}\label{i!c}
i^!\Lambda=\Lambda(-c)[-2c].
\end{equation}
By~\ref{dual}, one has
$$i_*\Rhom(i^!\Lambda,i^!A)=\Rhom(i_!i^!\Lambda,A)$$
which by adjunction for $i_*$ gives a map
$$i^*\Rhom(i_!i^!\Lambda,A)\ra \Rhom(i^!\Lambda,i^!A).$$
But the adjunction map (for $i_!$) $i_!i^!\Lambda\ra \Lambda$
dualizes to
$$\Rhom(\Lambda,A)\ra\Rhom(i_!i^!\Lambda,A)$$ which gives by composition a
morphism $$i^*A=i^*\Rhom(\Lambda,A)\ra
\Rhom(i^!\Lambda,i^!A)=i^!A(c)[2c]$$ which is the usual morphism
for closed immersion of schemes. This morphism is compatible with
the duality in an obvious sense. The usual purity theorem gives
then the proposition, at least for $A\in D^+_c(\Y)$. By duality,
one gets the proposition for $A\in D^-_c(\Y)$, and therefore for
$A\in D_c(\Y)$ using the distinguished triangle $\tau_{>0}A\ra
A\ra\tau_{\leq 0}A$.\hfill\hfill$\square$\end{proof}

\section{Base change} We start with a cartesian diagram of
 stacks
 $$\xymatrix{\X'\ar[r]^\pi\ar[d]_\phi\ar@{}[rd]|\Box&\X\ar[d]^f\\\Y'\ar[r]^p&\Y}$$
and we would like to prove a natural base change isomorphism
\begin{equation}\label{bc!}
  p^*Rf_!=R\phi_!\pi^*
\end{equation} of functors $D_c(\X)\ra D_c(\Y')$.
Though technically not needed, before proving the general base change Theorem we consider first some simpler cases
where one can prove a dual version:
\begin{equation}\label{bc}
    p^!Rf_*=R\phi_*\pi^!.
\end{equation}
 \subsection{Smooth base change}\label{5.1}
In this subsection we prove the base change isomorphism in the case when $p$ (and hence also $\pi $) is smooth.

\begin{proof}Because the relative dimension of
 $p$ and $\pi$ are the same, by~\ref{f!lisse}, one reduces the formula~\ref{bc} to
$$p^*Rf_*=R\phi_*\pi^*.$$ By adjunction, one has a morphism
$p^*Rf_*\ra R\phi_*\pi^*$ which we claim is an isomorphism (for
complexes bounded below this follows immediately from the smooth
base change theorem).  To prove that this map is an isomorphism,
we consider first the case when $\mc Y'$ is algebraic space and
show that our morphism restricts to an isomorphism on $\mc Y'_{\et
}$.  Since $p$ is representable $\mc X'$ represents a sheaf on
$\X_{\le }$.  Let $\X_{\le |\X'}$ denote the localized topos,
$w:\X_{\le |\X'}\rightarrow \Y_\et $ the projection,  and let
$A\in \D_c(\X)$ be a complex.   Let $X\rightarrow \X$ be a smooth
surjection with $X$ a scheme, and let $X_\bullet \rightarrow \X$
denote the associated simplicial space.  Let $X'_\bullet $ denote
the base change of $X_\bullet $ to $\mc Y'$.  Then $X_\bullet '$
defines a hypercover of the initial object in the topos $\X_{\le
|\X'}$ and hence we have an equivalence of topos $\X _{\le ,
\X'}\simeq \X _{\le , X_{\bullet }'}$.  Let $w_\bullet :\X_{\le
|X_{\bullet '}}\rightarrow \mc Y'_\et $ be the projection. Since
the restriction functor from $\X_{\le }$ to $\X_{\le |X_\bullet
'}$ takes homotopically injective complexes  to homotopically
injective complexes (since it has an exact left adjoint),
$p^*Rf_*A|_{\mc Y'_\et }$ is equal to $Rw_{\bullet *}(A|_{\X _{\le
|X_\bullet '}})$. On the other hand, $w_\bullet $ factors as
\begin{equation}
\begin{CD}
\mc X_{\le |X'_{\bullet }}@>\alpha >> X'_{\bullet , \et } @>\phi _\bullet >> \mc Y'_{\et },
\end{CD}
\end{equation}
where $\phi _{\bullet }:X'_{\bullet }\rightarrow \mc Y'_\et $ is the projection. Since $\alpha _*$ is exact, we find that $Rw_{\bullet *}(A|_{\X _{\le |X_\bullet '}})$ is isomorphic to $R\phi _{\bullet *}(A|_{X_{\bullet , \et }'}).$ Similarly, factoring the morphism of topos $\mc X_{\le }'\simeq \mc X'_{\le |X'_{\bullet }}\rightarrow \mc Y'_\et $ as
\begin{equation}
\begin{CD}
\mc X'_{\le |X'_\bullet }@>>> X'_{\bullet, \et }@>\phi _\bullet >> \mc Y'_\et
\end{CD}
\end{equation}
we see that $R\phi_*(A|_{\mc X'_\le })$ is isomorphic to $R\phi
_{\bullet *}(A|_{X'_{\bullet \et }})$. We leave to the reader that
the resulting isomorphism $p^*Rf_*(A)|_{\mc Y'_\et }\rightarrow
R\phi_*\pi ^*A|_{\mc Y'_\et }$ agrees with the morphism defined
above.  Thus this proves the case when $p$ is representable.

For the general case, let $Y'\rightarrow \mc Y'$ denote a smooth surjection with $Y'$ a scheme, so we have a commutative diagram
\begin{equation}
\begin{CD}
X'@>\sigma >> \X'@>\pi >> \X\\
@VgVV @VV\phi V @VVfV \\
Y'@>q>> \Y'@>p>> \Y
\end{CD}
\end{equation}
with cartesian squares. Let $A\in D_c(\X)$.  To prove that the
morphism $p^*Rf_*A\rightarrow R\phi_*\pi ^*A$ is an isomorphism,
it suffices to show that the induced morphism
$q^*p^*Rf_*A\rightarrow q^*R\phi_*\pi ^*A$ is an isomorphism.  By
the representable case, $q^*R\phi_*\pi ^*A\simeq Rg_*\sigma ^*\pi
^*A$ so it suffices to prove that the composite
\begin{equation}
q^*p^*Rf_*A\rightarrow Rg_*\sigma ^*\pi ^*A
\end{equation}
is an isomorphism.  By the construction, this map is equal to the base change morphism for the diagram
\begin{equation*}
\begin{CD}
X'@>\pi \circ \sigma >> \X\\
@VgVV @VVfV \\
Y'@>p\circ q>> \mc Y
\end{CD}
\end{equation*}
and hence it is an isomorphism by the representable case.
\hfill\hfill$\square$\end{proof}

\subsection{Computation of $Rf_*$ for proper representable
morphisms}

\begin{prop}\label{f!=f*propre} Let $f:\X\rightarrow \Y$ be a proper representable morphism of  $S$--stacks. Then the functor $Rf_!:D_c^-(\X)\rightarrow D_c^-(\Y)$ is canonically isomorphic to $Rf_*:D_c^-(\X)\rightarrow D_c^-(\Y)$.
\end{prop}
\begin{proof}
The key point is the following lemma.
\begin{lem} There is a canonical morphism $Rf_*\Omega _\X\rightarrow \Omega _\Y$.
\end{lem}
\begin{proof}
Using~\ref{D1} and smooth base change, it suffices to construct a
functorial morphism in the case of schemes, and to show that $\ext
^i(Rf_*\Omega _{\X}, \Omega _\Y) = 0$ for $i<0$. Now if $\X$ and
$\Y$ are schemes, we have $\Omega _\X = f^!\Omega _{\Y}$ so we
obtain by adjunction and the fact that $Rf_! = Rf_*$ a morphism
$Rf_*\Omega _\X\rightarrow \Omega _\Y$. For the computation of
$\ext$'s note that
$$\Rhom(Rf_*\Omega_\X,\Omega_\Y)=\Rhom(\Omega_\X,f^!\Omega_\Y)=\Lambda.$$
\hfill\hfill$\square$\end{proof}

We define a  map $Rf_*\circ D_\X\ra D_\Y\circ f_*$ by taking the
composite
$$Rf_*\Rhom(-,\Omega_\X)\ra \Rhom(Rf_*(-),Rf_*\Omega_\X)\ra \Rhom (Rf_*(-), \Omega _\Y).$$
To verify that this map is an isomorphism we may work locally on
$\Y$.  This reduces the proof to the case when $\X$ and $\Y$ are
algebraic spaces in which case the result is standard.
\hfill\hfill$\square$\end{proof}

 \subsection{Base change by an immersion}\label{5.3}
In this subsection we consider the case when $p$ is an immersion.

By replacing $\Y$ by a suitable open substack, one is reduced to the case when $p$ is a closed immersion. Then, ~\ref{bc!} follows from the projection formula~\ref{corkunneth} as in
~\cite{Fre-Kie}, p.81. Let us recall the argument. Let $A\in D_c(\X)$. Because
$p$ is a closed immersion, one has $p^*p_*=\Id$.
One has
(projection formula~\ref{corkunneth} for $p$)
$$p_*p^*Rf_!A=p_*\Lambda\Otimes Rf_!A.$$
One has then $$ Rf_!A\Otimes p_*\Lambda=Rf_!(A\Otimes
f^*p_*\Lambda)$$ (projection formula~\ref{corkunneth} for $f$).
But, we have trivially the base change for $p$, namely
$$f^*p_*=\pi_*\phi^*.$$
Therefore, one gets\begin{eqnarray*}
  Rf_!(A\Otimes
f^*p_*\Lambda) &=& Rf_!(A\Otimes \pi_*\phi^*\Lambda) \\
  &=& Rf_!\pi_*(\pi^*A\Otimes \phi^*\Lambda)\textup{ projection for }\pi \\
  &=& p_*\phi_!\pi^*A\textup{ because }\pi_*=\pi_!\
  (\ref{f!=f*propre}).
\end{eqnarray*} Applying $p^*$ gives the base change isomorphism.

\begin{remark}
 One can prove~\ref{bc}, {at least for $A$ bounded below}, more directly
as follows. Start with $A$ on $\X$ an injective complex. Because
$R^0f_*A_i$ is flasque, it is $\Gamma_{\Y'}$-acyclic. Then,
$p^!Rf_*A$ can be computed using the
 complex  ${\mc H}^0_{\Y'}(R^0f_*A_i)$. On the other hand, $\pi^!A$ can
 be computed by the complex ${\mc H}^0_{\X'}(A_i)$ which is a flasque
 complex  (formal, or
\cite{SGA42}, V.4.11). Therefore, the direct image by $\phi$ is
just $R^0\phi_*{\mc H}^0_{\X'}(A_i)$. One is reduced to the
formula
$$R^0\phi_*{\mc H}^0_{\X'}={\mc H}^0_{\Y'}(R^0f_*).$$
\end{remark}

\subsection{Base change by a universal homeomorphism}\label{5.4} If $p$ is a universal homeomorphism,
then $p^! = p^*$ and $\pi ^! = \pi ^*$. Thus in this case \ref{bc}
is equivalent to an isomorphism $p^*Rf_*\rightarrow R\phi _*\pi
^*$.  We define such a morphism by taking the usual base change
morphism (adjunction).

Let $A\in \D _c(\mc X)$.  Using a hypercover of $\X$ as in
\ref{5.1}, one sees that to prove that the map
$p^*Rf_*A\rightarrow R\phi _*\pi ^*A$ is an isomorphism it
suffices to consider the case when $\X$ is a scheme. Furthermore,
by the smooth base change formula already shown, it suffices to
prove that this map is an isomorphism after making a smooth base
change $Y\rightarrow \Y$.  We may therefore assume that $\Y$ is
also a scheme in which case the result follows from the classical
corresponding result for \'etale topology (see~\cite{SGA1},
IV.4.10).

\subsection{Base change morphism in general}\label{5.5}

Before defining the base change morphism we need a general
construction of strictly simplicial schemes and algebraic spaces.

Fix an algebraic stack $\mc X$.  In the following construction all
 schemes and morphisms are assumed over $\mc  X$ (so in
particular products are taken over $\mc  X$).

Let $X_\bullet $ be a strictly simplicial scheme, $[n]\in \Delta ^+
$ an object, and $a:V\rightarrow X_n$ a surjective morphism. We
then construct a strictly simplicial scheme $M(X_\bullet , a)$ (sometimes
written $M_{\mc  X}(X_\bullet , a)$ if we want to make clear the reference
to $\mc  X$) with a morphism $M(X_\bullet , a)\rightarrow X_\bullet $ such that the
following hold:
\begin{enumerate}
\item [(i)] For $i<n$ the morphism $M(X_\bullet , a)_i\rightarrow X_i$ is an
isomorphism.
\item [(ii)] $M(X_\bullet , a)_n$ is equal to $V$ with the projection to
$X_n$ given by $a$.
\end{enumerate}

The construction of $M(X_\bullet , a)$ is a standard application of the
skeleton and coskeleton functors (\cite{SGA42}, exp. Vbis). Let us
review some of this because the standard references deal only with
simplicial spaces whereas we consider strictly simplicial spaces.

To construct $M(X_\bullet , a)$, let $\Delta ^+ _n\subset \Delta ^+ $
denote the full subcategory whose objects are the finite sets with
cardinality $\leq n$.  Denote by $\text{Sch}^{\Delta ^{+\opp} _n}$
the category of functors from $\Delta ^{+\opp} _n$ to schemes (so
$\text{Sch}^{\Delta ^{+\opp}}$ is the category of strictly
simplicial schemes). Restriction from $\Delta ^{+\opp}$ to $\Delta
^{+\opp}_n$ defines a functor (the \emph{$n$-skeleton functor})
\begin{equation}
\text{sq}_n:\text{Sch}^{\Delta ^{+\opp}}\rightarrow
\text{Sch}^{\Delta ^{+\opp}_n}
\end{equation}
which has a right adjoint
\begin{equation}
\text{cosq}_n:\text{Sch}^{\Delta ^{+\opp}_n}\rightarrow
\text{Sch}^{\Delta ^{+\opp}}
\end{equation}
called the \emph{$n$-th coskeleton functor}. For $X_\bullet \in
\text{Sch}^{\Delta ^{+\opp}_n}$, the coskeleton $\text{cosq}_nX$
in degree $i$ is equal to
\begin{equation}\label{cosqdef}
(\text{cosq}_nX)_i = \varprojlim _{\substack{[k]\rightarrow [i]\\
k\leq n}}X_k,
\end{equation}
where the limit is taken over the category of morphisms
$[k]\rightarrow [i]$ in $\Delta ^+$ with $k\leq n$.

 Note in
particular that for $i\leq n$ we have $(\text{cosq}_nX)_i = X_i$ since the category of morphisms $[k]\rightarrow [i]$ has an initial object $\text{id}:[i]\rightarrow [i]$.

\begin{lemma}\label{easylem} For any $X_\bullet \in \text{\rm Sch}^{\Delta ^{+\opp}_n}$ and
$i>n$ the morphism
\begin{equation}
(\text{\rm cosq}_nX)_i\rightarrow (\text{\rm cosq}_{i-1}\text{\rm
sq}_{i-1}\text{\rm cosq}_nX)_i
\end{equation}
is an isomorphism.
\end{lemma}
\begin{proof} Using the formula \ref{cosqdef} the morphism can
be identified with the natural map
\begin{equation}
\varprojlim _{\substack{[k]\rightarrow [i]\\ k\leq
n}}X_k\rightarrow \varprojlim _{\substack{[k]\rightarrow [i]\\
k\leq i-1}}(\varprojlim _{\substack{[w]\rightarrow [k]\\ w\leq
n}}X_w)
\end{equation}
which is clearly an isomorphism. \hfill\hfill$\square$\end{proof}

\begin{lemma} The functors $\text{\rm sq}_n$ and $\text{\rm cosq}_n$
commute with fiber products.
\end{lemma}
\begin{proof} The functor $\text{sq}_n$ commutes with fiber
products by construction, and the functor $\text{cosq}_n$ commutes
with fiber products by adjunction.
\hfill\hfill$\square$\end{proof}

 To construct $M(X_\bullet , a)$, we first construct an object
$M'(X_\bullet , a)\in \text{Sch}^{\Delta ^{+\opp}_n}$.  The restriction of
$M'(X_\bullet , a)$ to $\Delta ^{+\opp} _{n-1}$ will be equal to
$\text{sq}_{n-1}X$, and $M'(X_\bullet , a)_n$ is defined to be $V$.  For
$0\leq j\leq n$ define $\delta _j:M'(X_\bullet , a)_n\rightarrow M'(X_\bullet ,
a)_{n-1} = X_{n-1}$ to be the composite
\begin{equation}
\begin{CD}
V@>a>> X_n@>\delta _{X, j}>> X_{n-1},
\end{CD}
\end{equation}
where $\delta _{j, X}$ denotes the map obtained from the strictly
simplicial structure on $X_\bullet $. There is an obvious morphism
$$M'(X_\bullet ,a)\ra\text{sq}_n(X_\bullet )\text{ inducing }\text{cosq}_nM'(X_\bullet ,
a))\ra{\text{cosq}_n\text{sq}_nX_\bullet }.$$

We then define
\begin{equation}
M(X_\bullet , a):= (\text{cosq}_nM'(X_\bullet , a))\times
_{\text{cosq}_n\text{sq}_nX_\bullet }X_\bullet ,
\end{equation}
where the map $X_\bullet \rightarrow \text{cosq}_n\text{sq}_nX_\bullet $ is the
adjunction morphism. The map $M(X_\bullet , a)\rightarrow X_\bullet $ is defined to
be the projection. The properties  (i) and (ii)
follow immediately from the construction.

\begin{proposition} Let $\mc  X$ be an algebraic stack and $X_\bullet
\rightarrow \mc  X$ a hypercover by schemes.  Let $n$ be a
natural number and $a:V\rightarrow X_n$ a surjection.  Then
$M_{\mc  X} (X_\bullet , a)\rightarrow \mc  X$ is also a hypercover. If $X_\bullet $ is a smooth hypercover and $a$ is smooth
and surjective, then  $M_{\mc  X}(X_\bullet , a)$ is also a smooth hypercover.
\end{proposition}
\begin{proof}
By definition of a  hypercover, we must verify that for all $i$
the map
\begin{equation}
M(X_\bullet , a)_i\rightarrow (\text{cosq}_{i-1}\text{sq}_{i-1}M(X_\bullet , a))_i
\end{equation}
is surjective.  Note that this is immediate for $i\leq n$.  For
$i>n$ we compute
\begin{eqnarray*}
(\text{cosq}_{i-1}\text{sq}_{i-1}M(X_\bullet , a))_i &\simeq &
(\text{cosq}_{i-1}\text{sq}_{i-1}(\text{cosq}_nM'(X_\bullet , a)\times
_{\text{cosq}_n\text{sq}_nX_\bullet }X_\bullet ))_i\\
&\simeq &  (\text{cosq}_{i-1}\text{sq}_{i-1}(\text{cosq}_nM'(X_\bullet ,
a)))_i\times _{(\text{cosq}_i\text{sq}_{i-1}\text{cosq}_n\text{sq}_nX_\bullet )_i}(\text{cosq}_{i-1}\text{sq}_{i-1}X_\bullet )_i\\
&\simeq &  (\text{cosq}_nM'(X_\bullet , a))_i\times
_{(\text{cosq}_n\text{sq}_nX_\bullet )_i}(\text{cosq}_{i-1}\text{sq}_{i-1}X_\bullet )_i.
\end{eqnarray*}
Here the second isomorphism is because $\text{sq}_n$ and
$\text{cosq}_n$ commute with products, and the third isomorphism is
by \ref{easylem}. Hence it suffices to show that the natural map
\begin{equation}
X_i\rightarrow (\text{cosq}_{i-1}\text{sq}_{i-1}X_\bullet )_i
\end{equation}
is surjective, which is true since $X_\bullet $ is a hypercover.
This also proves that if $X_\bullet $ is a smooth hypercover and
$a$ is smooth, then $M_{\mc  X}(X_\bullet , a)$ is a smooth
hypercover. \hfill\hfill$\square$\end{proof}

 The construction of $M_{\mc  X}(X_\bullet , a)$ is functorial.
Precisely, let $f:\mc  X\rightarrow \mc  Y$ be a morphism of
algebraic stacks, $X_\bullet \rightarrow \mc  X$ a strictly
simplicial scheme over $\mc  X$, $Y_\bullet \rightarrow \mc  Y$ a
strictly simplicial scheme over $\mc  Y$, and $f_\bullet :X_\bullet
\rightarrow Y_\bullet $ a morphism over $f$.  Then for any
commutative diagram of schemes
\begin{equation}
\begin{CD}
V@>\tilde f>> W\\
@VaVV @VVbV \\
X_n@>f_n>> Y_n
\end{CD}
\end{equation}
there is an induced morphism of strictly simplicial schemes
$M_{\mc  X}(X_\bullet , a)\rightarrow M_{\mc  Y}(Y_\bullet , b)$ over $f_\bullet $.

\begin{proposition}\label{5.4.4} Let $f:\mc  X\rightarrow \mc  Y$ be a morphism of
finite type between algebraic $S$-stacks locally of finite type.
Then there exists smooth hypercovers $p:X_\bullet \rightarrow \mc X$
and $q:Y_\bullet \rightarrow \mc  Y$ by schemes and a commutative
diagram
\begin{equation}\label{simpdiagram}
\begin{CD}
X_\bullet @>f_\bullet >> Y_\bullet \\
@VpVV @VVqV \\
\mc  X@>f>> \mc  Y,
\end{CD}
\end{equation}
where each morphism $f_n:X_n\rightarrow Y_n$ is a closed
immersion.
\end{proposition}
\begin{proof}
We construct inductively hypercovers $X_\bullet ^{(n)}\rightarrow
\mc  X$ and $Y_\bullet ^{(n)}\rightarrow \mc  Y$ and a commutative
diagram
\begin{equation}
\begin{CD}
X_\bullet ^{(n)}@>>> Y_\bullet ^{(n)}\\
@VVV @VVV \\
\mc  X@>>> \mc  Y
\end{CD}
\end{equation}
together with  a commutative diagram
\begin{equation}
\begin{CD}
X_\bullet ^{(n)}@>>> Y_\bullet ^{(n)}\\
@VVV @VVV \\
X_\bullet ^{(n-1)}@>>> Y_\bullet ^{(n-1)}
\end{CD}
\end{equation}
over $f$.  We further arrange so that the following hold:
\begin{enumerate}
\item [(i)] For $i<n$ the maps $X^{(n)}_i\rightarrow X^{(n-1)}_i$
and $Y^{(n)}_i\rightarrow Y^{(n-1)}_i$ are isomorphisms.
\item [(ii)] For $i\leq n$ the maps $X^{(n)}_i\rightarrow
Y^{(n)}_i$ are closed immersions.
\end{enumerate}
This suffices for we can then take $X_\bullet = \varprojlim X^{(n)}$ and $Y_\bullet  = \varprojlim Y^{(n)}$.

For the base case $n=0$, choose any $2$--commutative diagram
\begin{equation}
\begin{CD}
X=\sqcup X_i@>\tilde f=\sqcup\tilde f_i>> Y=\sqcup Y_i\\
@V{p=\sqcup p_i}VV @VV{q=\sqcup q_i}V \\
\mc  X@>>> \mc  Y,
\end{CD}
\end{equation}
with $p_i$ and $q_i$ smooth, surjective, and of finite type, and
$X_i$ and $Y_i$ affine schemes. Then $\tilde f_i$ are also of
finite type, so there exists a closed immersion
$X_i\hookrightarrow \mathbb{A}^{r_i}_{Y_i}$ for some integer $r$
over $X_i\rightarrow Y_i$. Replacing $Y_i$ by
$\mathbb{A}^{r_i}_{Y_i}$ we may assume that $\tilde f$ is a closed
immersions.  We then obtain $X_\bullet ^{(0)}\rightarrow Y_\bullet
^{(0)}$ by taking the coskeletons of $p$ and $q$.

Now assume that $X_\bullet ^{(n-1)}\rightarrow Y_\bullet ^{(n-1)}$ has been
constructed.  Choose a commutative diagram
\begin{equation}
\begin{CD}
V@>j>> W\\
@VaVV @VVbV \\
X_n^{(n-1)}@>>> Y_n^{(n-1)},
\end{CD}
\end{equation}
with $a$ and $b$ smooth and surjective, and $j$ a closed
immersion. Then define $X_\bullet ^{(n)}\rightarrow Y_\bullet ^{(n)}$ to be
\begin{equation}
M_{\mc  X}(X^{(n-1)}_\bullet , a)\rightarrow M_{\mc  Y}(Y_\bullet ^{(n-1)}, b).
\end{equation}
\hfill\hfill$\square$\end{proof}

\begin{rem}\label{filteringremark}
The same argument used in the proof shows that for any commutative
diagram
\begin{equation}
\begin{CD}
X_{\bullet }@>f_\bullet >> Y_\bullet \\
@VpVV @VVqV \\
\mc X@>f>> \mc Y,
\end{CD}
\end{equation}
where $p$ and $q$ are smooth hypercovers, there exists a morphism
of simplicial schemes $g:\widetilde X_\bullet \rightarrow \widetilde
Y_\bullet $ over $f_\bullet $ with each $g_n:\widetilde X_n\rightarrow
\widetilde Y_n$ an immersion such that $\widetilde X_\bullet $
(resp. $\widetilde Y_\bullet $) is a hypercover of $\mc X$ (resp.
$\mc Y$).  In other words, the category of diagrams
\ref{simpdiagram} is connected.
\end{rem}

Let $f:\mc  X\rightarrow \mc  Y$ be a morphism of  algebraic stacks
over $S$.  For $F\in \D_c^-(\mc  X)$ we can compute $Rf_!F$ as
follows. Let $Y_\bullet \rightarrow \mc  Y$ be a smooth hypercover,
and let $\pi :\mc  X_{Y_\bullet }\rightarrow \mc  X$ be the base
change of $\mc  X$ to $Y_\bullet $.  Let $f_\bullet :\mc X_{Y_\bullet
}\rightarrow Y_\bullet $ be the projection. Let $\omega _{\mc
X_{Y_\bullet }}$ denote the pullback of the dualizing sheaf $\Omega
_{\mc  X}$ to $\mc  X_{Y_\bullet }$, and let $D_{\mc X_{Y_\bullet }}$
denote the functor $\Rhom(-, \omega _{\mc X_{Y_\bullet }})$.
Similarly let $\omega _{Y_\bullet }$ denote the pullback of $\Omega
_{\mc  Y}$ to $Y_\bullet $, and let $D_{Y_\bullet }$ denote $\Rhom  (-, \omega _{Y_\bullet })$.

If $d_n$ (resp. $d_n'$) denotes the relative dimension of $Y_n$
over $\Y$ (resp. $Y_n'$ over $\Y'$), then $d_n$ (resp. $d_n'$) is
also equal to the relative dimension of $\X_{Y_n}$ over $\X$
(resp. $\X'_{Y_n'}$ over $\X'$).  From \ref{f!lisse} it follows that
the restriction of $\omega _{\X_{Y_\bullet }}$ to $\X_{Y_n}$ is
canonically isomorphic to $\Omega _{\X_{Y_n}}\langle -d_n\rangle
$. Similarly the restriction of $\omega _{Y_\bullet }$ to $Y_n$ is
canonically isomorphic to $\Omega _{Y_n}\langle -d_n\rangle $.
Note that this combined with \ref{Dinvolutif} shows that
$D_{Y_\bullet }\circ D_{Y_\bullet } = \text{id}$ (resp.
$D_{\X_{Y_\bullet }}\circ D_{\X_{Y_\bullet }} = \text{id}$) on the
category $\D_c(Y_\bullet )$ (resp. $\D_c(\X_{Y_\bullet })$).

For $F\in D_c(\X)$,
we can then consider
\begin{equation}
D_{Y_\bullet }Rf_{\bullet *}D_{\mc  X_{Y_\bullet }}(\pi ^*F)\in
\D(Y_{\bullet , \et}).
\end{equation}
The sheaf $D_{\mc  X_{Y_\bullet }}(\pi ^*F)$ is just the restriction
of $D_{\mc  X}(F)$ to $\mc  X_{Y_\bullet }$. It follows from this
that $Rf_{\bullet *}D_{\mc X_{Y_\bullet }}(\pi ^*F)$ is equal to the
restriction of $Rf_*D_{\mc  X}(F)$ to $Y_\bullet $, and this in turn
implies that $D_{Y_\bullet }Rf_{\bullet *}D_{\mc X_{Y_\bullet }}(\pi
^*F)$ is isomorphic to the restriction of $Rf_!F$ to $Y_{\bullet ,
\et}$. From this we conclude that $Rf_!F$ is equal to the sheaf
obtained from $D_{Y_\bullet }Rf_{\bullet *}D_{\mc  X_{Y_\bullet }}(\pi
^*F)$ and the equivalence of categories (\ref{etexample})
$\D_c(\mc  Y)\simeq \D_c(Y_\bullet ).$

\begin{theorem}\label{5.5.6} Let
\begin{equation}\label{5.5.6.1}
\begin{CD}
\mc  X'@>a>> \mc  X\\
@Vf'VV @VVfV \\
\mc  Y'@>b>> \mc  Y
\end{CD}
\end{equation}
be a cartesian square of  stacks over $S$. Then there is a
natural isomorphism of functors
\begin{equation}
b^*Rf_!\rightarrow Rf'_!a^*.
\end{equation}
\end{theorem}
\begin{proof}
By \ref{5.4.4}, there exists a commutative diagram
\begin{equation}\label{5.5.6.3}
\begin{CD}
Y_\bullet '@>j>> Y_\bullet \\
@VpVV @VVqV \\
\mc  Y'@>\rho >> \mc  Y,
\end{CD}
\end{equation}
where $p$ and $q$ are smooth hypercovers and $j$ is a closed
immersion.

Let $\mc  X'_{Y'_\bullet }$ denote the base change $\mc  X'\times
_{\mc  Y'}Y_{\bullet }'$ and $\mc  X_{Y_\bullet }$ the base change
$\mc  X\times _{\mc  Y}Y_\bullet $. Then there is a cartesian
diagram
\begin{equation}
\begin{CD}
\mc  X'_{Y'_\bullet }@>i>> \mc  X_{Y_\bullet }\\
@Vg'VV @VVgV \\
Y'_\bullet @>j>> Y_\bullet ,
\end{CD}
\end{equation}
where $i$ and $j$ are closed immersions.

As before let $\omega _{\X ' _{Y'_\bullet }}$ (resp. $\omega
_{\X_{Y_\bullet }}$, $\omega _{Y'_\bullet }$, $\omega _{Y_\bullet }$)
denote the pullback of $\Omega _{\X'}$ (resp. $\Omega _\X$,
$\Omega _{\Y'}$, $\Omega _{\Y}$) to $\X'_{Y'_\bullet }$ (resp.
$\X_{Y_\bullet }$, $Y_\bullet '$, $Y_\bullet $), and let
$D_{\X'_{Y'_\bullet }}$ (resp. $D_{\X_{Y_\bullet }}$, $D_{Y'_\bullet }$,
$D_{Y_\bullet }$) denote the functor $\Rhom (-, \omega
_{X'_{Y'_\bullet }})$ (resp. $\Rhom (-, \omega _{\X_{Y_\bullet }})$,
$\Rhom (-, \omega _{Y'_\bullet })$, $\Rhom (-, \omega _{Y_\bullet
})$).

\begin{lem}\label{5.5.7} Let $\T$ be a topos and $\Lambda $ a sheaf of rings in $\T$.  Then for any $A, B, C\in \D(\T, \Lambda )$ there is a canonical morphism
\begin{equation}
A{\Otimes}\Rhom (B, C)\rightarrow \Rhom (\Rhom (A, B), C).
\end{equation}
\end{lem}
\begin{proof}
We have
\begin{equation}\label{adeq}
\RHom (A{\Otimes}\Rhom (B, C), \Rhom (\Rhom (A, B), C)) \simeq
\RHom (A{\Otimes}\Rhom (B, C){\Otimes}\Rhom (A, B), C).
\end{equation}
Let
\begin{equation*}
a:A{\Otimes}\Rhom (A, B)\rightarrow B, \ \ b:B {\Otimes}\Rhom (B,
C)\rightarrow C
\end{equation*}
be the evaluation morphisms.  Then the morphism
\begin{equation*}
A\Otimes\Rhom (B, C)\Otimes\Rhom (A, B)\xrightarrow{a}
B\Otimes\Rhom (B,C)\xrightarrow{b}C
\end{equation*}
and the isomorphism \ref{adeq} give the Lemma.
\hfill\hfill$\square$\end{proof}

Let $\mc F$ denote the functor
\begin{equation}
D_{Y'_\bullet }j^*D_{Y_\bullet }g_*D_{\X_{Y_\bullet
}}i_*D_{\X'_{Y'_\bullet }}:\D_c(\X'_{Y'_\bullet })\rightarrow
\D(Y'_\bullet ).
\end{equation}

\begin{prop}\label{5.4.8}
There is an isomorphism of functors $\mc F\simeq Rg'_*$.
\end{prop}
\begin{proof}
Consider first the functor $\mc F':= \mc F\circ g^{\prime *}$ and
let $A\in \D_c^-(Y'_\bullet )$.  Then
\begin{eqnarray*}
\mc F'(A) & = & D_{Y'_\bullet }j^*D_{Y_\bullet }g_*\Rhom (i_*\Rhom (g^{\prime *}A, \omega _{\X '_{Y'_\bullet }}), \omega _{\X_{Y_\bullet }}) \s (\text{definition})\\
& \simeq & D_{Y'_\bullet }j^*D_{Y_\bullet }g_*\Rhom (i_*\Rhom (i^*i_*g^{\prime *}A, \omega _{\X'_{Y'_\bullet }}), \omega _{\X_{Y_\bullet }}) \s (i^*i_* = \text{id})\\
& \simeq & D_{Y'_\bullet }j^*D_{Y_\bullet }g_*\Rhom (\Rhom (i_*g^{\prime *}A, i_*\omega _{\X '_{Y'_\bullet }}), \omega _{\X_{Y_\bullet }}) \s (\text{adjunction for $(i^*, i_*)$})\\
& \leftarrow & D_{Y'_\bullet }j^*D_{Y_\bullet }g_*(i_*g^{\prime *}A{\Otimes}\Rhom (i_*\omega _{\X'_{Y'_\bullet }}, \omega _{\X_{Y_\bullet }}))\s (\ref{5.5.7})\\
& \simeq & D_{Y'_\bullet }j^*D_{Y_\bullet }g_*(g^*j_*A{\Otimes}\Rhom (i_*\omega _{\X'_{Y'_\bullet }}, \omega _{\X_{Y_\bullet }}))\s (i_*g'^* = g^*j_* \text{ by proper base change})\\
& \simeq &  D_{Y'_\bullet }j^*D_{Y_\bullet }(j_*A{\Otimes}Rg_*\Rhom (i_*\omega _{\X'_{Y'_\bullet }}, \omega _{\X_{Y_\bullet }})) \s (\text{projection formula})\\
& \simeq & \Rhom (j^*\Rhom (j_*A{\Otimes}Rg_*\Rhom (i_*\omega _{\X'_{Y'_\bullet }}, \omega _{\X_{Y_\bullet }}), \omega _{Y_\bullet }), \omega _{Y'_\bullet })\s (\text{definition})\\
& \simeq & j^*j_*\Rhom (j^*\Rhom (j_*A{\Otimes}Rg_*\Rhom (i_*\omega _{\X'_{Y'_\bullet }}, \omega _{\X_{Y_\bullet }}), \omega _{Y_\bullet }), \omega _{Y'_\bullet })\s (j^*j_* = \text{id})\\
& \simeq & j^*\Rhom (\Rhom (j_*A{\Otimes}Rg_*\Rhom (i_*\omega _{\X'_{Y'_\bullet }}, \omega _{\X_{Y_\bullet }}), \omega _{Y_\bullet }), j_*\omega _{Y'_\bullet })\s (\text{adjunction for $(j^*, j_*)$})\\
& \simeq & j^*\Rhom (\Rhom (j_*A, \Rhom (Rg_*\Rhom (i_*\omega _{\X'_{Y'_\bullet }}, \omega _{\X_{Y_\bullet }}), \omega _{Y_\bullet })), j_*\omega _{Y'_{\bullet }})\\
& \leftarrow & A{\Otimes}j^*\Rhom (\Rhom (Rg_*\Rhom (i_*\omega
_{\X'_{Y'_\bullet }}, \omega _{\X_{Y_\bullet }}), \omega _{Y_\bullet }),
j_*\omega _{Y'_{\bullet }}) \s (\text{\ref{5.5.7}}).
\end{eqnarray*}
The following Lemma therefore shows that there is a canonical
morphism $A\rightarrow \F'(A)$ functorial in $A$.
\begin{lem} For all $s\in \mathbb{Z}$ there is a canonical isomorphism
$$
\mc H^s(j^*\Rhom (\Rhom (Rg_*\Rhom (i_*\omega _{\X'_{Y'_\bullet }},
\omega _{\X_{Y_\bullet }}), \omega _{Y_\bullet }), j_*\omega
_{Y'_{\bullet }})\simeq R^sg'_*\Lambda .
$$
In particular, $\tau _{\leq 0}j^*\Rhom (\Rhom (Rg_*\Rhom
(i_*\omega _{\X'_{Y'_\bullet }}, \omega _{\X_{Y_\bullet }}), \omega
_{Y_\bullet }), j_*\omega _{Y'_{\bullet }})\simeq g'_*\Lambda $, so
the composite
$$
\Lambda \rightarrow g'_*\Lambda \simeq \tau _{\leq 0}j^*\Rhom
(\Rhom (Rg_*\Rhom (i_*\omega _{\X'_{Y'_\bullet }}, \omega
_{\X_{Y_\bullet }}), \omega _{Y_\bullet }), j_*\omega _{Y'_{\bullet }})
$$
induces a canonical morphism
$$
\Lambda \rightarrow j^*\Rhom (\Rhom (Rg_*\Rhom (i_*\omega
_{\X'_{Y'_\bullet }}, \omega _{\X_{Y_\bullet }}), \omega _{Y_\bullet }),
j_*\omega _{Y'_{\bullet }}).
$$
\end{lem}
\begin{proof}
It suffices to construct such a canonical isomorphism over each
$Y'_n$.  Let $d_n$ (resp. $d_n'$) denote the relative dimension of
$Y_n$ (resp. $Y_n'$) over $\Y$ (resp. $\Y'$).  Note that $d_n$
(resp. $d_n'$) is also equal to the relative dimension of
$\X_{Y_n}$ (resp. $\X'_{Y'_n}$) over $\X$ (resp. $\X'$).  As
mentioned above we therefore have
$$
\omega _{\X'_{Y'_n}}\simeq \Omega _{\X'_{Y'_n}}\langle
-d_n'\rangle, \ \ \omega _{\X_{Y_n}}\simeq \Omega
_{\X_{Y_n}}\langle -d_n\rangle , \ \ \omega _{Y_n}\simeq \Omega
_{Y_n}\langle -d_n\rangle , \ \ \omega _{Y'_n}\simeq \Omega
_{Y'_n}\langle -d_n'\rangle.
$$
From this and an elementary manipulation using the identity
$$\Rhom (A\langle n\rangle , B\langle m\rangle )\simeq \Rhom (A,
B)\langle m-n\rangle$$ we get
\begin{equation}\label{bigmess}
\begin{matrix}
&& j^*\Rhom (\Rhom (Rg_*\Rhom (i_*\omega _{\X'_{Y'_n }}, \omega _{\X_{Y_n }}), \omega _{Y_n }), j_*\omega _{Y'_{n }})\\
&\simeq & j^*\Rhom (\Rhom (Rg_*\Rhom (i_*\Omega _{\X'_{Y'_n }}\langle -d_n'\rangle , \Omega _{\X_{Y_n }}\langle -d_n\rangle ), \Omega _{Y_n }\langle -d_n\rangle ), j_*\Omega _{Y'_{n }}\langle -d_n'\rangle )\\
& \simeq &  j^*\Rhom (\Rhom (Rg_*\Rhom (i_*\Omega _{\X'_{Y'_n }},
\Omega _{\X_{Y_n }}), \Omega _{Y_n }), j_*\Omega _{Y'_{n }}).
\end{matrix}
\end{equation}
We then get
\begin{eqnarray*}
\Rhom (i_*\Omega _{\X'_{Y'_n}}, \Omega _{\X_{Y_n}})& \simeq & \Rhom (i_!\Omega _{\X'_{Y'_n}}, \Omega _{\X_{Y_n}})\s (\text{\ref{prop-dual-imm-ferme}})\\
& \simeq & i_*\Rhom (\Omega _{\X'_{Y'_n}}, i^!\Omega _{\X_{Y_n}})\s (\text{\ref{dual}})\\
& \simeq & i_*\Rhom (\Omega _{\X'_{Y'_n}}, \Omega _{\X'_{Y'_n}})\s (i^!\Omega _{\X_{Y_n}} = \Omega _{\X'_{Y'_n}})\\
& \simeq & i_*\Lambda \s (\text{\ref{Dinvolutif}}).
\end{eqnarray*}
Therefore \ref{bigmess} is equal to
$$
j^*\Rhom (\Rhom (Rg_*i_*\Lambda , \Omega _{Y_n}), j_*\Omega
_{Y'_n})\simeq j^*\Rhom (\Rhom (j_*Rg'_*\Lambda , \Omega _{Y_n}),
j_*\Omega _{Y'_n}).
$$
Then
\begin{eqnarray*}
j^*\Rhom (\Rhom (j_*Rg'_*\Lambda , \Omega _{Y_n}), j_*\Omega _{Y'_n})& \simeq & j^*\Rhom (j_*\Rhom (Rg'_*\Lambda , j^!\Omega _{Y_n}), j_*\Omega _{Y'_n}) \s (\text{\ref{dual}})\\
& \simeq & j^*j_*\Rhom (j^*j_*\Rhom (Rg'_*\Lambda , \Omega _{Y'_n}), \Omega _{Y'_n}) \\
& \simeq & \Rhom (\Rhom (Rg'_*\Lambda , \Omega _{Y'_n}), \Omega _{Y'_n})\s (j^*j_* = \text{id})\\
& \simeq & Rg'_*\Lambda \s (\text{\ref{Dinvolutif}}).
\end{eqnarray*}

\hfill\hfill$\square$\end{proof} The functor $\text{id}\rightarrow
\mc F'$ induces for any $A\in \D_c(Y'_\bullet )$ and $B\in
\D_c(\X'_{Y'_\bullet })$ a morphism
$$
\Rhom (A, Rg'_*B)\rightarrow \Rhom (A, \mc F'(Rg'_*B))\simeq \Rhom
(A, \mc F(g^{\prime *}Rg'_*B))\rightarrow \Rhom (A, \mc F(B)),
$$
where the last morphism is induced by adjunction $g^{\prime
*}Rg'_*B\rightarrow B$. This map is functorial in $A$, so by
Yoneda's Lemma we get a canonical morphism $Rg'_*B\rightarrow \mc
F(B)$.    To prove \ref{5.4.8} we show that this map is an
isomorphism for all $B\in \D_c(\X'_{Y'_\bullet })$.

For this we can restrict the map to any $\X'_{Y'_n}$. Noting that
the shifts and Tate twists cancel as in \ref{bigmess}, we get
\begin{eqnarray*}
\mc F(B) |_{Y'_n} & \simeq & \Rhom (j^*\Rhom (Rg_*\Rhom (i_*\Rhom (B, \Omega _{\X'_{Y'_n}}), \Omega _{\X_{Y_n}}), \Omega _{Y_n}), \Omega _{Y'_n})\\
& \simeq & \Rhom (j^*\Rhom (Rg_*i_*\Rhom (\Rhom (B, \Omega _{\X'_{Y'_n}}), Ri^!\Omega _{\X_{Y_n}}),\Omega _{Y_n}),\Omega _{Y'_n})\s (\ref{dual})\\
& \simeq & \Rhom (j^*\Rhom (j_*Rg'_*B, \Omega _{Y_n}), \Omega _{Y'_n})\s (Ri^!\Omega _{\X_{Y_n}} = \Omega _{\X'_{Y'_n}}, \ref{Dinvolutif}, \text{and} j_*Rg'_* = Rg_*i_*)\\
& \simeq & \Rhom (j^*j_*\Rhom (Rg'_*B, Rj^!\Omega _{Y_n}), \Omega _{Y_n'})\s (\ref{dual})\\
& \simeq & \Rhom (\Rhom (Rg'_*B, \Omega _{Y'_n}), \Omega _{Y'_n})\s (j^*j_* = \text{id}, j^!\Omega _{Y_n} = \Omega _{Y'_n})\\
& \simeq & Rg'_*B \s (\ref{Dinvolutif}).
\end{eqnarray*}
We leave to the reader the task of verifying that this isomorphism
agrees with the map obtained by restriction from the morphism $\mc
F(B)\rightarrow Rg'_*B$ constructed above, thereby completing the
proof of \ref{5.4.8}. \hfill\hfill$\square$\end{proof}

Let $\pi :\X_{Y_\bullet }\rightarrow \X$ (resp. $\pi ':\X'_{Y'_\bullet
}\rightarrow \X'$) denote the projection. The isomorphism $\mc
F\simeq Rg'_*$ induces a morphism of functors
\begin{equation}\label{5.5.9.2}
\begin{matrix}
j^*D_{Y_\bullet }Rg_*D_{\X_{Y_\bullet }} & \rightarrow & j^*D_{Y_\bullet }Rg_*D_{\X_{Y_\bullet }}i_*i^* \s (\text{id}\rightarrow i_*i^*)\\
& \simeq & D_{Y'_\bullet }D_{Y'_\bullet }j^*D_{Y_\bullet }Rg_*D_{\X_{Y_\bullet }}i_*D_{\X'_{Y'_\bullet }}D_{\X'_{Y'_\bullet }}i^*\s (\ref{Dinvolutif})\\
& \simeq & D_{Y'_\bullet }\mc FD_{\X'_{Y'_\bullet }}i^*\s (\text{definition})\\
& \simeq & D_{Y'_\bullet }Rg'_*D_{\X'_{Y'_\bullet }}i^*\s (\ref{5.4.8}).
\end{matrix}
\end{equation}
This map induces a morphism
\begin{equation}\label{basearrow2}
\begin{matrix}
\rho ^*Rf_! & \simeq & \rho ^*Rq_*D_{Y_\bullet }Rg_*D_{\X_{Y_\bullet }}\pi ^*\s (\text{cohomological descent})\\
& \rightarrow & Rp_*j^*D_{Y_\bullet }Rg_*D_{\X_{Y_\bullet }}\pi ^*\s (\text{base change morphism})\\
&\rightarrow & Rp_*D_{Y'_\bullet }Rg'_*D_{\X'_{Y'_\bullet }}i^*\pi ^*\s (\ref{5.5.9.2})\\
& \simeq & Rp_*D_{Y'_\bullet }Rg'_*D_{\X'_{Y'_\bullet }}\pi ^{\prime *}a^*\s (i^*\pi ^* = \pi ^{\prime *}a^*)\\
& \simeq & Rf'_!a^*\s (\text{cohomological descent}).
\end{matrix}
\end{equation}
which we call the \emph{base change morphism}.

By construction this morphism is compatible with smooth base
change on $\mc  Y$ and $\mc  Y'$.  It follows that in order to
verify that \ref{basearrow2} is an isomorphism it suffices to
consider the case when $\mc  Y'$ and $\mc  Y$ are schemes.
Furthermore, by construction if $X_\bullet \rightarrow \mc  X$ is a
smooth hypercover and $X_\bullet '$ the base change to $\mc  Y'$,
then the base change arrow \ref{basearrow2} is compatible with the
spectral sequences \ref{5.16}.  It follows that to verify that
\ref{basearrow2} is an isomorphism it suffices to consider the
case of schemes which is \cite{SGA43}, XVII, 5.2.6. Finally the independence of the
choices follows by a standard argument from \ref{filteringremark}.
  This completes the proof of
\ref{5.5.6}. \hfill\hfill$\square$\end{proof}

\subsection{Equivalence of different definitions of base change morphism}

In this subsection we show that the base changed morphism defined in the previous subsection agrees with the morphism defined earlier for smooth morphisms, immersions, and universal homeomorphisms.

\subsubsection{The case when $\rho $ is smooth.}\label{5.6.1}

Choose a diagram as in \ref{5.5.6.3}, and
let $d$ denote the locally constant function on $\Y'$ which is the relative dimension of $\rho $. For any morphism $\mc Z\rightarrow \Y'$ we also write $d$ for the pullback of the function $d$ to $\mc Z$. Note  that
\begin{equation}\label{twisteq}
j^*\omega _{Y_\bullet }\simeq \omega _{Y'_\bullet }\langle -d\rangle, \ \ i^*\omega _{\X_{Y_\bullet }}\simeq \omega _{\X'_{Y_\bullet '}}\langle -d\rangle .
\end{equation}

\begin{lem} For any $A\in D_c(\X_{Y_\bullet })$ (resp. $B\in D_c(Y_\bullet )$) there is a natural isomorphism $D_{\X'_{Y'_\bullet }}(i^*A\langle d\rangle )\simeq i^*D_{\X_{Y_\bullet }}(A)$ (resp. $D_{Y'_\bullet }j^*(B\langle d\rangle )\simeq j^*D_{Y_\bullet }(B)$).
\end{lem}
\begin{proof} Consider the natural map
\begin{equation}\label{pullmorphism}
\begin{matrix}
i^*\Rhom (A, \omega _{\X_{Y_\bullet }})&\rightarrow & \Rhom (i^*A, i^*\omega _{\X_{Y_\bullet }})\\
&\simeq & \Rhom (i^*A, \omega _{\X'_{Y'_\bullet }})\langle -d \rangle\\
& \simeq & \Rhom (i^*A\langle d\rangle, \omega _{\X'_{Y'_\bullet }}).
\end{matrix}
\end{equation}
We claim that this map is an isomorphism.  This can be verified over each $\X'_{Y'_n}$.  Let $\pi _n:\X_{Y_n}\rightarrow \X$ (resp. $\pi _n':\X'_{Y'_n}\rightarrow \X'$) be the projection.  By the equivalence of triangulated categories $D_c(\X_{Y_\bullet })\simeq D_c(\X)$, there exists an object $A'\in D_c(\X)$ so that the restriction of $A$ to $\X_{Y_n}$ is isomorphic to $\pi _n^*A'$.  The morphism \ref{pullmorphism} is then identified with the isomorphism
\begin{equation*}
\begin{matrix}
i^*\Rhom (A, \omega _{\X_{Y_\bullet }})& \simeq & i^*\pi _n^*\Rhom (A', \Omega _{\X}) & { \ \ \ \ } (\ref{Rhom-im-inv})\\
& \simeq & \pi _n^{\prime *}\rho ^*\Rhom (A', \Omega _{\X}) & \\
& \simeq & \pi _n^{\prime *}\Rhom (\rho ^*A', \rho ^*\Omega _{\X})& { \ \ \ \ } (\ref{Rhom-im-inv})\\
& \simeq & \Rhom (\pi _n^{\prime *}\rho ^*A', \pi _n^*\rho ^*\Omega _{\X}) & { \ \ \ \ } (\ref{Rhom-im-inv})\\
& \simeq & \Rhom (i^*A\langle d\rangle, \omega _{\X'_{Y'_\bullet }}).
\end{matrix}
\end{equation*}
The same argument proves the statement $D_{Y'_\bullet
}j^*(B\langle d\rangle )\simeq j^*D_{Y_\bullet }(B)$.
\hfill\hfill$\square$\end{proof}

For any $A\in D_c(Y_\bullet )$, let $\alpha _A$ denote the isomorphism
\begin{eqnarray*}
j^*A\langle d\rangle & \simeq & j^*D_{Y_\bullet }D_{Y_\bullet }(A)\langle d\rangle \\
& \simeq & D_{Y'_\bullet }j^*D_{Y_\bullet }(A).
\end{eqnarray*}
For $B\in \D_c(\X_{Y_\bullet })$ let $\beta _B$ denote the isomorphism
\begin{eqnarray*}
i^*B\langle d \rangle & \simeq & i^*D_{\X_{Y_\bullet }}D_{\X_{Y_\bullet }}(B) \\
& \simeq & D_{\X'_{Y'_\bullet }}i^*D_{\X_{Y_\bullet }}(B).
\end{eqnarray*}
Also for $C\in \D_c(\X'_{Y'_\bullet })$ let $\gamma _C$ be the isomorphism
\begin{eqnarray*}
C\langle -d\rangle & \simeq & D_{\X'_{Y'_\bullet }}D_{\X'_{Y'_\bullet }}(C\langle -d\rangle )\\
& \simeq & D_{\X'_{Y'_\bullet }}i^*i_*D_{\X'_{Y'_\bullet }}(C\langle -d\rangle )\\
& \simeq  & i^*D_{\X_{Y_\bullet }}i_*D_{\X'_{Y'_\bullet }}(C),
\end{eqnarray*}
and let $\gamma ':D_{\X_{Y_\bullet }}i_*D_{\X'_{Y'_\bullet }}(C)\rightarrow i_*C\langle -d\rangle $ denote the map obtained by adjunction. This map also induces for every $E\in \D_c(\X_{Y_\bullet })$ a morphism  $\delta _E$ given by
\begin{equation*}
\begin{CD}
D_{\X_{Y_\bullet }}i_*i^*D_{\X_{Y_\bullet }}(E) @>>>  D_{\X_{Y_\bullet }}i_*D_{\X'_{Y'_\bullet }}(i^*E)\langle d\rangle @>\gamma '>> i_*i^*E.
\end{CD}
\end{equation*}

The map $\alpha _A$ is a special case of a more general class of morphisms. For  $A, M\in \D_c(Y_\bullet )$ let $S_{A, M}:j^*A\Otimes D_{Y'_\bullet }j^*D_{Y_\bullet }(M)\rightarrow D_{Y'_\bullet }j^*D_{Y_\bullet }(A\Otimes M)$ denote the composite
\begin{eqnarray*}
j^*A\Otimes\Rhom (j^*Rhom (M, \omega _{Y_\bullet }), \omega _{Y'_\bullet })& \simeq & j^*(A\Otimes\Rhom (\Rhom (M, \omega _{Y_\bullet }), j_*\omega _{Y'_\bullet })\\
& \rightarrow & j^*\Rhom (\Rhom (A, \Rhom (M, \omega _{Y_\bullet })), j_*\omega _{Y'_\bullet })\\
& \simeq & j^*\Rhom (\Rhom (A\Otimes M, \omega _{Y_\bullet }), j_*\omega _{Y'_\bullet })\\
& \simeq & D_{Y'_\bullet }j^*D_{Y_\bullet }(A\Otimes M).
\end{eqnarray*}
Here the second morphism is given by \ref{5.5.7} and the third morphism is by the adjunction property of $\otimes $.

\begin{lem} For any $A\in \D_c(Y_\bullet )$ the map $\alpha _A$ is equal to the composite
\begin{equation*}
\begin{CD}
j^*A\langle d\rangle \simeq  j^*A\Otimes j^*\Lambda \langle d\rangle @>\alpha _{\Lambda }>> j^*A\Otimes \Rhom (j^*Rhom (\Lambda , \omega _{Y_\bullet }), \omega _{Y'_\bullet })@>S_{A, \Lambda }>> D_{Y'_\bullet }j^*D_{Y_\bullet }(A).
\end{CD}
\end{equation*}
\end{lem}
\begin{proof} This follows from the definitions.
\hfill\hfill$\square$\end{proof}

\begin{lem} For any $A, B, M\in \D_c(Y_\bullet )$, the diagram
\begin{equation*}
\begin{CD}
j^*A\Otimes j^*B\Otimes D_{Y'_\bullet }j^*D_{Y_\bullet }(M)@>S_{B, M}>> j^*A\Otimes D_{Y'_\bullet }j^*D_{Y_\bullet }(B\Otimes M)\\
@V\simeq VV @VVS_{A, B\Otimes M}V \\
j^*(A\Otimes B)\Otimes D_{Y'_\bullet }j^*D_{Y_\bullet }(M)@>S_{A\Otimes B, M}>> D_{Y'_\bullet }j^*D_{Y_\bullet }(A\Otimes B\Otimes M)
\end{CD}
\end{equation*}
commutes.
\end{lem}
\begin{proof}
Consider the diagram
$$
\begin{CD}
A\Otimes B\Otimes j_*D_{Y'_\bullet }j^*D_{Y_\bullet }(M)@>>> A\Otimes[[B, [M, \omega _{Y_\bullet }]], j_*\omega _{Y'_\bullet }]@>>> A\Otimes [[B\Otimes M, \omega _{Y_\bullet }], j_*\omega _{Y'_\bullet }]\\
@| @VVV @VVV \\
A\Otimes B\Otimes j_*D_{Y'_\bullet }j^*D_{Y_\bullet }(M)@. [[A, [B, [M, \omega _{Y_\bullet }]]], j_*\omega _{Y'_\bullet }]@>>> [[A, [B\Otimes M, \omega _{Y_\bullet }]], j_*\omega _{Y_\bullet '}]\\
@| @VVV @VVV \\
A\Otimes B\Otimes j_*D_{Y'_\bullet }j^*D_{Y_\bullet }(M)@>>> [[A\Otimes B, [M, \omega _{Y_\bullet }]], j_*\omega _{Y'_\bullet }]@>>> [[A\Otimes B\Otimes M, \omega _{Y_\bullet }], j_*\omega _{Y'_\bullet }]
\end{CD}$$
where to ease the notation we write simply $[-,-]$ for $\Rhom (-,
-)$. An elementary verification shows that each of the small
inside diagrams commute, and hence the big outside rectangle also
commutes.  Applying $j^*$ we obtain the lemma.
\hfill\hfill$\square$\end{proof}

Similarly, for $A, M\in \D_c(\X'_{Y'_\bullet })$, let $R_{A, M}:i_*A\Otimes D_{\X_{Y_\bullet }}i_*D_{\X'_{Y'_\bullet }}(M)\rightarrow D_{\X_{Y_\bullet }}i_*D_{\X'_{Y'_\bullet }}(A\Otimes M)$ be the map
\begin{eqnarray*}
i_*A\Otimes D_{\X_{Y_\bullet }}i_*D_{\X'_{Y'_\bullet }}(M) & \simeq & i_*A\Otimes \Rhom (i_*\Rhom (M, \omega _{\X'_{Y'_\bullet }}), \omega _{\X_{Y_\bullet }})\\
& \simeq & i_*A\Otimes \Rhom (\Rhom (i_*M, i_*\omega _{\X'_{Y'_\bullet }}), \omega _{\X_{Y_\bullet }})\\
& \rightarrow & \Rhom (\Rhom (i_*A, \Rhom (i_*M, i_*\omega _{\X'_{Y'_\bullet }})), \omega _{\X_{Y_\bullet }})\\
& \simeq &\Rhom (\Rhom (i_*(A\Otimes M), i_*\omega _{\X'_{Y'_\bullet }}),\omega _{\X_{Y_\bullet }})\\
& \simeq &\Rhom (i_*\Rhom (A\Otimes M, \omega _{\X'_{Y'_\bullet }}), \omega _{\X_{Y_\bullet }})\\
& = &  D_{\X_{Y_\bullet }}i_*D_{\X'_{Y'_\bullet }}(A\Otimes M),
\end{eqnarray*}
where the third morphism is provided by \ref{5.5.7}. As above, one verifies that for $A, B, M\in \D_c(\X'_{Y'_\bullet })$ the diagram
\begin{equation*}
\begin{CD}
i_*A\Otimes i_*B\Otimes D_{\X_{Y_\bullet }}i_*D_{\X'_{Y'_\bullet }}(M)@>R_{B, M}>> i_*A\Otimes D_{\X_{Y_\bullet }}i_*D_{\X'_{Y'_\bullet }}(B\Otimes M)\\
@V\simeq VV @VVR_{A, B\Otimes M}V\\
i_*(A\Otimes B)\Otimes D_{\X_{Y_\bullet }}i_*D_{\X'_{Y'_\bullet }}(M)@>R_{A\Otimes B, M}>> D_{\X_{Y_\bullet }}i_*D_{\X'_{Y'_\bullet }}(A\Otimes B\Otimes M)
\end{CD}
\end{equation*}
commutes.

From this it follows that if $\varphi _A:A\rightarrow \mc F'(A)$ denotes the morphism constructed in the proof of \ref{5.4.8}, then the diagram
\begin{equation*}
\begin{CD}
A@= A\\
@V\varphi _AVV @VV\text{adjunction}V \\
D_{Y'_\bullet }j^*D_{Y_\bullet }g_*D_{\X_{Y_\bullet }}i_*D_{\X'_{\bullet }}g^{\prime *}A @. g^{\prime }_*g^{\prime *} A\\
@V\alpha _{g_*D_{\X_{Y_\bullet }}i_*D_{\X'_{\bullet }}g^{\prime *}}VV @AA\gamma A \\
j^*g_*D_{\X_{Y_\bullet }}i_*D_{\X_{Y'_\bullet }'}g^{\prime *}A @>\text{base change} >> g'_*i^*D_{\X_{Y_\bullet }}i_*D_{\X'_{\bullet }}g^{\prime *}A
\end{CD}
\end{equation*}
commutes.

Note that the base change morphism in the above diagram is an isomorphism.  This can be verified over each $\X'_{Y'_n}$.  Here the functor $D_{\X_{Y_n}}i_*D_{\X'_{Y'_n}}$ is up to shift and Tate torsion isomorphic to $i_!=i_*$.  The base change morphism is therefore induced by the isomorphism
\begin{equation*}
j^*g_*i_*\simeq j^*j_*g'_*\simeq g'_*i^*i_*.
\end{equation*}

 By the definition of the morphism in \ref{5.4.8} this implies that for any $A\in \D_c({\X'_{Y'_\bullet }})$ the diagram
\begin{equation*}
\begin{CD}
D_{Y'_\bullet }j^*D_{Y_\bullet }g_*D_{\X_{Y_\bullet }}i_*D_{\X'_{Y'_\bullet }}(A)@>\alpha _{g_*D_{\X_{Y_\bullet }}i_*D_{\X'_{Y'_\bullet }}(A)}>> j^*g_*D_{\X_{Y_\bullet }}i_*D_{\X'_{\bullet }}(A)\langle d\rangle \\
@V\ref{5.4.8}VV @VV\text{base change} V \\
g'_*(A)@<\gamma _A<< g'_*i^*D_{\X_{Y_\bullet }}i_*D_{\X'_{\bullet }}(A)\langle d\rangle
\end{CD}
\end{equation*}
commutes. Combining the commutativity of this diagram with the commutativity of the diagram (verification left to the reader)
\begin{equation*}
\begin{CD}
D_{\X_{Y_\bullet }}i_*D_{\X'_{Y'_\bullet }}D_{\X'_{Y'_\bullet }}i^*D_{\X_{Y_\bullet }}@>D_{\X'_{Y'_\bullet }}^2 = \text{id}>> D_{\X_{Y_\bullet }}i_*i^*D_{\X_{Y_\bullet }}\\
@V\gamma 'VV @VV\delta V \\
i_*D_{\X'_{Y'_\bullet }}i^*D_{\X_{Y_\bullet }}\langle -d\rangle @>\beta>> i_*i^*.
\end{CD}
\end{equation*}
one sees that the diagram
\begin{equation}\label{square}
\begin{CD}
D_{Y'_\bullet }j^*D_{Y_\bullet }g_*D_{\X_{Y_\bullet }}D_{\X_{Y_\bullet }}@>\ref{5.5.9.2}>> D_{Y'_\bullet }D_{Y'_\bullet }g'_*D_{\X_{Y_\bullet '}'}i^*D_{\X_{Y_\bullet }}\\
@V\alpha VV @VVD_{Y'_\bullet }^2 = \text{id} V \\
j^*g_*D_{\X_{Y_\bullet }}D_{\X_{Y_\bullet }}@. g'_*D_{\X'_{Y'_\bullet }}i^*D_{\X_{Y_\bullet }}\\
@VD_{\X_{Y_\bullet }}^2 = \text{id} VV @VV\beta V \\
j^*g_*\langle d\rangle @>\text{base change}>> g'_*i^*\langle d\rangle
\end{CD}
\end{equation}
commutes.

We are now ready to prove the equivalences of the two definitions of the base change morphism.  The morphism constructed in \ref{5.5} is the composite
\begin{equation*}
\begin{matrix}
\rho ^*f_! & = & \rho ^*q_*D_{Y_\bullet }g_*D_{\X_{Y_\bullet }}\pi ^*&\\
& \simeq & p_*j^*D_{Y_\bullet }g_*D_{\X_{Y_\bullet }}\pi ^*& \\
& \rightarrow  & p_*D_{Y'_\bullet }g'_*D_{\X'_{Y'_\bullet }}i^*\pi ^* & { \ \ \ \ }(\ref{5.5.9.2})\\
& \simeq & p_*D_{Y'_\bullet }g'_*D_{\X'_{Y'_\bullet }}\pi ^{\prime *}a^*&\\
& \simeq & f'_!a^*.
\end{matrix}
\end{equation*}
The dual version of this morphism is given by
\begin{equation*}
\begin{matrix}
\rho ^!f_*& = & D_{\Y'}\rho ^*q_*D_{Y_\bullet }g_*D_{\X_{Y_\bullet }}\pi ^*D_{\X}& \\
& \simeq & p_*D_{Y'_\bullet }j^*D_{Y_\bullet }g_*D_{\X_{Y_\bullet }}D_{\X_{Y_\bullet }}\pi ^*& \\
& \rightarrow & p_*D_{Y'_\bullet }D_{Y'_\bullet }g'_*D_{\X'_{Y'_\bullet }}i^*D_{\X_{Y_\bullet }}\pi ^* & { \ \ \ \ } (\ref{5.5.9.2}) \\
& \simeq & p_*g'_*i^*\pi ^*\langle d\rangle . & \\
& \simeq & f'_*a^!.&
\end{matrix}
\end{equation*}
By the commutativity of \ref{square} this is the same as the composite
\begin{equation*}
\begin{matrix}
\rho ^!f_* & = & D_{\Y'}\rho ^*q_*D_{Y_\bullet }g_*D_{\X_{Y_\bullet }}\pi ^*D_{\X}& \\
& \simeq &  p_*D_{Y'_\bullet }j^*D_{Y_\bullet }g_*D_{\X_{Y_\bullet }}D_{\X_{Y_\bullet }}\pi ^*& \\
& \simeq & p_*D_{Y'_\bullet }j^*D_{Y_\bullet }g_*\pi ^*& { \ \ \ \ } D_{\X_{Y_\bullet }}^2 = \text{id} \\
& \simeq & p_*j^*g_*\pi ^* \langle d\rangle & { \ \ \ \ } D_{Y'_\bullet }j^*D_{Y_\bullet }\simeq j^*\langle d\rangle  \\
& \rightarrow & p_*g'_*i^* \pi ^* \langle d \rangle & { \ \ \ \ } (\text{base change morphism }) \\
& \simeq &  f'_*a^!.&
\end{matrix}
\end{equation*}
From this it follows that the morphism defined in \ref{5.1} agrees with the one defined in \ref{5.5}.

\subsubsection{The case when $\rho $ is a universal homeomorphism}

The same argument used in the previous section shows the agreement of the base change morphism in \ref{5.5} with the base change morphism in \ref{5.4}.  Indeed the only property of smooth morphisms used in the previous section is that the dualizing sheaves can be described as in \ref{twisteq}.  This also holds when $\rho $ is a universal homeomorphism (with $d = 0$).

\subsubsection{The case when $\rho $ is an immersion}

With notation as in \ref{5.3}, note first that to prove that the two base change morphisms agree it suffices to show that they agree on sheaves of the form $\pi _*A$ with $A\in D_c(\X')$.  Indeed for any $B\in D_c(\X)$ either base change isomorphism factors as
\begin{equation*}
\begin{CD}
p^*f_!B@>\text{id}\rightarrow \pi _*\pi ^*>> p^*f_!\pi _*\pi ^*B@>>> \phi _!\pi ^*\pi _*\pi ^*B = \phi _!\pi ^*B.
\end{CD}
\end{equation*}

In order to prove that the two base change morphisms agree, it is useful to first give an alternate description of the morphism defined in \ref{5.3}.

With notation as in \ref{5.3},  there is for any $A\in D_c(\X')$ a canonical isomorphism
\begin{eqnarray*}
D_{\X'}(A) & \simeq & \Rhom (A, \pi ^!\Omega _{\X})\\
& \simeq & \pi ^*\pi _*\Rhom (A, \pi ^!\Omega _{\X})\\
& \simeq & \pi ^*\Rhom (\pi _*A, \Omega _{\X})\\
& \simeq & \pi ^*D_{\X}(\pi _*A),
\end{eqnarray*}
and similarly $D_{\Y'}\simeq p^*D_{\Y}p_*$.  We can also write these isomorphisms as $\pi _*D_{\X'} \simeq  D_\X\pi _*$ and $p_*D_{\Y'} \simeq D_{\Y}p_*$.

We therefore obtain a morphism
\begin{equation}\label{alternate2}
\begin{matrix}
p^*Rf_!(\pi _*A) &=& p^*D_{\Y}f_*D_{\X}(\pi _*A) & \\
& \simeq & p^*D_{\Y}f_*\pi _*D_{\X'}(A) & \\
& \simeq & p^*D_{\Y}p_*\phi _*D_{\X'}(A) & \\
& \simeq & p^*p_*D_{\Y'}\phi _*D_{\X'}(A)& \\
& \simeq & \phi _!A. &
\end{matrix}
\end{equation}

\begin{lem} This morphism agrees with the one defined in \ref{5.3}. In particular the morphism \ref{alternate2} is an isomorphism.
\end{lem}
\begin{proof} Chasing through the definitions this amounts to the commutativity of the following diagram
\begin{equation*}
\begin{CD}
p^*f_!(\pi _*A)@>\simeq >> p^*(p_*\Lambda \Otimes f_!\pi _*A)\\
@AAA @VV\text{projection formula}V \\
p^*D_{\Y}p_*p^*f_*D_{\X}(\pi _*A) @. p^*(f_!(f^*p_*\Lambda \Otimes \pi _*A))\\
@V\simeq VV @VV\simeq V \\
D_{\Y'}p^*f_*D_{\X}(\pi _*A) @. p^*f_!(\pi _*A\Otimes \pi _*\phi ^*\Lambda )\\
@A\text{base change}AA @VV\simeq V \\
D_{\Y'}\phi _*\pi ^*D_{\X}(\pi _*A) @. p^*f_!\pi _*(A\Otimes \phi ^*\Lambda )\\
@A\simeq AA @VV\simeq V \\
D_{\Y'}\phi _*D_{\X'}A@>\simeq >> \phi _!\pi ^*(A).
\end{CD}
\end{equation*}
We leave to the reader this verification.
\hfill\hfill$\square$\end{proof}

In particular, since the map \ref{alternate2} is an isomorphism we can define the base change morphism for $B\in D_c(\X)$ as the composite
\begin{equation}\label{alternate}
\begin{CD}
p^*f_!B @>>> p^*f_!(\pi _*\pi ^*B)@>\ref{alternate2}>> \phi _!\pi ^*B.
\end{CD}
\end{equation}

Using this alternate description of the base change morphism in \ref{5.3}, we can prove the equivalence with that given in \ref{5.5}.  By a standard reduction it suffices to consider the case of a closed immersion.  So fix the diagram \ref{5.5.6.1} with $\rho $  a closed immersion, and choose a diagram as in \ref{5.5.6.3}.  Since $\rho $ is a closed immersion we may without loss of generality assume that \ref{5.5.6.3} is cartesian.

\begin{lem}\label{pulllem}
The functors $j_*:D(Y_\bullet ')\rightarrow D(Y'_\bullet )$ and $i_*:D({\X'_{Y'_\bullet }})\rightarrow D(\X_{Y_\bullet })$ have right adjoints $j^!$ and $i^!$ respectively.
\end{lem}
\begin{proof} In fact $j^! = j^*R\Gamma _{Y'_\bullet }$ and $i^! = i^*R\Gamma _{\X'_{Y'_\bullet }}$.

\hfill\hfill$\square$\end{proof}

Note that for any $[n]\in \Delta $, the restriction of $j^!$ (resp. $i^!$) to a functor $D(Y_n)\rightarrow D(Y'_n)$ (resp. $D(\X_{Y_n})\rightarrow D(\X'_{Y'_n})$) agrees with the usual extraordinary inverse image.  This follows for example from the explicit description of these functors in the proof of \ref{pulllem}.

\begin{lem} There are canonical isomorphisms $\omega _{Y'_\bullet }\simeq j^!\omega _{Y_\bullet }$ and $\omega _{\X'_{Y'_\bullet}}\simeq i^!\omega _{\X_{Y_\bullet }}$.
\end{lem}
\begin{proof} By the glueing lemma \ref{1.2}, it suffices to construct an isomorphism over each $Y_n'$ (resp. $\X'_{Y'_n}$).
Let $d$ denote the relative dimension of $Y_n$ over $\Y$.  Then $d$ is also equal to the relative dimension of $Y'_n$ over $\Y'$, the relative dimension of $\X_{Y_n}$ over $\X$, and the relative dimension of $\X'_{Y_n'}$ over $\X'$.  We therefore have
\begin{equation*}
i^!\omega _{\X_{Y_\bullet }}|_{\X'_{Y'_n}} = i^!\Omega _{\X_{Y_n}}\langle -d\rangle \simeq \Omega _{\X'_{Y'_n}}\langle -d\rangle \simeq \omega _{\X'_{Y_\bullet }}|_{\X'_{Y'_n}}
\end{equation*}
and
\begin{equation*}
j^!\omega _{Y_\bullet }|_{Y_n.} = j^!\Omega _{Y_n}\langle -d\rangle \simeq \Omega _{Y'_n}\langle -d\rangle \simeq \omega _{Y'_\bullet }|_{Y_n'}.
\end{equation*}

\hfill\hfill$\square$\end{proof}

\begin{lem} For any $A\in D(Y'_\bullet )$ and $B\in D(Y_\bullet )$ (resp. $C\in D(\X'_{Y'_\bullet })$ and $E\in D(\X_{Y_\bullet })$) we have
\begin{equation*}
j_*\Rhom (A, j^!B)\simeq \Rhom (j_*A, B), \ \ i_*\Rhom (C, i^!E)\simeq \Rhom (i_*C, E).
\end{equation*}
\end{lem}
\begin{proof} Since $j^!$ is right adjoint to $j_*$, there is an adjunction morphism $j^!j_*\rightarrow \text{id}$.  This map induces a morphism
\begin{equation*}
j_*\Rhom (A, j^!B)\simeq \Rhom (j_*A, j_*j_!B)\rightarrow \Rhom (j_*A, B).
\end{equation*}
That this map is an isomorphism can be verified after restricting
to each $Y_n$ in which case it follows from the theory for schemes
\cite{SGA43}, XVIII, 3.1.10. The same argument gives the second
isomorphism in the Lemma. \hfill\hfill$\square$\end{proof}

\begin{cor} For any $A\in D_c(Y'_\bullet )$ there is a natural isomorphism $D_{Y_\bullet }j_*A\simeq j_*D_{Y'}(A)$, and for $B\in D_c(\X'_{Y'_\bullet })$ there is a canonical isomorphism $D_{\X_{Y_\bullet }}i_*B\simeq i_*D_{\X'_{Y'_\bullet }}B$.
\end{cor}

For $A\in \D_c(Y'_\bullet )$, let $\alpha _A$ denote the isomorphism
\begin{eqnarray*}
j_*A & \simeq & D_{Y'_\bullet }j^*j_*D_{Y'_\bullet }A\\
& \simeq & D_{Y'_\bullet }j^*D_{Y_\bullet }j_*A,
\end{eqnarray*}
and for $B\in D_{\X'_{Y'_\bullet }}$ let $\beta_B$ denote the isomorphism
\begin{eqnarray*}
i_*B & \simeq & i_*D_{\X'_{Y'_\bullet }}D_{\X'_{Y'_\bullet }}(B)\\
& \simeq & i_*D_{\X'_{Y'_\bullet }}i^*i_*D_{\X'_{Y'_\bullet }}(B)\\
& \simeq & i_*D_{\X'_{Y'_\bullet }}i^*D_{\X_{Y_\bullet }}(i_*B).
\end{eqnarray*}

Define $\gamma _B'$ to be the isomorphism
\begin{equation*}
\begin{matrix}
i_*B & \simeq & D_{\X_{Y_\bullet }}D_{\X_{Y_\bullet }}i_*B\\
& \simeq & D_{\X_{Y_\bullet }}i_*D_{\X'_{Y'_\bullet }}B,
\end{matrix}
\end{equation*}
and let $\gamma _B:i^*D_{\X_{Y_\bullet }}i_*D_{\X'_{Y'_\bullet }}(B)\rightarrow B$ be the isomorphism obtained by adjunction.

Following the same outline used in \ref{5.6.1} (replacing the $\alpha $'s, $\beta $'s, and $\gamma $'s by the above defined morphisms), one sees that the morphism \ref{5.5.9.2} in the case of a closed immersion is given by the composite
\begin{equation*}
\begin{matrix}
j^*D_{Y_\bullet }g_*D_{\X_{Y_\bullet }}i_* & \simeq & j^*D_{Y_\bullet }g_*i_*D_{\X'_{Y'_\bullet }}\\
& \simeq & j^*D_{Y_\bullet }j_*g'_*D_{\X'_{Y'_\bullet }}\\
& \simeq & j^*j_*D_{Y_\bullet '}g'_*D_{\X'_{Y'_\bullet }}\\
& \simeq & D_{Y_\bullet '}g'_*D_{\X'_{Y'_\bullet }}.
\end{matrix}
\end{equation*}
From this it follows that the sequence of morphisms in \ref{basearrow2} is identified via cohomological descent with the sequence of morphisms \ref{alternate}, and hence the two base change morphisms are the same.

\subsection{Kunneth formula}

Let $\mc Y_1$ and $\mc Y_2$ be  stacks, and set $\mc Y:= \mc
Y_1\times \mc Y_2$.  Let $p_i:\mc Y\rightarrow \mc Y_i$ ($i=1,2$)
be the projection and for two complexes $L_i\in \D^-_c(\mc Y_i)$
let $L_1{\Otimes_S}L_2\in \D(\mc Y)$ denote
$p_1^*L_1{\Otimes}_{\Lambda }p_2^*L_2$.

\begin{lemma}\label{dualdescription}
There is a natural isomorphism $K_{\mc Y}\simeq K_{\mc
Y_1}{\Otimes_S}K_{\mc Y_2}$.
\end{lemma}
\begin{proof}
By (\cite{SGA5}, III.1.7.6) there is for any smooth morphisms
$U_i\rightarrow \mc Y_i$ ($i=1,2$) with $U_i$ a scheme, a
canonical isomorphism
\begin{equation}
K_{\mc Y}|_{U_1\times _SU_2}\simeq K_{\mc
Y_1}|_{U_1}{\Otimes_S}K_{\mc Y_2}|_{U_2}.
\end{equation}
Furthermore, this isomorphism is functorial with respect to
morphisms $V_i\rightarrow U_i$.  It follows that the sheaf $K_{\mc
Y_1}{\Otimes_S}K_{\mc Y_2}$ also satisfies the $\mc
Ext$--condition (\ref{1.2}), and hence to give an isomorphism as in
the Lemma it suffices to give an isomorphism in the derived
category of $U_1\times _SU_2$ for all smooth morphisms
$U_i\rightarrow \mc Y_i$. \hfill\hfill$\square$\end{proof}

\begin{lemma}\label{hommap}
Let $(\T, \Lambda )$ be a ringed topos.  Then for any $P_1, P_2,
M_1, M_2\in \D(\T, \Lambda )$, there is a canonical morphism
$$
\Rhom (P_1, M_1){\Otimes}\Rhom (P_2, M_2)\rightarrow \Rhom
(P_1{\Otimes}P_2, M_1{\Otimes}M_2).
$$
\end{lemma}
\begin{proof}
It suffices to give a morphism
$$
\Rhom (P_1, M_1){\Otimes}\Rhom (P_2,
M_2){\Otimes}P_1{\Otimes}P_2\rightarrow M_1{\Otimes}M_2.
$$
This we get by tensoring the two evaluation morphisms
$$
\Rhom (P_i, M_i){\Otimes}P_i\rightarrow M_i.
$$
\hfill\hfill$\square$\end{proof}
 For the definition and standard properties of homotopy colimits
 we refer to \cite{Bok-Nee93}.

\begin{lemma}\label{reduc-tenseur} Let $A,B\in\D(\X)$.Then we have\begin{enumerate}
    \item $\hocolim\tau_{\leq n}A=A$;
    \item $A\Otimes B=\hocolim\tau_{\leq n}A\Otimes\tau_{\leq
    n}B$.
\end{enumerate}
\end{lemma}

\begin{proof}Consider the triangle
$$\oplus\tau_{\leq n}A\xrightarrow{1-\text{shift}}\oplus\tau_{\leq n}A\ra A\leqno(*)$$

If $C=\hocolim\tau_{\leq n}A$ is the cone of $1-\text{shift}$, one
gets a morphism $C\ra A$. By construction, one has $$\mc
H(C)=\varinjlim\mc H(\tau_{\leq n}A)=\mc H(A)$$ proving that $C\ra
A$ is an isomorphism. Tensoring (*) by $B$ we get therefore a
distinguished triangle
$$\oplus\tau_{\leq n}A\Otimes B\xrightarrow{1-\text{shift}}\oplus\tau_{\leq n}A
\Otimes B\ra A\Otimes B$$ proving $$\hocolim\tau_{\leq n}A\Otimes
B=A\Otimes B.$$ Applying this process again we find
$$\hocolim\tau_{\leq n}\Otimes\tau_{\leq m}B=A\Otimes B.$$ Because
the diagonal is cofinal in $\NN\times \NN$, the lemma follows.
\hfill\hfill$\square$\end{proof}
\begin{prop} For $L_i\in \D^-_c(\mc Y_i)$ ($i=1,2$), there is a canonical isomorphism
\begin{equation}
D_{\mc Y_1}(L_1){\Otimes}_SD_{\mc Y_2}(L_2)\simeq D_{\mc
Y}(L_1{\Otimes_S}L_2).
\end{equation}
\end{prop}
\begin{proof}
By \ref{dualdescription} and \ref{hommap} there is a canonical
morphism (note here we also use that $K_{\mc Y_i}$ has finite
injective dimension)
\begin{equation}
D_{\mc Y_1}(L_1){\Otimes}_SD_{\mc Y_2}(L_2)\rightarrow D_{\mc
Y}(L_1{\Otimes_S}L_2).
\end{equation}
To verify that this map is an isomorphism, it suffices to show
that for every $j\in Z$ the map
\begin{equation}
\mc H^j(D_{\mc Y_1}(L_1){\Otimes}_SD_{\mc Y_2}(L_2))\simeq \mc
H^j(D_{\mc Y}(L_1{\Otimes_S}L_2)).
\end{equation}
Because $\Otimes$ commutes with homotopy colimits
(\ref{reduc-tenseur}), we deduce from $D(A)=\hocolim D(\tau_\geq m
A)$  (use~\ref{reduc-tenseur}) that to prove this we may replace
$L_i$ by $\tau _{\geq m}L_i$ for $m$ sufficiently negative, and
therefore it suffices to consider the case when $L_i\in \D^b_c(\mc
Y_i)$. Furthermore, we may work locally in the smooth topology on
$\mc Y_1$ and $\mc Y_2$, and therefore it suffices to consider the
case when the stacks $\mc Y_i$ are schemes.  In this case the
result is \cite{SGA43}, XVII, 5.4.3.
\hfill\hfill$\square$\end{proof}

Now consider morphisms of  $S$--stacks $f_i:\mc X_i\rightarrow
\mc Y_i$ ($i=1,2$), and let $f:\mc X:= \mc X_1\times _S\mc
X_2\rightarrow \mc Y:= \mc Y_1\times _S\mc Y_2$ be the morphism
obtained by taking fiber products.  Let $L_i\in \D^-_c(\mc X)$.

\begin{theorem} There is a canonical isomorphism in $\D_c(\mc Y)$
\begin{equation}\label{Kunnethmap}
Rf_!(L_1{\Otimes_S}L_2)\rightarrow
Rf_{1!}(L_1){\Otimes_S}Rf_{2!}(L_2).
\end{equation}
\end{theorem}
\begin{proof}
We define the morphism \ref{Kunnethmap} as the composite
\begin{equation*}
\begin{CD}
Rf_!(L_1{\Otimes_S}L_2)@>\simeq >> D_{\mc Y}(f_*D_{\mc X}(L_1{\Otimes_S}L_2))\\
@>\simeq >> D_{\mc Y}(f_*(D_{\mc X_1}(L_1){\Otimes_S}D_{\mc X_2}(L_2)))\\
@>>>  D_{\mc Y}(f_{1*}D_{\mc X_1}(L_1){\Otimes_S}(f_{2*}D_{\mc X_2}(L_2)))\\
@>\simeq >> D_{\mc Y_1}(f_{1*}D_{\mc X_1}(L_1)){\Otimes_S}D_{\mc Y_2}(f_{2*}D_{\mc X_2}(L_2))\\
@>\simeq >> Rf_{1!}(L_1){\Otimes_S}Rf_{2!}(L_2).
\end{CD}
\end{equation*}
That this map is an  isomorphism follows from a standard reduction
to the case of schemes using hypercovers of $\mc X_i$, biduality,
and the spectral sequences \ref{5.16}.
\hfill\hfill$\square$\end{proof}

\end{document}